\documentclass[12pt, leqno]{amsart}

\usepackage{graphicx}
\usepackage{amssymb,amsmath,amsthm,amscd}
\usepackage{mathrsfs}
\usepackage{enumerate}
\usepackage{enumitem}
\usepackage{verbatim}
\usepackage[usenames,dvipsnames]{color}
\usepackage[colorlinks=true, pdfstartview=FitV,
linkcolor=blue,citecolor=blue,urlcolor=blue]{hyperref}
\usepackage[all]{xy}
\usepackage{tikz}
\usetikzlibrary{shapes, arrows}
\usepackage{tikz-cd}
\usetikzlibrary{matrix, arrows}
\usepackage{stmaryrd} 
\usepackage{dsfont}
\usepackage{bm} %bold font for greek letters
\usepackage{scalerel}[2016/12/29] 
\usepackage{amsthm}
\usepackage{adjustbox}
\usepackage{float}
\usepackage{array}
\usepackage{multirow} 
\usepackage{ytableau} 
\usepackage{mathdots}

\allowdisplaybreaks[4]

\DeclareMathSymbol{\shortminus}{\mathbin}{AMSa}{"39}

\newcommand{\nc}{\newcommand}

\numberwithin{equation}{section}

\usepackage[colorinlistoftodos]{todonotes}
\usepackage{mathtools}
\mathtoolsset{showonlyrefs,showmanualtags}
\usepackage{empheq}

\makeatletter
\newcommand{\reqnomode}{\tagsleft@false\let\veqno\@@eqno}

\theoremstyle{plain}
\newtheorem{lem}{Lemma}[section]
\newtheorem{pro}[lem]{Proposition}
\newtheorem{thm}[lem]{Theorem}
\newtheorem{cor}[lem]{Corollary}
\newtheorem{defi}[lem]{Definition}

\newcommand{\Pro}{\begin{pro}}
	\newcommand{\enpro}{\end{pro}}
\newcommand{\Lem}{\begin{lem}}
	\newcommand{\enlem}{\end{lem}}
\newcommand{\Thm}{\begin{thm}}
	\newcommand{\enthm}{\end{thm}}
\newcommand{\Cor}{\begin{cor}}
	\newcommand{\encor}{\end{cor}}
\newcommand{\Defi}{\begin{defi}}
	\newcommand{\enDefi}{\end{defi}}
\newcommand{\Proof}{\begin{proof}}
	\newcommand{\enproof}{\end{proof}}

\theoremstyle{definition} %?
\newtheorem{rem}[lem]{Remark}
\newtheorem{exam}[lem]{Example}
\newtheorem{Convention}[lem]{Convention}

\newcommand{\Conv}{\begin{Convention}}
	\newcommand{\enconv}{\end{Convention}}
\nc{\Rem}{\begin{rem}}
	\nc{\enrem}{\end{rem}}

\newcommand{\arxiv}[1]{\href{http://arxiv.org/abs/#1}{\tt arXiv:\nolinkurl{#1}}}

\newcommand{\monoto}{\hookrightarrow}
\nc{\epito}{\twoheadrightarrow}
\newcommand{\isoto}[1][]{\mathop{\xrightarrow%
		[{\raisebox{.3ex}[0ex][.3ex]{$\scriptstyle{#1}$}}]%
		{{\raisebox{-.6ex}[0ex][-.6ex]{$\mspace{2mu}\sim\mspace{2mu}$}}}}}

\nc{\rmkend}{\hfill$\triangledown$}
\nc{\defend}{\hfill$\triangle$}

\nc{\ccc}{\mathfrak{c}}
\nc{\CCC}{\mathfrak{C}}
\nc{\Ck}{\mathfrak{C}}

\nc{\kor}{\mathbb{C}}
\nc{\indx}{\mathbb{I}}

\nc{\CC}{C}
\nc{\cc}{c}
\nc{\sss}{s}
\nc{\ck}{\mathfrak{c}}

\nc{\Bg}{B}
\nc{\Ag}{A}
\nc{\As}{\pmb{\Ag}(z)}
\nc{\Asmone}{\pmb{\Ag}^{(-1)}(z)}
\nc{\Aps}{\pmb{\Ag}_+(z)}
\nc{\Apsone}{\pmb{\Ag}^{(1)}_+(z)}
\nc{\Apsmone}{\pmb{\Ag}^{(-1)}_+(z)}
\nc{\Ams}{\pmb{\Ag}_-(z)}
\nc{\Hg}{H}
\nc{\Thg}{\Theta}
\nc{\Thgs}{\pmb{\Theta}}
\nc{\Thgsr}{\pmb{\grave{\Theta}}}
\nc{\bTh}{\pmb{\grave{\Theta}}(z)}
\nc{\avar}[1]{\pmb{\Ag}_{#1, +}(z)} 
\nc{\ivar}[2]{\pmb{\Ag}_{#1, +}^{#2}(z)} 
\nc{\thvar}[1]{\pmb{\grave{\Theta}}_{#1}(z)}

\nc{\smin}{{\shortminus}}

\nc{\Ui}{\widetilde{\mathbf{U}}^\imath}

\nc{\fext}[2]{{#1}[\negthinspace[#2]\negthinspace]}

\nc{\gTh}{\grave{\Theta}}

\nc{\DD}{\pmb{D}}
\nc{\KK}{\mathbb{K}}

\nc{\xgp}{x^+}
\nc{\Xgp}{\pmb{x}^+}
\nc{\xgm}{x^-}
\nc{\Xgm}{\pmb{x}^-}
\nc{\xgpm}{x^{\pm}}
\nc{\psig}{\phi^+}
\nc{\phig}{\phi^-}
\nc{\phipm}{\phi^\pm}
\nc{\Psig}{\pmb{\phi^+}}
\nc{\Phig}{\pmb{\phi^-}}
\nc{\Phipm}{\pmb{\phi^\pm}}
\nc{\hg}{h}

\nc{\bY}{\mathbf{Y}}
\nc{\bA}{\mathbf{A}}
\nc{\bT}{\mathbf{T}}

\nc{\Eg}{e}
\nc{\Kg}{K}
\nc{\ev}{{\rm ev}}

\nc{\Rep}{\on{Rep}} 

\DeclareRobustCommand{\sqbin}{\genfrac{[}{]}{0pt}{}}

\nc{\qq}{(q-q^{-1})^{-1}}
\nc{\factor}{\Omega}

\nc{\chring}{\Z[Y_a^{\pm 1}]_{a \in \C^\times}}
\nc{\chmap}{\chi_q}
\nc{\chmod}{\Z[\mathbf{Y}_a^{\pm 1}]_{a \in \C^\times}}
\nc{\ourchmap}{\boldsymbol{\chi}_q}
\nc{\Yring}{\mathcal{Y}} 
\nc{\Ymod}{\pmb{\mathcal{Y}}}
\nc{\Xmod}{\pmb{\mathcal{X}}}
\nc{\Yg}{\mathbf{Y}_{i,a}}
\nc{\Ys}{\mathbf{Y}_{i,\mathbf{s}}}
\nc{\Kmat}{\mathcal{K}^0} 

\nc{\Tr}{\on{Tr}}
\nc{\id}{\on{id}}

\nc{\ourR}{\mathbf{R}}
\nc{\ourQ}{\mathbf{Q}}
\nc{\ourY}{\mathbf{Y}}
\nc{\ourX}{\mathbf{X}}
\nc{\ourA}{\mathbf{A}}
\nc{\ourB}{\mathbf{B}}

\nc{\Oq}{\widetilde{\mathcal{O}}_q}
\nc{\Oqi}{\widetilde{\mathcal{O}}_q^{[i]}} 
\nc{\Oqii}{\widetilde{\mathcal{O}}_q^{[i,i+1]}} 
\nc{\Oqiii}{\widetilde{\mathcal{O}}_q^{[i,i+1,i+2]}} 
\nc{\oq}{\mathcal{O}_q}
\nc{\Oqc}{\mathcal{O}_q^{\mathbf{c}}(\widehat{\g})}
\nc{\car}{\mathcal{H}}
\nc{\qaa}{U_q(\widehat{\mathfrak{sl}}_{n+1})}
\nc{\qla}{U_q(L\mathfrak{sl}_{n+1})}
\nc{\uqla}{\widetilde{U}_q(L\mathfrak{g})}
\nc{\drqaa}{{}^{\mathrm{Dr}}\qaa}
\nc{\uqsl}{U_q(\widehat{\mathfrak{sl}}_2)}
\nc{\uLsl}{U_qL\mathfrak{sl}_2}
\nc{\Serre}{\mathsf{Serre}}
\nc{\Sym}{\on{Sym}}
\nc{\UXp}{UX_+} 

\nc{\Tbr}{\mathbf{T}}
\nc{\degdr}{\on{deg}^{\mathrm{Dr}}}
\nc{\Uq}{\mathbf{U}}
\nc{\Uu}{\widetilde{\mathbf{U}}}
\nc{\gaf}{\widehat{\mathfrak{g}}}

\nc{\ada}{\ad_{\bhg_{-1}}}
\nc{\adb}{\ad_{\bhg_{1}}}
\nc{\adc}{\ad_{\bHg_{1}^{(2)}}}

\nc{\tAg}{\widetilde{A}}
\nc{\hAg}{\widehat{A}}

\nc{\bHg}{\overline{\Hg}}
\nc{\bhg}{\overline{\hg}}
\nc{\tHg}{\widetilde{H}}

\nc {\DrOq}{{}^{\mathrm{Dr}}\Oq}
\nc{\gup}[1]{^{(#1)}}
\nc{\gupp}{^{(1),+}}
\nc{\mi}{^{-1}}
\nc{\adh}{\operatorname{ad}_{\bar{\hg}_{-1}}}
\nc{\adhp}{\operatorname{ad}_{\bar{\hg}_{1}}}
\nc{\Omg}{\Omega^{-1}}
\nc{\Htwo}{\overline{H}_1^{(2)}}
\nc{\ad}{\operatorname{ad}}

\nc{\sspan}{\on{span}}

\newcommand{\commentout}[1]{}

\newcommand{\on}{\operatorname}
\nc{\be}{\begin{enumerate}}
	\nc{\ee}{\end{enumerate}}
\newcommand{\eq}{\begin{equation}}
	\newcommand{\eneq}{\end{equation}}
\nc{\bc}{\begin{cases}}
	\nc{\ec}{\end{cases}}
\newcommand{\eqn}{\begin{eqnarray*}}
	\newcommand{\eneqn}{\end{eqnarray*}}
\newcommand{\ba}{\begin{array}}
	\newcommand{\ea}{\end{array}}

\newcommand{\C}{{\mathbb C}}

\newcommand{\Z}{{\mathbb Z}}
 
\newcommand{\R}{{\mathbb R}} 

\newcommand{\g}{{\mathfrak{g}}}
\nc{\Ad}{\operatorname{Ad}}

\nc{\gr}{\on{gr}}

\newcommand{\End}{\operatorname{End}}
\nc{\Aut}{\operatorname{Aut}}

\nc{\coker}{\operatorname{coker}}

\nc{\Img}{\on{Im}}
\nc{\res}{\on{res}}

\nc{\modv}[1]{{#1}\operatorname{-mod}}

\nc{\bl}{\bigl(}
\nc{\br}{\bigr)}

\newlength{\mylength}
\DeclareRobustCommand{\SkipTocEntry}[5]{}

\AtBeginDocument{%
   \def\MR#1{}
}

\title[Evaluation modules for quantum symmetric pairs]{Boundary $q$-characters of evaluation modules for split quantum affine symmetric pairs}
\author[J.-R. Li]{Jian-Rong Li}
\address{Faculty of Mathematics, University of Vienna, Oskar Morgenstern Platz 1, 1090 Vienna, Austria, OrciD: 0000-0001-7896-7391}
\email{\href{mailto:jianrong.li@univie.ac.at}{jianrong.li@univie.ac.at}}

\author[T. Prze\'{z}dziecki]{Tomasz Prze\'{z}dziecki}
\address{Faculty of Mathematics, University of Vienna, Oskar Morgenstern Platz 1, 1090 Vienna, Austria, OrciD: 0000-0001-9700-1007}
\email{\href{mailto:tomasz.przezdziecki@univie.ac.at}{tomasz.przezdziecki@univie.ac.at}}



\thanks{The first author was supported by the Austrian Science Fund (FWF): P-34602, Grant DOI: 10.55776/P34602, and PAT 9039323, Grant-DOI 10.55776/PAT9039323. The second author was supported by the EPSRC grant No.\ EP/W022834/1 \emph{Kac--Moody quantum symmetric pairs, KLR algebras and generalized Schur--Weyl duality}.}

\begin{document}

\begin{abstract}
We study evaluation modules for quantum symmetric pair coideal subalgebras of affine type $\mathsf{AI}$. By computing the action of the generators in Lu and Wang's Drinfeld-type presentation on Gelfand--Tsetlin bases, we determine the spectrum of a large commutative subalgebra arising from the Lu--Wang presentation. This leads to an explicit formula for boundary analogues of $q$-characters in the setting of quantum affine symmetric pairs. We interpret this formula combinatorially in terms of semistandard Young tableaux. Our results imply that boundary $q$-characters share familiar features with ordinary $q$-characters -- such as a version of the highest weight property -- yet they also display new phenomena, including an extra symmetry. In particular, we provide the first examples of boundary $q$-characters for quantum affine symmetric pairs that do not arise from restriction of ordinary $q$-characters, thereby revealing genuinely new structures in this new setting. 
\end{abstract}

\maketitle

\setcounter{tocdepth}{1}
\tableofcontents

\section{Introduction}

\subsection{Quantum symmetric pairs of type $\mathsf{AI}$}

\nc{\aaui}{\mathbf{U}^\imath_{\mathsf{aff}}}
\nc{\ichmap}{\chi_q^\imath}
\nc{\Uhz}{\fext{U_q(\widetilde{\mathfrak{h}})}{z}} 
\nc{\auiv}[1]{\Ui_{\mathsf{aff},#1}}
\nc{\auivv}[1]{\mathbf{U}^\imath_{\mathsf{aff},#1}}
\nc{\eva}{\on{ev}_a}
\nc{\aui}{\Ui_{\mathsf{aff}}}

A non-standard $q$-deformation $U'_q(\mathfrak{so}_n)$\footnote{In the text we will typically abbreviate it as $\mathbf{U}^\imath$.} of the universal enveloping algebra of $\mathfrak{so}_n$ was first constructed by Gavrilik and Klimyk \cite{GavKl} in 1991. Unlike the usual Drinfeld--Jimbo quantum group $U_q(\mathfrak{so}_n)$, it does not possess a Hopf algebra structure, but can instead be realized as a coideal subalgebra inside $U_q(\mathfrak{sl}_n)$. Following the work of, \emph{i.a.},  Letzter \cite{letzter-99, letzter-02} and Kolb \cite{kolb-14}, this non-standard $q$-deformation is now recognized as an example of a more general pattern, involving simultaneous \emph{compatible} quantizations of a semisimple Lie algebra~$\g$, and a subalgebra $\g^\theta$ fixed under an involution. From this point of view, $U'_q(\mathfrak{so}_n)$ is a quantum symmetric pair (QSP) coideal subalgebra of type $\mathsf{AI}$. 

A complete classification of simple finite-dimensional $U'_q(\mathfrak{so}_n)$-modules was first established by Iorgov and Klimyk \cite{Ior-Kl} (when $q$ is not a root of unity), using Gelfand--Tsetlin methods. This result was later reproved by Wenzl \cite{Wenzl}, who proposed a Verma-type module approach. Parts of this classification were also later recovered by Watanabe \cite{Wata}, using a more general `classical weight module' approach. 

Interestingly, the representation theory of $U'_q(\mathfrak{so}_n)$ is rather more complicated than that of the usual quantum group $U_q(\mathfrak{so}_n)$. In particular, the absence of a triangular decomposition poses an obstacle to adapting classical methods from the theory of semisimple Lie algebras. This complication is also reflected by the extra richness  of the representation theory of $U'_q(\mathfrak{so}_n)$, namely the fact that it includes two types of representations: classical and \emph{non-classical} ones. The former are $q$-deformations of the corresponding irreducible representations of the Lie algebra $\mathfrak{so}_n$, while the latter do not have a well-defined classical limit.

\subsection{Evaluation modules} 

In this paper, we study an affine analogue of $U'_q(\mathfrak{so}_n)$, i.e, the QSP coideal subalgebra $\aaui$ of \emph{affine} type $\mathsf{AI}$. In contrast to the finite case, relatively little is known about its representation theory. A classification of finite-dimensional simple modules exists only in rank one, i.e., in the case of the $q$-Onsager algebra \cite{ito-ter-10}. 
Incidentally, more information is available for some other families of coideal subalgebras. 
In particular, classifications of finite-dimensional irreducibles have already been obtained for QSP coideal subalgebras in affine type $\mathsf{AII}$ \cite{gow-molev}, and for orthogonal ($\mathsf{AI}$) and symplectic  ($\mathsf{AII}$) twisted Yangians \cite{Molev-Yan}, using the RTT presentation. 

The representation theories of QSP coideal subalgebras of \emph{finite} and \emph{affine} $\mathsf{AI}$ types are connected by an evaluation homomorphism. Indeed, one of the original motivations for the introduction of twisted Yangians and affine QSP coideal subalgebras was to generalize the type $\mathsf{A}$ evaluation homomorphisms 
\[
Y(\mathfrak{sl}_n) \longrightarrow U(\mathfrak{sl}_n), \qquad U_q(L\mathfrak{sl}_n) \longrightarrow U_q(\mathfrak{sl}_n) 
\] 
to other classical types. 
Such analogues 
\[
Y^{tw}(\mathfrak{so}_n) \longrightarrow U(\mathfrak{so}_n), \qquad \aaui \longrightarrow U'_q(\mathfrak{so}_n)
\] 
were constructed in \cite{MRS-03}, using the RTT presentation (c.f.\ \cite{Ols, MNO}). Our first result (Theorem \ref{pro: ev homo}) is an alternative construction of an evaluation homomorphism, based purely on the braid group action from \cite{lu-wang-21}. 

Evaluation modules are known to play a key role in the representation theory of the quantum loop algebra  $U_q(L\mathfrak{sl}_n)$. They are precisely the minimal affinizations and include Kirillov--Reshetikhin modules. 
By \cite{chari-pressley-94}, every finite-dimensional $U_q(L\mathfrak{sl}_n)$-module can be obtained as a subquotient of a tensor product of evaluation modules. 
Also, by \cite{CP-small}, evaluation modules are precisely the small modules of $U_q(L\mathfrak{sl}_n)$, i.e., those which remain irreducible as representations of $U_q(\mathfrak{sl}_n)$. 
Therefore, it is natural to expect that evaluation modules will also play an important role in the representation theory of the  QSP coideal subalgebra $\aaui$.

\subsection{Boundary $q$-characters} 

The $q$-characters of quantum affine algebras were introduced by Frenkel and Reshetikhin \cite{FrenRes} (and by Knight \cite{Knight} in the Yangian case). They proposed two equivalent definitions: via the universal $R$-matrix, and via the spectrum of Drinfeld--Cartan operators. A third definition (along with a $t$-deformation and an axiomatic characterization), via quiver varieties, was later given by Nakajima \cite{Naka-qt}. Since then, $q$-characters have found numerous applications, not only in representation theory, but also in other fields, such as integrable systems and cluster theory, see, \emph{e.g.}, \cite{Hernandez-kr, Her-Lec-10, Frenkel-Hernandez, HerLec-KR}. 

Recent advances in understanding the structure of quantum affine symmetric pairs, in particular on the algebraic side, have provided the tools to generalize the theory of $q$-characters to this new setting. A crucial development was the construction of a Drinfeld-type presentation for split and certain quasi-split affine quantum symmetric pair coideal subalgebras by Lu and Wang \cite{lu-wang-21}, building on the earlier work of Baseilhac and Kolb \cite{bas-kol-20} in rank one (see also \cite{ZhangDr, LWZ-quasi}). In analogy to the usual Drinfeld presentation for quantum affine algebras, the Lu--Wang presentation also exhibits a large commutative subalgebra, 
generated by certain distinguished elements~$\Theta_{i,m}$. 
In \cite{Przez-23, LP24, LP26}, we proposed to define $q$-characters for quantum symmetric pairs (which we call `\emph{boundary $q$-characters}') using the generalized eigenvalues of the operators $\Theta_{i,m}$. The \emph{factorization} and \emph{coproduct formulae}: 
\begin{align*}
\Thgsr_i(z) &\equiv \pmb{\phi}_i^-(z\mi)\pmb{\phi}_i^+(C z) \quad \mod \ \fext{U_q(L\g)_+}{z}, \\ 
\Delta(\Thgsr_i(z)) &\equiv \Thgsr_i(z) \otimes \Thgsr_i(z) \quad  \ \  \mod \ \fext{\aaui  \otimes U_q(L\g)_{+}}{z}, 
\end{align*}
established in \emph{op.\ cit.}, imply that the boundary $q$-character map $\ichmap$, defined in this manner, is compatible with the usual $q$-character homomorphism. More precisely, the diagram 
\[
\begin{tikzcd}[ row sep = 0.2cm]
{[}\Rep U_q(L\g)] \arrow[r, "\chmap"] & \Z[Y_{i,a}^{\pm 1}]_{i \in \indx_0, a \in \C^{\times}}  \\
 \curvearrowright  & \curvearrowright  \\
{[}\Rep \aaui] \arrow[r, "\ichmap"] & \Uhz 
\end{tikzcd}
\]
commutes if $\Uhz$ is endowed with an appropriate `twisting' action of $\Z[Y_{i,a}^{\pm 1}]_{i \in \indx_0, a \in \C^{\times}}$. 
This is an analogue of the fact that $\chmap$ is a ring homomorphism, or Axiom 3 in Nakajima's abstract characterization \cite{Naka-qt}. 

So far very little is known about the properties and behaviour of boundary $q$-characters. In the case of restriction representations (from $U_q(L\mathfrak{sl}_n)$), they can be computed directly using the coproduct and factorization formulae above. However, since most simple $\aaui$-modules do not arise via restriction, the usefulness of this approach is limited. 

In the present paper, we take advantage of the fact that the finite-dimensional representation theory of $U'_q(\mathfrak{so}_n)$ is well understood, and study the boundary $q$-characters of $\aaui$-modules which arise via evaluation. Using Gelfand--Tsetlin theory, we are able to produce explicit character formulae, interpret them in terms of combinatorics of semistandard Young tableaux, and deduce analogues of several key properties enjoyed by $q$-characters in the quantum affine algebra case. 
We expect these results to be an important step towards a systematic theory of boundary $q$-characters. 

Nevertheless, it should be emphasized that the quantum symmetric pair setting is more complex than its quantum group analogue, with several major challenges to further progress. These include, for example, the absence of a triangular decomposition, no classification of finite-dimensional simple modules, and the lack of an alternative description of boundary $q$-characters via the universal $K$-matrix. Some progress towards the latter has recently been made in \cite{AppelVlaar2, AppelVlaar, AppelVlaar3, AV4}, but it remains unclear how the constructions in \emph{op.\ cit.} relate to our definition of boundary $q$-characters.

\subsection{Structure and main results}  

In Section \ref{sec: QSP}, we recall the definitions of QSP coideal subalgebras of split type, the associated braid group action and the Lu--Wang presentation. In Section \ref{sec: boundary q}, we recall the construction of $q$-characters from  \cite{FrenRes, FrenMuk-comb}, as well as the proposed definition of boundary $q$-characters from \cite{Przez-23, LP24, LP26}. Section \ref{subsec: rep th ui} collects the main results from the finite-dimensional representation theory of $U'_q(\mathfrak{so}_n)$. 

Let us now summarize our main results. \emph{Firstly}, in Theorem \ref{pro: ev homo}, we construct an evaluation homomorphism
\[
\eva \colon \aaui \to \mathbf{U}^\imath
\]
in terms of the Kolb--Letzter presentation, using the braid group action. 
\emph{Secondly}, in Section \ref{sec: GT action}, we compute the action of the Lu--Wang generators on the Gelfand--Tsetlin bases of evaluation modules (Theorems \ref{thm: A-action on GT}--\ref{thm: Theta-action on GT} and \ref{thm: A-action on GT nc}--\ref{thm: Theta-action on GT nc}). In particular, we show that the operators $\Theta_{i,m}$ act semisimply and compute their spectra. 
This is the technical core of the paper. Our proof relies on the Lu--Wang presentation, the evaluation homomorphism, and the Gelfand--Tsetlin-type formulae due Gavrilik, Iorgov and Klimyk \cite{GavKl, GI97, IK00}. In particular, we do not use the RTT presentation at all, which distinguishes our proof from the original proofs in the quantum affine algebra and (non-twisted) Yangian case \cite{FrMukHopf, NazTar}.  

\emph{Thirdly}, in Section \ref{sec: app to bq char}, we deduce an explicit formula for the boundary $q$-character of an evaluation module (Corollaries \ref{cor: iq char formula classical}--\ref{cor: iq char formula non-classical}), and provide a combinatorial interpretation in terms of semistandard Young tableaux (Corollary \ref{cor: tableaux formula for ich}). We also prove an analogue (Corollary \ref{cor: root property}) of the `highest weight property' (or Nakajima's Axiom 1) satisfied by ordinary $q$-characters, namely the fact that 
\eq \label{eq: into hwp}
\textstyle \chi_q(V) = M_+(1 + \sum_{p}M_p),
\eneq 
where $M_+$ is a dominant monomial, and each $M_p$ is a product of factors $A_{i,c}\mi$, which are $q$-analogues of negative simple roots. 

There are two important caveats regarding Corollary \ref{cor: root property}. Firstly, in the QSP setting, it is customary to use the concept of weight proposed by Letzter \cite{Letzter-cartan}. More precisely, she showed that every quantum symmetric pair coideal subalgebra of finite type admits a polynomial subalgebra which specializes to the corresponding classical Cartan subalgebra as $q \mapsto 1$. In the $\mathsf{AI}$ case, this Cartan subalgebra is generated by the subset of the canonical generators with odd subscripts, i.e., $B_1, B_3, B_5$, etc., and weights are defined as eigenvalues of these generators. Nevertheless, the Gelfand--Tsetlin basis is generally \emph{not} a weight basis in this sense. 
Secondly, our result is, in a sense, a dual version of \eqref{eq: into hwp}, since it involves an anti-dominant monomial $M_+$ and $M_p$'s which are products of factors $\ourA_{i,c}$. 

An interesting new feature of boundary $q$-characters is the presence of an extra symmetry (Corollary \ref{cor: symmetry cnst orbits}). The source of this symmetry is a natural $(\Z/2\Z)^{\lfloor n/2 \rfloor}$-action on the set of Gelfand--Tsetlin patterns by flipping the signs of entries in the last column. The eigenvalues of the operators $\Theta_{i,m}$ are constant along the orbits of this action, which is reflected by the multiplicity with which the corresponding monomials occur in the boundary $q$-character. The fact that the spectrum is not simple, even on  evaluation modules, presents new challenges and suggests that some refinement or deformation of the boundary $q$-character may be needed to build a more robust theory. 

Finally, we remark that the formulae for the action of the Lu--Wang generators obtained in this paper bear an obvious resemblance to the formulae describing the GKLO (Gerasimov--Kharchev--Lebedev--Oblezin) representations of shifted  quantum symmetric pairs discovered recently in \cite{LPgklo, LWW}. We expect a precise relationship between the two can be established via parabolic Verma modules, following the method of \cite[\S 2.7]{FPT}.

\addtocontents{toc}{\SkipTocEntry} 
	
\section*{Acknowledgements} 
We would like to thank Andrea Appel, Stefan Kolb, Catharina Stroppel and Bart Vlaar for helpful discussions and suggestions. 
We are also grateful to the International Centre for Mathematical Sciences (ICMS) in Edinburgh for hosting the programme \emph{``Towards q­-characters for quantum affine symmetric pairs''}, during which part of this work was carried out.

\section{Split quantum symmetric pairs} 
\label{sec: QSP}
\nc{\sln}{\mathfrak{sl}_{n+1}}
\nc{\asln}{\widehat{\mathfrak{sl}}_{n+1}}

We work over the field of complex numbers and assume that $q \in \C^\times$ is not a root of unity throughout. We also fix a square root $q^{\frac{1}{2}}$ of $q$, and a square root ${\rm i}$ of $-1$. Moreover, we assume that the rank is at least two, i.e., $n \geq 2$. 

For $k \in \frac{1}{2}\Z$, set 
\[
[k] = \frac{q^k-q^{-k}}{q-q\mi}, \qquad \{k\} = q^k + q^{-k}, \qquad 
[k]_+ = {\rm i} \frac{q^k+q^{-k}}{q-q\mi}. 
\]
Note that $[k]_+$ does not have a well-defined limit as $q \mapsto 1$.

\subsection{Quantum affine $\sln$} 

Let $\indx_0 = \{ 1,\cdots, n \}$ and $\indx = \indx_0 \cup \{0\}$. We consider the Lie algebras $\sln$ and $\asln$, with corresponding Cartan matrices  $(a_{ij})_{i,j \in \indx_0}$ and $(a_{ij})_{i,j \in \indx}$, respectively. 
We use standard notations and conventions regarding root systems, weight lattices, Weyl groups, etc., as in, e.g., \cite[\S 3.1]{lu-wang-21}. In particular, we 
let $\alpha_i$ ($i \in \indx$) denote the simple roots of $\asln$; let $\theta$ denote the highest root, and $\delta$ the basic imaginary root; let $P$ and $Q$ denote the weight and root  lattices of $\sln$, respectively; and let $\omega_i \in P$ ($i \in \indx_0$) be the fundamental weights of~$\sln$.

The \emph{quantum affine algebra} $\qaa$ is the algebra with generators $\Eg_i^\pm, \Kg_i^{\pm 1}$ $(i \in \indx)$ and relations: 
\begin{align}
\Kg_i\Kg_i^{-1} &= \Kg_i^{-1}\Kg_i = 1, & 
\Kg_i\Kg_j &= \Kg_j \Kg_i, \\
\Kg_i\Eg_j^{\pm} &= q^{\pm a_{ji}} \Eg_j^\pm\Kg_i, & 
[\Eg_i^{+}, \Eg_j^{-}] &= \delta_{ij} \frac{\Kg_i - \Kg_i^{-1}}{q - q^{-1}}, 
\end{align}
\[
\Serre_{ij}(\Eg_i^{\pm}, \Eg_j^\pm) = 0 \qquad (i \neq j),
\] 
where
\[
\Serre_{ij}(x,y) = \begin{cases}
x^2 y - [2] xyx + yx^2 & \mbox{if} \ a_{ji} = -1, \\
xy - yx & \mbox{if} \ a_{ji} = 0.  
\end{cases}
\]

\nc{\uaff}{\mathbf{U}_{\mathsf{aff}}}

The algebra $\qaa$ is a Hopf algebra, with the coproduct
\eq \label{eq: coprod on U}
\Delta(\Eg_i^+) = \Eg_i^+ \otimes 1 + \Kg_i \otimes \Eg_i,  \quad \Delta(\Eg_i^-) = \Eg_i^- \otimes \Kg_i^{-1} + 1 \otimes \Eg_i^-, \quad \Delta(\Kg_i^{\pm 1}) = \Kg_i^{\pm 1} \otimes \Kg_i^{\pm 1},  
\eneq
and the counit given by $\varepsilon(\Eg_i^{\pm}) = 0$ and $\varepsilon(\Kg_i^{\pm 1}) = 1$.

The \emph{quantum loop algebra} $\uaff = \qla$ is the quotient of $\qaa$ by the ideal generated by the central element $\prod_{i \in \indx} \Kg_i - 1$. 
Let us recall the ``new" Drinfeld presentation of the quantum loop algebra $\qla$. By \cite{drinfeld-dp, beck-94}, $\qla$ is isomorphic to the algebra generated by $\xgpm_{i,k}, \hg_{i,l}, \Kg^{\pm 1}_{i}$, where $k \in \Z$, $l \in \Z - \{0\}$ and $i \in \indx_0$, subject to the following relations:
\begin{align}
\Kg_i\Kg_i^{-1} =& \ \Kg_i^{-1}\Kg_i = 1, & 
\Kg_i \Kg_j =& \ \Kg_j \Kg_i, \\ 
[\hg_{i,k}, \hg_{j,l}] =& \ 0, & 
\Kg_i\hg_{j,k} =& \ \hg_{j,k}\Kg_i, \\ 
\Kg_i \xgpm_{j,k} =& \ q^{\pm a_{ji}} \xgpm_{j,k}\Kg_i,  &
[\hg_{i,k},\xgpm_{j,l}] =& \pm \textstyle \frac{[k\cdot a_{ji}]}{k}\xgpm_{j,k+l},\\
[\xgpm_{i,k+1}, \xgpm_{j,l}]_{q^{\pm a_{ji}}}  =& \ q^{\pm a_{ji}}[\xgpm_{i,k},\xgpm_{j,l+1}]_{q^{\mp a_{ji}}} &
[\xgp_{i,k}, \xgm_{j,l}] =& \ \delta_{ij} \textstyle\frac{1}{q-q^{-1}}(\psig_{i, k+l} - \phig_{i, k+l}), 
\end{align}
\begin{align}
\Sym_{k_1, k_2} \big( \xgpm_{j,l} \xgpm_{i,k_1} \xgpm_{i,k_2} - [2]  \xgpm_{i,k_1} \xgpm_{j,l} \xgpm_{i,k_2} +  \xgpm_{i,k_1} \xgpm_{i,k_2} \xgpm_{j,l} \big) =& \ 0 \quad (a_{ji} = -1), \\ 
[\xgpm_{i,k}, \xgpm_{j,l}] =& \ 0 \quad (a_{ji} = 0), 
\end{align} 
where $\Sym_{k_1, k_2} $ denotes symmetrization with respect to the indices $k_1, k_2$ and 
\[
\pmb{\phi}^{\pm}_i(z) = \sum_{k=0}^\infty \phipm_{i,\pm k} z^{\pm k} = \Kg_i^{\pm 1} \exp\left( \pm (q-q^{-1}) \sum_{k=1}^\infty \hg_{i,\pm k} z^{\pm k} \right). 
\] 

%%%%%%%%%%%%%%%%

\subsection{Quantum symmetric pairs of split type}

\nc{\Uiv}[1]{\Ui_{#1}}
\nc{\auic}{\mathbf{U}^\imath_{\mathsf{aff},\mathbf{c}}}
\nc{\uic}{\mathbf{U}^\imath_{\mathbf{c}}}
\nc{\DrUi}{{}^{\mathrm{Dr}}\aui}
\nc{\ui}{\mathbf{U}^\imath}
\nc{\MRSUi}{{}^{\mathrm{MRS}}\Ui}
\nc{\MRSaui}{{}^{\mathrm{MRS}}\aui}

In this paper, we consider quantum symmetric pair coideal subalgebras of $U_q(\sln)$ and $U_q(\asln)$ of split type, i.e., corresponding to Satake diagrams with no black nodes and no involution. They are also known as quantum symmetric pairs of type $\mathsf{AI}$. Assume $n \geq 2$. 

\Defi
Let $\aui = \auiv{n+1} = \Ui(\asln)$ be the algebra generated by $\Bg_i$ and central invertible elements $\KK_i$ $(i \in \indx)$ subject to relations 
\eq \label{eq: QSP Serre}
-q \KK_i\mi\Serre_{ij}(B_i,B_j) = 
\left\{
\begin{array}{r l l}
0 & \mbox{if} & a_{ji} = 0, \\
\Bg_j  & \mbox{if} & a_{ji} = -1. 
\end{array}
\right.
\eneq
We denote the subalgebra generated by $\Bg_i, \KK_i$ $(i \in \indx_0)$ as $\Ui = \Uiv{n+1} = \Ui(\sln)$. 
\enDefi

For $\mu = \sum_{i \in \indx} d_i \alpha_i \in \Z\indx = \oplus_{i \in \indx} \Z \alpha_i$, we will also, as usual, write $\KK_\mu = \prod_{i \in \indx} \KK_i^{d_i}$. 

Given $\mathbf{c} = (c_0, \cdots, c_n) \in (\C^\times)^{n+1}$, let $\auic$ (resp.\ $\uic$) be the quotient of $\aui$ (resp.\ $\Ui$) by the two-sided ideal generated by $q^{-2}\KK_i - c_i$, for $i \in \indx$ (resp.\ $i \in \indx_0$). 
By \cite[Theorem 7.1]{kolb-14}, for any $\mathbf{s} = (s_0, \cdots, s_n) \in \C^{n+1}$, there exists an injective algebra homomorphism 
\eq
\label{eq: Kolb emb}
\eta \colon \auic \monoto \qla, \quad \Bg_i \mapsto \Eg_i^- - c_i \Eg_i^+ \Kg_i^{-1} + \sss_i \Kg_i^{-1} \quad 
(i \in \indx), 
\eneq 
which restricts to a homomorphism $\uic \hookrightarrow U_q(\sln)$. The monomorphism $\eta$ endows $\auic$ and $\uic$ with the structures of coideal subalgebras. Explicitly, 
\eq \label{eq: coproduct explicit} 
\Delta(B_i) = 1 \otimes \eta(B_i) + \eta(B_i) \otimes K_i\mi. 
\eneq 

\subsection{Braid group action} 

Let $W_0 = \langle s_i \mid i \in \indx_0 \rangle$ be the Weyl group of $\sln$, and $W = \langle s_i \mid i \in \indx \rangle = W_0 \ltimes Q$ the affine Weyl group. The extended affine Weyl group is
\[
\widetilde{W} = W_0 \ltimes P = \Lambda \ltimes W, 
\] 
where $\Lambda = P/Q = \langle \pi \rangle$ is the cyclic group of automorphisms of the Dynkin diagram of $\asln$, generated by $\pi(i) = i+1 \mod n + 1$. 
Let $\mathcal{B}$ and $\widetilde{\mathcal{B}}$ be the braid groups associated to $W_0$ and $\widetilde{W}$, respectively. 

By \cite[Lemma 3.4]{lu-wang-21}, for each $i \in \indx$, there exists an automorphism $\Tbr_i$ of $\aui$ such that $\Tbr_i(\KK_\mu) = \KK_{s_i\mu}$ and 
\[
\Tbr_i(B_j) = 
\left\{ \begin{array}{l l}
\KK_j\mi B_j & \mbox{if } i=j, \\[2pt]
B_j & \mbox{if } a_{ji} = 0, \\[2pt]
B_jB_i - q_iB_iB_j  & \mbox{if } a_{ji} = -1, \\[2pt]
\end{array} \right. 
\] 
for $\mu \in \Z\indx$ and $j \in \indx$. 
For each $\lambda \in \Lambda$, there is also an automorphism $\Tbr_\lambda$ such that $\Tbr_\lambda(B_i) = B_{\lambda(i)}$ and $\Tbr_\lambda(\KK_i) = \KK_{\lambda(i)}$.   
These automorphisms define an action of the braid group $\widetilde{\mathcal{B}}$ on $\aui$, which restricts to an action of $\mathcal{B}$ on $\Ui$. Given a reduced expression $w = \lambda s_{i_1} \cdots s_{i_k} \in \widetilde{W}$, we set $\Tbr_{w} = \Tbr_\lambda \Tbr_{i_1} \cdots \Tbr_{i_k}$. 

\subsection{The Lu--Wang presentation} 

Let us recall the ``Drinfeld-type" presentation of $\aui$ due to Lu and Wang \cite{lu-wang-21}. 
We will use the following shorthand notations: 
\begin{align}
S(r_1, r_2 \mid s; i,j) =& \ \Ag_{i,r_1} \Ag_{i,r_2} \Ag_{j,s} - [2] \Ag_{i,r_1} \Ag_{j,s} \Ag_{i,r_2} + \Ag_{j,s}\Ag_{i,r_1} \Ag_{i,r_2}, \\
\mathbb{S}(r_1, r_2 \mid s; i,j) =& \ S(r_1, r_2 \mid s; i,j) + S(r_2, r_1 \mid s; i,j), \\ 
R(r_1, r_2 \mid s; i,j) =& \ \KK_i\Ck^{r_1}\big( - \sum_{p \geq 0} q^{2p} [2] [\Thg_{i, r_2-r_1-2p-1}, \Ag_{j,s-1}]_{q^{-2}}  \Ck^{p+1} \\
& \ - \sum_{p \geq 1} q^{2p-1} [2] [\Ag_{j,s}, \Thg_{i, r_2-r_1-2p}]_{q^{-2}} \Ck^{p} - [\Ag_{j,s}, \Thg_{i,r_2-r_1}]_{q^{-2}}\big), \\
\mathbb{R}(r_1, r_2 \mid s; i,j) =& \ R(r_1, r_2 \mid s; i,j) + R(r_2, r_1 \mid s; i,j). 
\end{align}

\Defi 
Let $\DrUi$ be the $\kor$-algebra generated by $\Hg_{i,m}$ and $\Ag_{i,r}$, where $m\geq1$, $r\in\Z$, and invertible central elements $\KK_i$ ($i \in \indx_0$), $\CCC$, subject to the following relations:
\begin{align} 
[\Hg_{i,m},\Hg_{j,l}] &=0, \label{eq: grel1} \\
[\Hg_{i,m}, \Ag_{j,r}] &= \textstyle \frac{[m \cdot a_{ji}]_{q}}{m} (\Ag_{j,r+m}- \Ag_{j,r-m}\CCC^m), \label{eq: grel2} \\ 
[\Ag_{i,r}, \Ag_{j,s}] &= 0 \quad \mbox{if} \quad a_{ji} = 0, \label{eq: grel3} \\ 
[\Ag_{i,r}, \Ag_{j,s+1}]_{q^{-a_{ji}}}   &= q^{-a_{ji}} [\Ag_{i,r+1}, \Ag_{j,s}]_{q^{a_{ji}}}  \quad \mbox{if} \quad i\neq j, \label{eq: grel4} \\ 
\label{eq: grel5}
[\Ag_{i,r}, \Ag_{i,s+1}]_{q^{-2}}  -q^{-2} [\Ag_{i,r+1}, \Ag_{i,s}]_{q^{2}}
&= q^{-2}\KK_i\CCC^r\Thg_{i,s-r+1} - q^{-4}\KK_i\CCC^{r+1}\Thg_{i,s-r-1}  \\ \notag
&\quad  + q^{-2}\KK_i\CCC^s\Thg_{i,r-s+1}  -q^{-4}\KK_i \CCC^{s+1}\Thg_{i,r-s-1}, \\
\mathbb{S}(r_1, r_2 \mid s; i,j) &= \mathbb{R}(r_1, r_2 \mid s; i,j) \quad \mbox{if} \quad a_{ji} = -1, \label{eq: grel6}
\end{align}
where $m,n\geq1$; $r,s, r_1, r_2\in \Z$ and 
\[
1+ \sum_{m=1}^\infty (q-q^{-1})\Thg_{i,m} z^m  =  \exp \left( (q-q^{-1}) \sum_{m=1}^\infty \Hg_{i,m} z^m \right).
\]
By convention, $\Thg_{i,0} = (q-q^{-1})^{-1}$ and $\Thg_{i,m} = 0$ for $m\leq -1$. Moreover, we define the following two commutative subalgebras of $\aui$:
\begin{equation} \label{eq: HT Cartans}
\mathcal{H} = \langle \Theta_{i,r} \mid i \in \indx_0, r \in \Z_{\geq0} \rangle, \qquad \mathcal{T} = \langle B_i \mid i \mbox{ is odd } \rangle. 
\end{equation}
We propose to call $\mathcal{H}$ the \emph{Lu--Wang Cartan subalgebra}, and $\mathcal{T}$ \emph{Letzter's Cartan aubalgebra}. 
\enDefi 

We choose the sign function 
\[
o(\cdot) \colon \indx_0 \to \{\pm 1\}, \qquad i \mapsto (-1)^{i+1}. 
\]
By \cite[Theorem 3.13]{lu-wang-21}, there is an algebra isomorphism 
\eq
\label{eq: droq->oq}
\DrUi \isoto{} \aui, \qquad A_{i,r} \mapsto o(i)^r \Tbr_{\omega_i}^{-r}(B_i), \quad
\KK_i \mapsto \KK_i, \quad  \Ck \mapsto \KK_\delta, 
\eneq
for $i \in \indx_0$. 

We will use the following notation for generating series: 
\[
\avar{i} = \sum_{r \geq 0} A_{i,r} z^r, \qquad \thvar{i} = \sum_{r \geq 0} \grave{\Theta}_{i,r} z^r = \frac{(q-q\mi)(1-q^{-2}\CCC z^2)}{1 - \CCC z^2} \sum_{r \geq 0} \Theta_{i,r} z^r. 
\]
By \cite[Proposition 3.1]{Przez-23}, 
\begin{equation} \label{eq: Theta from A}
\thvar{i} = 1 + \frac{q^2(q-q\mi)\CCC \KK_i\mi z^2}{1 - \CCC z^2}\left(z\mi[A_{i,-1}, \avar{i}]_{q^{-2}} - q^{-2}[A_{i,0}, \avar{i}]_{q^2} \right). 
\end{equation} 

\subsection{The $q$-Onsager algebra} 

The $q$-Onsager algebra is a split quantum symmetric pair coideal subalgebra of $U_q(\widehat{\mathfrak{sl}}_2)$. It was first introduced by Terwilliger in \cite{Ter01}.  
As in \cite{lu-wang-21}, we consider a central extension of the $q$-Onsager algebra called the \emph{universal} $q$-Onsager algebra. It has the following definition by generators and relations. 

\Defi
The universal $q$-Onsager algebra $\Oq$ is the algebra generated by $B_0, B_1$ and invertible central elements $\KK_0, \KK_1$ subject to relations:
\eq
\label{eq: Ons rel}
\sum_{r=0}^3 (-1)^r \sqbin{3}{r} \Bg_i^{3-r}\Bg_j\Bg_i^r = - q\mi \KK_i [2]^2[\Bg_i, \Bg_j] \qquad (0 \leq i \neq j \leq 1).
\eneq 
\enDefi

For $i \in \indx_0$, let $\Oqi$ be the subalgebra of $\aui$ generated by $B_i$, $\Tbr_{\omega'_i}(B_i)$, $\KK_i^{\pm1}$ and $(\CCC \KK_i^{-1})^{\pm1}$. By \cite[Proposition 3.9]{lu-wang-21}, there is an algebra isomorphism $\iota \colon \Oq \to \Oqi$ sending 
\eq \label{eq: rank 1 q onsager}
B_1 \mapsto B_i, \quad B_0 \mapsto \Tbr_{\omega'_i}(B_i), \quad \KK_1 \mapsto \KK_i, \quad \KK_0 \mapsto \KK_\delta \KK_i^{-1}, 
\eneq 
where $\omega'_i = \omega_i s_i$.

\subsection{Assumptions about parameters}
\label{subsec: parameters}

Most of the time, we will work with a specific set of parameters, in accordance with the conventions of \cite{Wenzl}. 
We abbreviate $\ui = \ui_{\mathbf{c}}$, for the particular choice of $\mathbf{c} = (-q\mi, \cdots, -q\mi) \in (\C^\times)^n$, which is equivalent to specializing $\KK_i = -q$ for each $i \in \indx_0$. In the affine case, we additionally fix a non-zero complex number $a \in \C^\times$, and let $\aaui$ be the central reduction of $\aui$ by the relations: $\KK_i = -q$ ($i \in \indx_0$) and $\mathfrak{C} = a^2$. 

%%%%%%%%%%%%%%%%

\section{Boundary $q$-characters}
\label{sec: boundary q}

\nc{\Uic}{\mathbf{U}^\imath_{\mathbf{c}}}
\nc{\Uh}{U_q(\widetilde{\mathfrak{h}})}

\subsection{$q$-characters of quantum affine algebras}

Let $\Rep \uaff$ be the monoidal category of finite dimensional representations of $\uaff = \qla$. It acts on $\Rep \aaui$, the category of finite dimensional representations of $\aaui$, via the coproduct \eqref{eq: coproduct explicit}. 
Passing to the Grothendieck groups of the respective categories, $[\Rep \uaff]$ is a ring acting on $[\Rep \aaui]$. 
Let $\Uh$ be the subalgebra of $\uaff$ generated by $h_{i,k}$ ($i \in \indx$, $k < 0$). Following \cite{FrenRes, FrenMuk-comb}, set 
\begin{align} \label{eq: defi of Yia}
Y_{i,b} =& \  K_{\omega_i}^{-1} \exp \left( -(q-q\mi) \sum_{k > 0} \tilde{h}_{i,-k} b^k z^k  \right) \in \Uhz \qquad (b \in \C^\times), 
\end{align}
where 
\eq
\tilde{h}_{i,-k} = \sum_{j \in \indx_0} \widetilde{C}_{ji}(q^k) \hg_{j,-k},
\eneq 
and $\widetilde{C}(q)$ is the inverse of the $q$-Cartan matrix (the natural quantization of the type $\mathsf{A}$ Cartan matrix with the diagonal $2$'s replaced by $[2]_q$). 
Let $\Yring = \Z[Y_{i,b}^{\pm 1}]_{i \in \indx_0, b \in \C^{\times}}$. 
By \cite[Theorem 3]{FrenRes}, there exists an injective ring homomorphism 
\[
\chmap \colon [\Rep \uaff] \to \Yring \subset \Uhz,
\]
called the \emph{$q$-character map}, given by 
\[
[V] \mapsto \Tr_V \left[ \exp \left( -(q-q\mi) \sum_{i \in \indx_0} \sum_{k >0} \frac{k}{[k]_{q}} \pi_V(h_{i,k}) \otimes \tilde{h}_{i,-k} z^k \right) \cdot (\pi_V \otimes 1)(T) \right], 
\]
where $T$ is as in \cite[(3.8)]{FrenRes}, and $\pi_V \colon \uaff \to \End(V)$ is the representation. The expression above derives from the Khoroshkin-Tolstoy-Levendorsky-Soibelman-Stukopkin-Damiani (KTLSSD) factorization of the universal $R$-matrix \cite{khor-tol, lev-sob-st, damiani-r}. 

For future use, we also set
\begin{align*} 
X_{i,b} &= K_{\omega_i}^{-1/2} \exp \left( -(q-q\mi) \sum_{k > 0} \frac{\tilde{h}_{i,-k}}{q^{-k/2}+q^{k/2}} b^k z^k  \right), \\
\widetilde{X}_{i,b} &= K_{\omega_i}^{-1/2} \exp \left( -(q-q\mi) \sum_{k > 0} \frac{\tilde{h}_{i,-k}}{q^{-k/2}+(-1)^kq^{k/2}} b^k z^k  \right). 
\end{align*}
Note that $Y_{i,b} = X_{i,bq^{-1/2}} X_{i,bq^{1/2}} = \widetilde{X}_{i,bq^{-1/2}} \widetilde{X}_{i,-bq^{1/2}}$.

\subsection{Algebraic interpretation} 

By \cite[Proposition 2.4]{FrenMuk-comb}, the $q$-character map can equivalently be defined more explicitly in terms of the joint spectrum of the Drinfeld--Cartan operators (i.e., the coefficients of the series $\pmb{\phi}^{\pm}_i(z)$). More precisely, there is a one-to-one correspondence between the monomials occurring in $\chmap(V)$ and the common eigenvalues of $\pmb{\phi}^{\pm}_i(z)$ on $V$. 
Let 
$V = \bigoplus_{\gamma} V_{\gamma}$, with $\gamma = (\gamma^\pm_{i,\pm m})_{i \in \indx_0, m \in \Z_{\ge 0}}$, where 
\begin{align*}
V_{\gamma} = \{ v \in V \mid \exists p \ \forall i \in \indx_0 \ \forall m \in \Z_{\ge 0}:  (\psi^\pm_{i,\pm m} - \gamma^\pm_{i,\pm m})^p \cdot v = 0 \}, 
\end{align*}
be the generalized eigenspace decomposition of $V$. 
Collect the eigenvalues into generating series $\gamma^\pm_i(z) = \sum_{m \ge 0} \gamma^\pm_{i,\pm m} z^{\pm m}$. 
By \cite[Proposition 2.4]{FrenMuk-comb}, the series $\gamma^\pm_i(z)$ are expansions (at $0$ and $\infty$, respectively) of  the same rational function of the form
\[
q^{\deg Q_i - \deg R_i} \frac{Q_i(q^{-1}z)R_i(qz)}{Q_i(qz)R_i(q^{-1}z)}, 
\]
for some polynomials $Q_i(z), R_i(z)$ with constant term $1$. Writing  
\eq \label{eq: QR notation}
Q_i(z) = \prod_{r=1}^{k_i} (1 - z a_{i,r}), \qquad R_i(z) = \prod_{s=1}^{l_i} (1 - z b_{i,s}), 
\eneq
the $q$-character of $V$ can now be expressed as 
\[
\chi_q(V) = \sum_{\gamma} \dim(V_\gamma) M_\gamma, \qquad 
M_\gamma = \prod_{i \in \indx_0} \prod_{r=1}^{k_i} Y_{i,a_{i,r}} \prod_{s=1}^{l_i} Y_{i,b_{i,s}}\mi. 
\]

\subsection{Boundary $q$-characters} 

In \cite{Przez-23, LP24, LP26}, we defined the notion of a \emph{boundary $q$-character} of a $\aaui$-module in terms of the spectra of the Drinfeld--Cartan operators in the Lu--Wang (Drinfeld-type) presentation. Let us recall the corresponding  definitions and results. We need the following normalization: 
\[
\mathbf{\grave{H}}_i(z) = \sum_{r \geq 1} \grave{H}_{i,r} z^r = (q-q\mi)\mi \on{log} \left( \thvar{i} \right). 
\] 

\Defi
Let  
\[
\mathcal{K}^0 =  \exp \left( -(q-q\mi) \sum_{i \in \indx_0} \sum_{k >0} \frac{k}{[k]_{q}} \grave{H}_{i,k} \otimes \tilde{h}_{i,-k} z^k \right) \in \fext{\aaui \otimes \Uh}{z}. 
\]
We define the \emph{boundary $q$-character map} to be 
\eq \label{eq:iq char defin}
\ichmap \colon [\Rep \aaui] \ \to \ \Uhz, \qquad [V] \mapsto \Tr_V((\pi_V \otimes 1)(\mathcal{K}^0)). 
\eneq 
\enDefi

Consider $\Uhz$ as a $\Yring$-module via the ring homomorphism
\eq \label{eq: Yring module}
\Yring \to \Yring \hookrightarrow \Uhz, \qquad Y_{i,b} \mapsto Y_{i,a^2b}Y_{i,b\mi}\mi. 
\eneq
To state the next theorem, we first need to introduce some notation. Given a polynomial $P(z) \in \C[z]$ with constant term $1$, let $P^\dag(z)$ be the polynomial with constant term~$1$ whose roots are obtained from those of $P(z)$ via the transformation $x \mapsto a^{-2} x\mi$; and let $P^*(z)$ be the polynomial with constant term $1$ whose roots are the inverses of the roots of $P(z)$. Here $a$ is the same as the non-zero scalar fixed in \S \ref{subsec: parameters}. 

\Thm
\label{cor: FR thm Oq} 
\begin{enumerate} 
\item
Let $W \in \Rep \uaff$. 
Then the generalized eigenvalues of $\pmb{\grave{\Thg}}_i(z)$ on the restricted representation $\eta^*(W)$   are of the form: 
\eq \label{eq: FR formula Oq}
\gamma_i^\iota(z) = \frac{\ourQ_i(q^{-1}z)}{\ourQ_i(qz)} \frac{\ourQ_i^\dag(qz)}{\ourQ_i^\dag(q^{-1}z)}, 
\eneq
where $\ourQ_i(z)$ is a polynomial with constant term $1$. Explicitly, 
\[
\ourQ_i(z) = {Q_i(a^2 z)R_i^*(z)}, \quad  \ourQ_i^\dag(z) = {R_i(a^2 z)Q_i^*(z)}. 
\]
There is a one-to-one correspondence between the common eigenvalues of $\pmb{\grave{\Thg}}_i(z)$ and the monomials in $Y_{i,b}^{\pm1}$ occurring in $\ichmap(\eta^*(W))$. 
\item The following diagram commutes: 
\[
\begin{tikzcd}[ row sep = 0.2cm]
{[}\Rep \uaff]  \arrow[r, "\chmap"] & \Yring  \\
 \curvearrowright  & \curvearrowright  \\
{[}\Rep \aaui] \arrow[r, "\ichmap"] & \Uhz. 
\end{tikzcd}
\]
\end{enumerate} 
\enthm

\Proof 
See \cite[Corollary 5.1]{Przez-23}   and \cite[Corollary 9.3]{LP24}. 
\enproof

\subsection{Dominant monomials} 

The correspondence between the common eigenvalues of $\pmb{\grave{\Thg}}_i(z)$ and the monomials occurring in the boundary $q$-character of a $\aaui$-module can be written more compactly if we switch to the following variables: 
\eq \label{eq: Y, X defin} 
\ourY_{i,q^k} = Y_{i,aq^k}Y_{i,aq^{-k}}\mi, \quad \ourX_{i,q^{k/2}} = X_{i,aq^{k/2}}X_{i,aq^{-k/2}}\mi, \quad
\widetilde{\ourX}_{i,\pm q^{k/2}} = \widetilde{X}_{i,\pm aq^{k/2}}\widetilde{X}_{i,\mp aq^{-k/2}}\mi, 
\eneq 
for $k \in \Z$, and a fixed $a \in \C^\times$ as in \S \ref{subsec: parameters}. 
Note that 
\[\ourX_{i,q^{k/2}} = \ourX_{i,q^{-{k/2}}}\mi, \quad \ourY_{i,q^k} = \ourX_{i,q^{k-\frac{1}{2}}} \ourX_{i,q^{k+\frac{1}{2}}}, \quad \ourX_{i,1} = \ourY_{i,1} = 1.
\]
 We also have the relations $\widetilde{\ourX}_{i,q^{k/2}} = \widetilde{\ourX}_{i,-q^{-{k/2}}}\mi$ and $\ourY_{i,q^k} = \widetilde{\ourX}_{i,q^{k-\frac{1}{2}}} \widetilde{\ourX}_{i,-q^{k+\frac{1}{2}}}$. Therefore, we can (and will) rewrite any monomial in the variables $\ourX_{i,q^{k/2}}$ (resp.\ $\widetilde{\ourX}_{i,\pm q^{k/2}}$) uniquely as a monomial in the variables $\ourX_{i,q^{k/2}}^{\pm 1}$ (resp.\ $\widetilde{\ourX}_{i,\pm q^{k/2}}^{\pm1}$) for $k \in \Z_{\geq 1}$ (resp.\ $k \in \Z_{\geq 0}$). We call the latter monomials \emph{standardized}. 
Let 
\[
\Xmod_a = \Z[\ourX_{i,q^{k/2}}^{\pm1}]_{i \in \indx_0, k \geq 1}, \qquad \widetilde{\Xmod}_a = \Z[\widetilde{\ourX}_{i,\pm q^{k/2}}^{\pm1}]_{i \in \indx_0, k \geq 0}
\]

\Defi
We say that a standardized monomial is \emph{dominant} (resp.\ \emph{anti-dominant}) if it contains no negative (resp.\ positive) powers of $\ourX_{i,q^{k/2}}$ or $\widetilde{\ourX}_{i,\pm q^{k/2}}$. 
\enDefi

Recall 
\begin{align} \label{eq: defi of Yia}
A_{i,b} =& \  K_{i}^{-1} \exp \left( -(q-q\mi) \sum_{k > 0} h_{i,-k} b^k z^k  \right) = 
Y_{i, qb}Y_{i, q^{-1}b} \prod_{a_{ij} = -1} Y_{j, b}\mi 
\qquad (b \in \C^\times) 
\end{align}
from \cite{FrenRes, FrenMuk-comb}. 
Set  
\[
\ourA_{i,q^k} = A_{i,aq^k}A_{i,aq^{-k}}\mi 
 \qquad (k \in \Z). 
\]
Note that $\ourA_{i,q} = \ourY_{i,q^2} \prod_{a_{ij} = -1} \ourY_{j, q}\mi $ and $\ourA_{i,1} = 1$.  
We can use the variables $\ourA_{i,q^k}$ to define a partial order on the set of standardized monomials.

\Defi \label{defi: monomial ordering}
We define a partial order on the set of standardized monomials by the rule: $M \leq M'$ if and only if $M'/M$ can be expressed as a monomial in $\ourA_{i, q^{k}}$ for $k \geq 0$. 
\enDefi

We end this section by summarizing the relationship between eigenvalues and monomials. 

\Lem
\label{lem: monomial vs eigen} 
Let $V$ be a finite-dimensional $\aaui$-module, and $v \in V$ a generalized eigenvector for $\mathcal{H}$. 
Suppose that the series of generalized eigenvalues of $\pmb{\grave{\Thg}}_i(z)$ on $v$ is given by a rational function $\gamma_i^\imath(z)$, and let $M$ be the associated monomial in $\ichmap(V)$. Let $P_i(z) = (1-q^kaz)$ and $Q_i(z) = (1-q^{k/2}az)$ with $k \geq 1$. Then there is the following correspondence between factors of $\gamma_i^\imath(z)$ and factors of $M$:   
\begin{align*} 
\frac{P_{i}(q\mi z)}{P_{i}(q z)} \frac{P^\dag_{i}(q z)}{P^\dag_{i}(q\mi z)} &= \frac{(1-q^{k-1}az)(1-q^{1-k}az)}{(1-q^{k+1}az)(1-q^{-k-1}az)}&  &\longleftrightarrow&  \ourY_{i,q^k}, \\ 
\frac{P_{i}(q z)}{P_{i}(q\mi z)} \frac{P^\dag_{i}(q\mi z)}{P^\dag_{i}(q z)} &= \frac{(1-q^{k+1}az)(1-q^{-k-1}az)}{(1-q^{k-1}az)(1-q^{1-k}az)}&  &\longleftrightarrow&  \ourY_{i,q^k}\mi, 
\end{align*} 
%%%%%%%
as well as 
\begin{align*}
\frac{Q_{i}(q^{-\frac12} z)}{Q_{i}(q^{\frac12} z)} \frac{Q^\dag_{i}(q^{\frac12} z)}{Q^\dag_{i}(q^{-\frac12} z)} &= \frac{(1-q^{\frac{k-1}{2}}az)(1-q^{\frac{1-k}{2}}az)}{(1-q^{\frac{k+1}{2}}az)(1-q^{\frac{-k-1}{2}}az)}&  &\longleftrightarrow&  \ourX_{i,q^{\frac{k}{2}}}, \\ 
\frac{Q_{i}(q^{-\frac12} z)}{Q_{i}(-q^{\frac12} z)} \frac{Q^\dag_{i}(q^{\frac12} z)}{Q^\dag_{i}(-q^{-\frac12} z)} &= \frac{(1-q^{\frac{k-1}{2}}az)(1-q^{\frac{1-k}{2}}az)}{(1+q^{\frac{k+1}{2}}az)(1+q^{\frac{-k-1}{2}}az)}&  &\longleftrightarrow&  \widetilde{\ourX}_{i,q^{\frac{k}{2}}}. 
\end{align*} 
\enlem 

\Proof
The proof is analogous to the argument in the proof of \cite[Proposition 2.4]{FrenMuk-comb}. 
\enproof

%%%%%%%%%%%%%%%%%%%%%%%%%

\section{Representation theory of $\ui$} 
\label{subsec: rep th ui}

\subsection{Weights} 
\nc{\croot}{\Delta_{cl}}
\nc{\rank}{\on{rank}}
\nc{\calpha}{\dot{\alpha}} 
\nc{\cl}{\on{cl}}

Let $\croot$ be the root system of $\mathfrak{so}_{n+1}$, i.e., 
\[
\croot = \mbox{ root system of type } 
\left\{ \begin{array}{l l}
\mathsf{D}_{(n+1)/2} & \mbox{ if } n \mbox{ is odd}, \\[1pt]
\mathsf{B}_{n/2} & \mbox{ if } n \mbox{ is even}. 
\end{array} \right. 
\]
Set $\indx_{cl} = \{1, \cdots, \rank(\croot) \}$. Let 
$P_{cl}$ and $Q_{cl}$ denote the corresponding weight and root lattices, respectively; and let $\pi_i \in P_{cl}$ ($i \in \indx_{cl}$) be the fundamental weights.  

Let $\{\epsilon_i\}$ be the standard basis of $\R^k$, where $r = \rank(\croot) = \lceil n/2 \rceil$. We identify roots and weights of $\mathfrak{so}_{n+1}$ with vectors in $\R^r$ in the usual way. More precisely, the simple roots are 
\[
\calpha_i = \epsilon_i - \epsilon_{i+1} \quad (1 \leq i \leq r-1), \qquad \calpha_r =  
\left\{ \begin{array}{l l}
\epsilon_{r-1} + \epsilon_r & \mbox{ if } n \mbox{ is odd}, \\[1pt]
\epsilon_r & \mbox{ if } n \mbox{ is even}, 
\end{array} \right. 
\]
and the fundamental weights are 
\[
\pi_i = \sum_{j=1}^i \epsilon_i, \qquad \pi_{r-1} =  
\left\{ \begin{array}{l l}
\big(\frac{1}{2} \sum_{j=1}^{r-1} \epsilon_j\big) - \frac{1}{2} \epsilon_r  & \mbox{ if } n \mbox{ is odd}, \\[2pt]
\sum_{j=1}^{r-1} \epsilon_j & \mbox{ if } n \mbox{ is even}, 
\end{array} \right. \qquad
\pi_r =  
\frac{1}{2} \sum_{j=1}^{r} \epsilon_j,
\]
for $1 \leq i \leq r-2$.

We now explain the connection between classical weights and the representation theory of $\ui$. 
We refer to generalized eigenvalues with respect to $B_1, B_3, \cdots, B_{2r -1}$ 
as \emph{weights}. There are two types of weights of interest to us. Given a tuple $(\lambda_i) \in \R^r$, set 
\[
\boldsymbol\lambda = ([\lambda_1], [\lambda_2], \cdots, [\lambda_r]), \qquad 
\boldsymbol\lambda_+ = ([\lambda_1]_+, [\lambda_2]_+, \cdots, [\lambda_r]_+). 
\] 
We call $\cl(\boldsymbol\lambda) = \cl(\boldsymbol\lambda_+) = \sum_{j=1}^r \lambda_j \epsilon_j$ the \emph{classical weight} corresponding to $\boldsymbol\lambda$ or $\boldsymbol\lambda_+$. 
A weight $\boldsymbol\lambda$ is called \emph{regularly dominant} if $\cl(\boldsymbol\lambda)$ is an integral dominant weight in $P_{cl}$. 
Note that, in that case, either $\lambda_i \in \Z$ or $\lambda_i \in \frac12 + \Z$, for all $i \in \indx_{cl}$ (all $\lambda_i$ are either integral or half-integral at the same time, but not mixed). In the former case, we will say that $\boldsymbol\lambda$ is  \emph{integral regularly dominant} and, in the latter case, that it is \emph{half-integral regularly dominant}. 
On the other hand, we will call the second type of weight, $\boldsymbol\lambda_+$, \emph{regularly dominant} if $\cl(\boldsymbol\lambda_+)$ is an integral dominant weight in $P_{cl}$, and, additionally, $\lambda_i \in \frac{1}{2} + \Z$, for all $i \in \indx_{cl}$. We denote the three sets of weights introduced above as $\Lambda^{\Z}_n$, $\Lambda^{\frac12 + \Z}_n$ and $\Lambda^+_n$, respectively; and put $\Lambda_n = \Lambda^{\Z}_n \cup \Lambda^{\frac12 + \Z}_n$.

\subsection{Verma construction} 
\label{subsec: Verma cnstr}

Abbreviate $r' = \lfloor n/2 \rfloor$. 
Given a regularly dominant weight~$\boldsymbol\lambda$, let $M_{\boldsymbol\lambda}$ be the quotient of $\ui$ by the left-sided ideal $I_{\boldsymbol\lambda}$ generated by 
\begin{equation} \label{eq: Verma rels}
B_{2i-1} - [\lambda_i] \quad (1 \leq i \leq r), \qquad (B_{2i-1} - [\lambda_i - 1])B_{2i} \quad (1 \leq i \leq r'). 
\end{equation} 
By \cite[Theorems 5.6, 5.9]{Wenzl}, $M_{\boldsymbol\lambda}$ is non-trivial and has a unique simple finite-dimensional quotient $L_{\boldsymbol\lambda}$. Simple modules arising in this manner are referred to as \emph{classical simple modules}. 
If $\cl(\boldsymbol\lambda) = \pi_m$ for some $m \in \indx_{cl}$, we call
\eq \label{eq: cl fun rep}
W_{\pi_m} := L_{\boldsymbol\lambda}
\eneq
the $m$\emph{-th fundamental representation}.

For any regularly dominant weight of the form $\boldsymbol\lambda_+$ and a sequence \linebreak $\underline{\varepsilon} = (\varepsilon(1), \cdots, \varepsilon(r)) \in \{\pm 1\}^r$, let $M_{\boldsymbol\lambda_+,{\underline{\varepsilon}}}$ be the quotient of $\ui$ by the left-sided ideal $I_{\boldsymbol\lambda_+,{\underline{\varepsilon}}}$  generated by 
\begin{equation*} 
B_{2i-1} - \varepsilon(i) [\lambda_i]_+ \quad (1 \leq i \leq r), \qquad (B_{2i-1} - \varepsilon(i) [\lambda_i - 1]_+)B_{2i} \quad (1 \leq i \leq r'). 
\end{equation*} 
By \cite[Corollary 5.7]{Wenzl}, $M_{\boldsymbol\lambda_+,{\underline{\varepsilon}}}$ is non-trivial and has a unique maximal \emph{semisimple} finite-dimensional quotient $N_{\boldsymbol\lambda_+,{\underline{\varepsilon}}}$, which decomposes as a direct sum of $2^{r'}$ representations of equal dimension. 
Simple modules arising in this manner are referred to as \emph{non-classical simple modules}.

\subsection{Classification}

Finally, let us recall the classification results from \cite{Ior-Kl} and \cite[Theorem 5.9]{Wenzl}. 
Given a dominant integral weight $\lambda \in P_{cl}$, let $L_\lambda$ be the simple $\mathfrak{so}_{n+1}$-module with highest weight $\lambda$. 

\Thm[{\cite{Ior-Kl, Wenzl}}]
Let $L$ be a simple finite-dimensional $\ui$-module. Then $L$ is either:
\begin{enumerate}
\item the classical simple module $L_{\boldsymbol\lambda}$ with highest regularly dominant weight $\boldsymbol\lambda$; or 
\item one of the $2^{r'}$-many non-classical simple modules in $N_{\boldsymbol\lambda_+,{\underline{\varepsilon}}}$, for a regularly dominant weight $\boldsymbol\lambda_+$ and a sequence ${\underline{\varepsilon}}$. 
\end{enumerate}
In case (1), $\dim L = \dim L_{\cl(\boldsymbol\lambda)}$, and, in case (2), $\dim L = 2^{-r'} \dim L_{\cl(\boldsymbol\lambda_+)}$. 
\enthm

\subsection{Gelfand--Tsetlin patterns}

\nc{\GT}{\on{GT}}
\nc{\GTn}{\on{GT}^{\on{nc}}}
\nc{\row}{\on{row}}
\nc{\col}{\on{col}}

We will also need the explicit construction of the simple modules due to Gavrilik, Iorgov and Klimyk \cite{GavKl, GI97, IK00}. 

\Defi
A \emph{Gelfand--Tsetlin pattern} (of orthogonal type) is a combinatorial object consisting of a set $\mathcal{M} = \{ m_{i,j} \mid 2 \leq j \leq n + 1, 1 \leq \lfloor j/2 \rfloor \}$ subject to the following constraints:
\begin{enumerate}[itemsep=3pt]
\item $\big(\forall m_{i,j} \in \mathcal{M}: m_{i,j} \in \Z \big)$ or $\big(\forall m_{i,j} \in \mathcal{M}: m_{i,j} \in \frac{1}{2}\Z \big)$, 
\item $m_{1, 2j} \geq m_{1,2j-1} \geq m_{2, 2j} \geq \cdots \geq m_{j-1,2j-1} \geq |m_{j,2j}|$, 
\item $m_{1, 2j+1} \geq m_{1,2j} \geq m_{2, 2j+1} \geq \cdots \geq m_{j-1,2j+1} \geq |m_{j,2j}|$. 
\end{enumerate} 
\enDefi
Graphically, the conditions defining a Gelfand--Tsetlin pattern can be presented using the following diagram:  
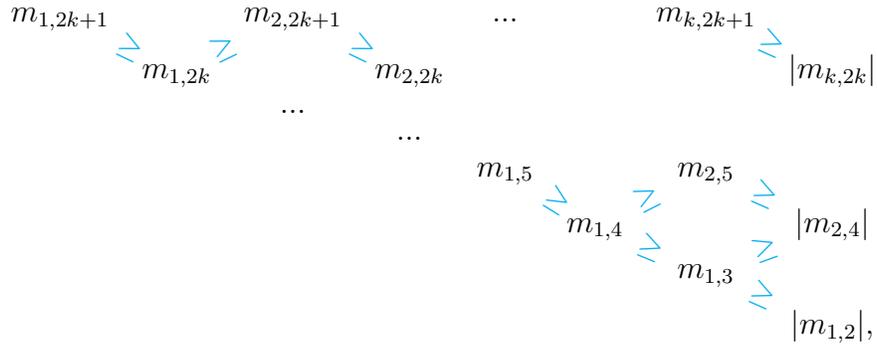
\[
\begin{tikzcd}[column sep=2pt,row sep=0pt]
m_{1, 2k+1} \arrow[rd, phantom, sloped, "\geq", cyan]  & & m_{2, 2k+1} \arrow[rd, phantom, sloped, "\geq", cyan]  &  & ... & &  m_{k,2k+1} \arrow[rd, phantom, sloped, "\geq", cyan]  \\
& m_{1, 2k} \arrow[ru, phantom, sloped, "\geq", cyan] & & m_{2,2k}   & & & & \lvert m_{k,2k} \rvert \\
& & ... \\
& & & ... \\
& & & & m_{1,5} \arrow[rd, phantom, sloped, "\geq", cyan]  & & m_{2,5} \arrow[rd, phantom, sloped, "\geq", cyan]  \\
& & & & & m_{1,4} \arrow[rd, phantom, sloped, "\geq", cyan] \arrow[ru, phantom, sloped, "\geq", cyan]   & & \lvert m_{2,4} \rvert \\
& & & & & & m_{1,3} \arrow[rd, phantom, sloped, "\geq", cyan] \arrow[ru, phantom, sloped, "\geq", cyan] \\
& & & & & & & \lvert m_{1,2} \rvert, 
\end{tikzcd}
\]
where the top row is omitted if $n+1 = 2k$ is even. 

It is convenient to assign coordinates to the entries in a Gelfand--Tsetlin pattern, counting the rows from the bottom, and the columns from the right of the diagram above (starting at $1$). Formally, we set
\[
\row(m_{i,j}) = j-1, \qquad \col(m_{i,j}) = j-2i +1. 
\]

\subsubsection{Classical modules} 
Given a regularly dominant weight $\boldsymbol\lambda = ([\lambda_1], \cdots, [\lambda_r])$, let $\GT(\boldsymbol\lambda)$ be the set of all Gelfand--Tsetlin patterns whose top row is $(m_{1,n+1}=\lambda_1, \cdots, \allowbreak  m_{r,n+1}=\lambda_r)$. 

%i-mult!
\Thm[{\cite{GavKl, GI97}}]
\label{thm: GIK B-action on GT}
The classical module $L_{\boldsymbol\lambda}$ has a $\C$-basis $\GT(\boldsymbol\lambda)$ such that 
\begin{align*} 
B_{2p} \cdot \sigma &= \sum^p_{j=1}P^j_{2p}(\sigma)\sigma^{2p}_{+j}
-\sum^p_{j=1}P^j_{2p}(\sigma^{2p}_{-j})\sigma^{2p}_{-j}, \\
B_{2p-1} \cdot \sigma &= \sum^{p-1}_{j=1}Q^j_{2p-1}(\sigma)
\sigma^{2p-1}_{+j} -\sum^{p-1}_{j=1}Q^j_{2p-1}({\sigma}^{2p-1}_{-j})
{\sigma}^{2p-1}_{-j} - R_{2p-1}(\sigma) \sigma, 
\end{align*} 
where the pattern $\sigma^k_{\pm j}$ is obtained from $\sigma$ by replacing
$m_{jk}$ by $m_{jk} \pm 1$, 
\eq \label{eq: l-coord} l_{j,2p}= m_{j,2p}+p-j, \qquad l_{j,2p-1}=m_{j,2p-1}+p-j, \eneq
and the scalar coefficients are given by 
\begin{align*}
P^j_{2p}(\sigma) &= {\rm i} \left( \frac{1}{\{ l'_j \} \{ l'_j + 1 \}}
\frac{\prod_{r=1}^p [l_r+l'_j][|l_r-l'_j-1|]
      \prod_{r=1}^{p-1}[l''_r +l'_j][|l''_r -l'_j-1|] }
     {\prod_{r\ne j}^p[l'_r+l'_j][|l'_r-l'_j|]
         [l'_r+l'_j+1][|l'_r-l'_j-1|]}\right)^{\frac12}, \\
%\label{A-coeff GT} \\ 
 Q^j_{2p-1}(\sigma) &= {\rm i} \left( \frac{1}{[\ell'_j]^2[2\ell'_j+1][2\ell'_j-1]}\frac{\prod_{r=1}^{p}
[\ell_r+\ell'_j][|\ell_r-\ell'_j|]
      \prod_{r=1}^{p-1}[\ell''_r +\ell'_j][|\ell''_r -\ell'_j|] }
     {\prod_{r\ne j}^{p-1} [\ell'_r+\ell'_j][|\ell'_r-\ell'_j|]
         [\ell'_r+\ell'_j-1] [|\ell'_r-\ell'_j-1|]} 
 \right)^{\frac12}, \\ 
R_{2p-1}(\sigma) &= \prod_{r=1}^{p} [l_{r,2p}] \prod_{r=1}^{p-1} [l_{r,2p-2}]
\left(\prod_{r=1}^{p-1} [l_{r,2p-1}] [l_{r,2p-1}-1]\right)^{-1}\ , 
\end{align*}
where we abbreviated $l_k =l_{k,2p+1}$, $l'_k = l_{k,2p}$ and $l''_k = l_{k,2p-1}$ in the formula for $P^j_{2p}(\sigma)$; as well as $\ell_k = l_{k,2p}$, $\ell'_k = l_{k,2p-1}$ and $\ell''_k = l_{k,2p-2}$ in the formula for $Q^j_{2p-1}(\sigma)$. 
\enthm

\Rem
In the formulae above, $\sigma^k_{\pm j}$ should be treated as zero whenever they violate the defining constraints of a Gelfand--Tsetlin pattern. \rmkend
\enrem

\Rem
One may worry that the formulae above are not well-defined because of the possibility of division by zero. However, it is easy to check that the Gelfand--Tsetlin constraints prevent most of the factors in the denominators from being equal to zero. The only exception is when $l_{p,2p-1} = 1$ in the formula for $R_{2p-1}(\sigma)$. In that case one sets $R_{2p-1}(\sigma) = 0$. \rmkend
\enrem 

\subsubsection{Non-classical modules}

If 
$\boldsymbol\lambda_+ = ([\lambda_1]_+, [\lambda_2]_+, \cdots, [\lambda_r]_+)$ is a regularly dominant weight, then each $\lambda_i \in \frac{1}{2} + \Z_{\geq 0}$. Let 
$\GTn(\boldsymbol\lambda_+)$ be the set of all Gelfand--Tsetlin patterns whose top row is 
$(m_{1,n+1}=\lambda_1, \cdots, \allowbreak  m_{r,n+1}=\lambda_r)$, satisfying the additional constraint that $m_{j,2j} \geq \frac{1}{2}$ for each $1 \leq j \leq r$. An element of $\GTn(\boldsymbol\lambda_+)$ may be presented graphically as 
\[
\begin{tikzcd}[column sep=2pt,row sep=0pt]
m_{1, 2k+1} \arrow[rd, phantom, sloped, "\geq", cyan]  & & m_{2, 2k+1} \arrow[rd, phantom, sloped, "\geq", cyan]  &  & ... & &  m_{k,2k+1} \arrow[rd, phantom, sloped, "\geq", cyan] & & \frac{1}{2}  \\
& m_{1, 2k} \arrow[ru, phantom, sloped, "\geq", cyan] & & m_{2,2k}   & & & &  m_{k,2k} \arrow[ru, phantom, sloped, "\geq", cyan] \\
& & ... \\
& & & ... \\
& & & & m_{1,5} \arrow[rd, phantom, sloped, "\geq", cyan]  & & m_{2,5} \arrow[rd, phantom, sloped, "\geq", cyan] & & \frac{1}{2}  \\
& & & & & m_{1,4} \arrow[rd, phantom, sloped, "\geq", cyan] \arrow[ru, phantom, sloped, "\geq", cyan]   & &  m_{2,4} \arrow[ru, phantom, sloped, "\geq", cyan] \arrow[rd, phantom, sloped, "\geq", cyan]\\
& & & & & & m_{1,3} \arrow[rd, phantom, sloped, "\geq", cyan] \arrow[ru, phantom, sloped, "\geq", cyan] & & \frac{1}{2}. \\
& & & & & & &  m_{1,2} \arrow[ru, phantom, sloped, "\geq", cyan]
\end{tikzcd}
\]

\Thm[{\cite{IK00, Ior-Kl}}] 
For each regularly dominant weight $\boldsymbol\lambda_+$ and a sequence $\vec{\varepsilon} = (\varepsilon_2, \cdots, \varepsilon_{n+1}) \in \{\pm1\}^{n}$, there exists a non-classical module $L_{\boldsymbol\lambda_+, {\vec{\varepsilon}}}$ with a $\C$-basis $\GTn(\boldsymbol\lambda_+)$ such that 
\begin{align*} 
B_{2p} \cdot \sigma &= \sum^p_{j=1}\dot{P}^j_{2p}(\sigma)\sigma^{2p}_{+j}
-\sum^p_{j=1}\dot{P}^j_{2p}(\sigma^{2p}_{-j})\sigma^{2p}_{-j} 
+ \delta_{m_{p,2p},\frac{1}{2}} \frac{\varepsilon_{2p+1}}{q^{\frac12} - q^{-\frac12}} \dot{S}_{2p}(\sigma) \sigma , \\
B_{2p-1} \cdot \sigma &= \sum^{p-1}_{j=1}\dot{Q}^j_{2p-1}(\sigma)
\sigma^{2p-1}_{+j} -\sum^{p-1}_{j=1}\dot{Q}^j_{2p-1}({\sigma}^{2p-1}_{-j})
{\sigma}^{2p-1}_{-j} + \varepsilon_{2p} \dot{R}_{2p-1}(\sigma) \sigma, 
\end{align*} 
where the scalar coefficients are given by 
\begin{align*}
\dot{P}^j_{2p}(\sigma) &= {\rm i} \left( 
\frac{(q-q\mi)^{-2}}{[l'_j] [ l'_j + 1 ]}
\frac{\prod_{r=1}^p [l_r+l'_j][|l_r-l'_j-1|]
      \prod_{r=1}^{p-1}[l''_r +l'_j][|l''_r -l'_j-1|] }
     {\prod_{r\ne j}^p[l'_r+l'_j][|l'_r-l'_j|]
         [l'_r+l'_j+1][|l'_r-l'_j-1|]}\right)^{\frac12}, \\
%\label{A-coeff GT nc} \\ 
 \dot{Q}^j_{2p-1}(\sigma) &=  \left(
 \frac{-1}{[\ell'_j]_+^2[2\ell'_j+1][2\ell'_j-1]}
  \frac{\prod_{r=1}^{p}
[\ell_r+\ell'_j][|\ell_r-\ell'_j|]
      \prod_{r=1}^{p-1}[\ell''_r +\ell'_j][|\ell''_r -\ell'_j|] }
     {\prod_{r\ne j}^{p-1}[\ell'_r+\ell'_j][|\ell'_r-\ell'_j|]
         [\ell'_r+\ell'_j-1] [|\ell'_r-\ell'_j-1|]} 
 \right)^{\frac12}, \label{B-coeff GT nc} \\ 
\dot{R}_{2p-1}(\sigma) &=  \prod_{r=1}^{p} [l_{r,2p}]_+ \prod_{r=1}^{p-1} [l_{r,2p-2}]_+
\left(\prod_{r=1}^{p-1} [l_{r,2p-1}]_+ [l_{r,2p-1}-1]_+ \right)^{-1}\ , \\
\dot{S}_{2p}(\sigma) &= {\rm i} \prod_{r=1}^{p} [l_{r,2p+1}-\tfrac12] \prod_{r=1}^{p-1} [l_{r,2p-1}-\tfrac12]
\left(\prod_{r=1}^{p-1} [l_{r,2p}+\tfrac{1}{2}] [l_{r,2p}-\tfrac12] \right)^{-1}, 
\end{align*}
where we abbreviated $l_k =l_{k,2p+1}$, $l'_k = l_{k,2p}$ and $l''_k = l_{k,2p-1}$ in the formula for $\dot{P}^j_{2p}(\sigma)$; as well as $\ell_k = l_{k,2p}$, $\ell'_k = l_{k,2p-1}$ and $\ell''_k = l_{k,2p-2}$ in the formula for $\dot{Q}^j_{2p-1}(\sigma)$. 
\enthm

In the formulae above, $\sigma^k_{\pm j}$ should be treated as zero whenever they are not in $\GTn(\boldsymbol\lambda_+)$.

%%%%%%%%%%%%%%
%%%%%%%%%%%%%%
%%%%%%%%%%%%%%

\section{Evaluation homomorphism} 
\label{sec: eval homo}

\nc{\LHS}{\on{LHS}}

\subsection{Construction via braid group action} 

In the theorem below, we define evaluation homomorphisms between affine and finite split quantum symmetric pair coideal subalgebras of type $\mathsf{AI}$. Formally, the definition is similar to the one for evaluation homomorphisms for quantum affine algebras \cite[Proposition 12.2.10]{chari-pressley}. The main difference is that a different braid group action is involved (i.e., not the usual Lusztig action). 

We also remark that an evaluation homomorphism for twisted $q$-Yangians, which are coideals in $U_q(\widehat{\mathfrak{gl}}_{n+1})$, was defined in \cite{MRS-03} using the RTT approach. 
It can be shown that our evaluation map coincides with a restriction of the evaluation map in \emph{op.\ cit.} to $\aui$, using the isomorphism from \cite[Theorem 11.7]{kolb-14}, and a similar argument to that in \cite{FrMukHopf}. 
However, we will not need this fact. 

\Thm \label{pro: ev homo}
For any $a \in \C^\times$, there exists a $\mathcal{B}$-equivariant algebra homomorphism 
\[
\eva \colon \aui \to \Ui, 
\]
given by 
\begin{alignat}{3}
B_i &\mapsto B_i, \qquad && \KK_i &\ \mapsto \ & \KK_i \qquad \qquad \qquad   (1 \leq i \leq n), \\ 
B_0 &\mapsto a \Tbr_{\pi\mi\omega_1}(B_1), \qquad && \KK_0 &\ \mapsto \ & a^2 \Tbr_{\pi\mi\omega_1}(\KK_1). 
\end{alignat} 
More explicitly, 
\begin{align} 
\eva(B_0) &= a \Tbr_{n} \cdots \Tbr_1 (B_1) = a \KK_1\mi \cdots \KK_n\mi [B_1, [B_{2}, \cdots, [B_{n-1}, B_n]_q \cdots ]_q ]_q, \\
\eva(\KK_0) &= a^2 \Tbr_{n} \cdots \Tbr_1 (\KK_1) = a^2 \KK_1\mi \cdots \KK_n\mi, \quad \eva(\CCC) = a^2.  
\end{align}
\enthm

\Proof
To prove that $\eva$ is a homomorphism, we only need to check that the relations involving $B_0$ and $\KK_0$ are preserved. \emph{Firstly}, consider the relation 
\eq \label{eq: ev rel 1}
\LHS_1 = B_i^2 B_0 - [2] B_iB_0B_i + B_0B_i^2 + q^{-1} \KK_iB_0 = 0, \qquad i \in \{1, n\}.
\eneq
In both cases the proof is the same, so let us take $i = 1$. Since $B_1$ commutes with $B_j$ for $3 \leq j \leq n$, it suffices to consider the $\mathsf{A}_2$ case. 
We have 
\begin{align}
(a\KK_1 \KK_2)\mi \eva (\LHS_1) =& \ B_1^2[B_1, B_2]_q - [2]B_1 [B_1, B_2]_q B_1 + [B_1, B_2]_q B_1^2 + q\mi \KK_1[B_1,B_2]_q \\ 
=& \  [B_1, (B_1^2B_2 - [2] B_1B_2B_1 + B_2 B_1^2 + q\mi \KK_1 B_2)]_q = 0, 
\end{align} 
since the second term inside the $q$-commutator vanishes by \eqref{eq: QSP Serre}. 

\emph{Secondly}, we check that the relation 
\[
\LHS_2 = B_0^2 B_j - [2] B_0B_jB_0 + B_jB_0^2 + q^{-1} \KK_0B_j = 0, \qquad j \in \{1, n\}
\]
is preserved. 
Again, we may assume, without loss of generality, that $i = n$. 
It is easy to calculate that $\Tbr_{1}\mi \cdots \Tbr_n\mi (B_n) = a^{-1} \eva(B_0)$, which implies that 
\[ 
\Tbr_{1}\mi \cdots \Tbr_n\mi (\eva(\LHS_2)) = a \eva(\LHS_1), 
\]
with $i=1$.  
Since $\Tbr_j$ are algebra isomorphisms, $\eva(\LHS_2) = 0$ follows from the fact that relation \eqref{eq: ev rel 1} is preserved by $\eva$. 

\emph{Thirdly}, we show that the relation 
\[
\LHS_3 = [B_0, B_j] = 0, \qquad j \notin \{1, n\}
\]
is also preserved. Since $B_j$ commutes with all $B_i$'s except $B_{j \pm 1}$, the calculation reduces to the $\mathsf{A}_3$ case (with $j = 2$). We have 
\begin{align}
(a \KK_1 \KK_2 \KK_3)\mi [2] \eva(\LHS_3) =& \ (B_2^2 B_1 - [2]B_2B_1B_2) B_3 + q^{2}(B_2^2 B_3 - [2]B_2B_3B_2) B_1 \\  \label{eq: ev proof} 
& \ + B_1 ([2]B_2B_3B_2 - B_3 B_2^2) + q^{2} B_3 ([2]B_2B_1B_2 - B_1 B_2^2). 
\end{align} 
The first two terms on the RHS above are equal to: 
\begin{align}
(B_2^2 B_1 - [2]B_2B_1B_2) B_3 =& \ - B_1(B_2^2B_3 + q^{-1} \KK_2 B_3), \\ 
q^{2}(B_2^2 B_3 - [2]B_2B_3B_2) B_1 =& \ -q^{2}B_3(B_2^2B_1 + q^{-1} \KK_2B_1). 
\end{align}
Adding the results to the first and second terms in \eqref{eq: ev proof}, respectively, we see that both sums vanish by \eqref{eq: QSP Serre}. 

\emph{Finally}, we prove the $\mathcal{B}$-equivariance. Again, it suffices to check that $\eva$ preserves the $\mathcal{B}$-action on $B_0$ and $\KK_0$. We will verify the $B_0$ case explicitly, the $\KK_0$ case being analogous. First take $2 \leq i \leq n-1$. Then $\Tbr_i(B_0) = B_0$. On the other hand, 
\[
\Tbr_i (\eva(B_0)) = a \Tbr_i\Tbr_{\pi\mi\omega_1}(B_1) = a \Tbr_{\pi\mi\omega_1} \Tbr_{i+1}(B_1) 
= a \Tbr_{\pi\mi\omega_1} (B_1) = \eva(B_0). 
\] 
Hence $\eva(\Tbr_i(B_0)) = \Tbr_i(\eva(B_0))$. Next, take $i=1$. Then 
\begin{align} 
\Tbr_1 (\eva(B_0)) &= a \Tbr_1\Tbr_n \cdots \Tbr_{1}(B_1) = a \Tbr_n \cdots \Tbr_{1} \Tbr_{2}(B_1) 
= a \Tbr_n \cdots \Tbr_{2} (B_2) \\ 
&= a \KK_2\mi \cdots \KK_n\mi [B_2, [B_{3}, \cdots, [B_{n-1}, B_n]_q \cdots ]_q ]_q
\end{align} 
On the other hand, 
\begin{align}
\eva(\Tbr_1(B_0)) = [\eva(B_0), B_1]_q &= a [\Tbr_n \cdots \Tbr_{1}(B_1), B_1]_q \\
&= a \KK_1\mi \cdots \KK_n\mi [ \cdots [ [[B_1, B_{2}]_q, B_1]_q, B_3]_q \cdots, B_n]_q \\
&= a \KK_2\mi \cdots \KK_n\mi [B_2, [B_{3}, \cdots, [B_{n-1}, B_n]_q \cdots ]_q ]_q, 
\end{align} 
since $[[B_1, B_{2}]_q, B_1]_q = \KK_1 B_2$ by \eqref{eq: QSP Serre}. The proof for $i = n$ is similar. 
\enproof

\Rem
Let $\tilde{\sigma}$ be the automorphism of $\aui$ defined by 
\[
\tilde{\sigma}(A_{i,r}) = A_{n+1-i,r}, \quad \tilde{\sigma}(H_{i,r}) = H_{n+1-i,r}, \quad \tilde{\sigma}(\KK_i) = \KK_{n+1-i},
\]
and $\sigma$ its restriction to $\Ui$. Then $\on{ev}^a = \sigma \circ \eva \circ \tilde{\sigma}$ is another homomorphism 
$\aui \to \Ui$, explicitly given by 
\begin{alignat}{3} 
B_i &\mapsto B_i, \qquad && \KK_i &\ \mapsto \ & \KK_i \qquad \qquad \qquad   (1 \leq i \leq n), \\ 
B_0 &\mapsto a \Tbr_{\pi\omega_n}(B_n), \qquad && \KK_0 &\ \mapsto \ & a^2 \Tbr_{\pi\omega_n}(\KK_n). 
\end{alignat} 
\enrem 

Going forward, we will also use the same notation $\eva$ for the composition of the map from Theorem \ref{pro: ev homo}  with the central reduction $\Ui \to \ui$ from \S \ref{subsec: parameters}.
Given a $\ui$-module $M$, we will denote its pullback via $\eva$ as 
\[
M(a) := \eva^*(M). 
\]

%%%%%%%%%%%%%%%%%%%%

\subsection{Compatibilities} 

We next show that the evaluation homomorphisms are compatible with a family of natural inclusions between coideal subalgebras of increasing ranks. Consider the map 
\[
\iota \colon \auiv{n} \to \auiv{n+1}, \qquad B_i \mapsto B_i, \quad \KK_i \mapsto \KK_i, \quad B_0 \mapsto [B_n, B_0]_q, \quad \KK_0 \mapsto \KK_0 \KK_n,
\]
for $1 \leq i \leq n-1$. 

\Lem
The map $\iota$ induces an injective algebra homomorphism, which fits into the following commutative diagram 
\[
\begin{tikzcd}
\auiv{n-1} \arrow[r, "\iota", hookrightarrow] \arrow[d, "\eva", swap]  & \auiv{n} \arrow[d, "\eva"] \\
\Uiv{n-1} \arrow[r, "\iota", hookrightarrow] & \Uiv{n}. 
\end{tikzcd}
\] 
In terms of the Lu--Wang presentation, it is given by 
\[
A_{i,r} \mapsto A_{i,r}, \quad H_{i,r} \mapsto H_{i,r}, \quad \KK_i^{\pm1} \mapsto \KK_i^{\pm1}, \quad \Ck \mapsto \Ck,
\]
for $1 \leq i \leq n-1$. 
\enlem

\Proof
The proof is an easy relation check and diagram chase. 
\enproof

\subsection{Explicit formulae for some generators} 

In the lemma below, we give explicit formulae for the images of the generators $A_{i,1}$ ($i \in \indx_0$) under the evaluation homomorphism.  
Set 
\[
P(y_1, \cdots, y_k) = [y_1, [y_2, [y_3, \cdots [y_{k-1}, y_k]_q \cdots ]_q ]_q. 
\]

\Lem \label{lem: exp formula ev A}
Let $i \in \indx_0$. Then 
\begin{align}
\eva(A_{i,1}) &= o(i) a \Tbr_i\mi \cdots \Tbr_1\mi \Tbr_i \cdots \Tbr_1(B_1) \\
&= o(i) a \KK_1\mi \cdots \KK_{i-1}\mi P(P(B_i, \cdots, B_1), B_1, \cdots, B_{i-1}). 
\end{align}
In particular, 
\[
\eva(A_{i+1,1}) = -\Tbr_{i+1}\mi \Tbr_i\mi \Tbr_{i+1}\Tbr_i \eva(A_{i,1}). 
\]
\enlem

\Proof
By definition, we have $A_{i,1} = o(i) \Tbr_{\omega_i}\mi(B_i)$. 
By \cite[Proposition 9.2]{Przez-23}, 
\[
\omega_i = \pi^i (s_{n-i+1}\cdots s_n) \cdots(s_1\cdots s_i).
\]
Accordingly, 
\[
A_{i,1} = o(i) \Tbr_i\mi \cdots \Tbr_1\mi \Tbr_{i+1}\mi\cdots\Tbr_n\mi(B_0). 
\]
Therefore, 
\begin{align}
o(i) \eva(A_{i,1}) &= \Tbr_i\mi \cdots \Tbr_1\mi \Tbr_{i+1}\mi\cdots\Tbr_n\mi(\eva(B_0)) \\
&= \Tbr_i\mi \cdots \Tbr_1\mi \Tbr_{i+1}\mi\cdots\Tbr_n\mi(\eva(B_0)) \\
&= a \Tbr_i\mi \cdots \Tbr_1\mi \Tbr_i \cdots \Tbr_1(B_1) \\
&= a \KK_i \Tbr_i\mi \cdots \Tbr_1\mi P(B_1, \cdots, B_i) \\ 
&= a\KK_1\mi \cdots \KK_{i-1}\mi P( P(B_i, \cdots, B_1) , B_1, \cdots, B_{i-1}) 
\end{align}
where the first equality follows from $\mathcal{B}$-equivariance (Theorem \ref{pro: ev homo}). 
The second statement of the lemma follows immediately from the first. 
\enproof

%%%%%%%%%%%%%%
%%%%%%%%%%%%%%
%%%%%%%%%%%%%%

\section{Action on the GT basis} 
\label{sec: GT action}

In this section we compute the action of the Lu--Wang generators on the Gelfand--Tsetlin basis. 
In subsections \ref{subsec: hor ind}--\ref{subsec: Theta on GT}, we treat the classical module case. 
Let us fix a regularly dominant weight $\boldsymbol\lambda$, and work with the evaluation $L_{\boldsymbol\lambda}(a) = \eva^*(L_{\boldsymbol\lambda})$ of the corresponding classical simple module for~$\ui$. 
The following is the first main result of this section. It will be proven in \S \ref{subsec: hor ind}--\ref{subsec: ver ind}. 

\Thm 
\label{thm: A-action on GT}
Let $\sigma \in \GT(\boldsymbol\lambda)$. Then 
\begin{align} 
A_{2p,s} \cdot \sigma &= \sum^p_{j=1} a^sq^{s \cdot \zeta_{j,2p}} P^j_{2p}(\sigma)\sigma^{2p}_{+j}
-\sum^p_{j=1} a^sq^{s \cdot (2-\zeta_{j,2p})} P^j_{2p}(\sigma^{2p}_{-j})\sigma^{2p}_{-j}, \\
A_{2p-1,s} \cdot \sigma &= \sum^{p-1}_{j=1} a^sq^{s \cdot \zeta_{j,2p-1}} Q^j_{2p-1}(\sigma) \sigma^{2p-1}_{+j}  - \sum^{p-1}_{j=1} a^s q^{s \cdot (2-\zeta_{j,2p-1})} Q^j_{2p-1}({\sigma}^{2p-1}_{-j}) {\sigma}^{2p-1}_{-j} \\
&- a^s R_{2p-1}(\sigma) \sigma, 
\end{align} 
for each $s \in \Z$, where 
\[
\zeta_{j,k} = \col(m_{j,k}) + 2m_{j,k} = 
\left\{ \begin{array}{l l}
2l_{j,2p} + 1 & \mbox{if } k=2p, \\[2pt]
2l_{j,2p-1} & \mbox{if } k = 2p-1.  
\end{array} \right. 
\]
\enthm

The proof is by induction on $k$ (`horizontal' induction) and $s$ (`vertical induction') in $A_{k,s}$, and consists of three parts. In the first two parts, we establish the theorem for $s=-1$, separately handling the cases where $k=2p$ and  $k=2p-1$. In the third part we generalize the result to arbitrary $s$. 

Assuming Theorem \ref{thm: A-action on GT}, let us introduce some notation to be used in the proof. We will write $A_{k,s}$ as the sum of operators 
\begin{align*}
A_{2p,s} &=  \sum_{j=1}^p  A_{2p,s}^{+j} - \sum_{j=1}^p  A_{2p,s}^{-j}, \qquad
A_{2p-1,s} = \sum_{j=1}^{p-1}  A_{2p-1,s}^{+j} - \sum_{j=1}^{p-1}  A_{2p-1,s}^{-j} -  A_{2p-1,s}^{(0)}.  
\end{align*}

The following is the second main result of this section. It will be proven in \S \ref{subsec: Theta on GT}. 

\Thm 
\label{thm: Theta-action on GT}
Let $\sigma \in \GT(\boldsymbol\lambda)$. Then $\sigma$ is an eigenvector for the action of the commutative algebra $\mathcal{H}$. The corresponding eigenvalues are given by 
\begin{align*}
\thvar{2p} 
\cdot \sigma &= 
\prod_{j=1}^p \frac{(1\smin q^{2l_{j,2p+1}\smin 1}az)(1\smin q^{1\smin 2l_{j,2p+1}}az)}{(1\smin q^{2l_{j,2p}+1}az)(1\smin q^{\smin 1\smin 2l_{j,2p}}az)} \frac{(1\smin q^{2l_{j,2p\smin 1}\smin 1}az)(1\smin q^{1\smin 2l_{j,2p\smin 1}}az)}{(1\smin q^{2l_{j,2p}\smin 1}az)(1\smin q^{1\smin 2l_{j,2p}}az)} \sigma, \\
\thvar{2p\smin 1} 
\cdot \sigma &=\prod_{j=1}^{p-1} 
\frac{(1-q^{2l_{j,2p}}az) (1-q^{-2l_{j,2p}}az)}{(1-q^{2l_{j,2p-1}}az)(1-q^{-2l_{j,2p-1}}az)} 
\frac{(1-q^{2l_{j,2p-2}}az)(1-q^{-2l_{j,2p-2}}az)}{(1-q^{2l_{j,2p-1}-2}az)(1-q^{2-2l_{j,2p-1}}az)} \cdot  \\
&\cdot \frac{(1-q^{2l_{p,2p}}az) (1-q^{-2l_{p,2p}}az)}{(1-az)^2} \sigma, 
\end{align*}
where $l_{p,2p-1} = 0$ by convention. 
\enthm

\Rem \label{rem: nice formula Theta}
We can express the formulae above in the following uniform way, using $\zeta$-coordinates: 
\[
\thvar{k} 
\cdot \sigma = 
\prod_{j=1}^{\lceil k/2 \rceil} 
\frac{(1-q^{\zeta_{j,k+1}-1}az)(1-q^{1-\zeta_{j,k+1}}az)}{(1-q^{\zeta_{j,k}}az)(1-q^{-\zeta_{j,k}}az)}
\frac{(1-q^{\zeta_{j,k-1}-1}az)(1-q^{1-\zeta_{j,k-1}}az)}{(1-q^{\zeta_{j,k}-2}az)(1-q^{2-\zeta_{j,k}}az)} \sigma, 
\]
where, for each $p > 0$, by convention, we set $\zeta_{p,2p-1}=0$ and $\zeta_{p,2p-2} = -1$. Note that when $k$ is odd, one of the fractions on the RHS reduces to $1$. \rmkend
\enrem

The analogues of Theorems \ref{thm: A-action on GT} and \ref{thm: Theta-action on GT} for non-classical modules will be presented in \S \ref{subsec: nc Theta on GT}. 

\subsection{Auxiliary formulae} 
\nc{\trip}{\mathbf{\underline{l}} } 

Below we collect several auxiliary technical results which will be used in the proofs of Theorems \ref{thm: A-action on GT} and \ref{thm: Theta-action on GT}. First of all, we will use the following easy consequence of the Lu--Wang relations throughout.  

\Lem \label{lem: Serrcor}
For any $1 \leq j \leq n-1$, we have  
\begin{equation} \label{eq: Serrecor}
 A_{j ,\smin1} = q\mi[B_{j \smin1}, [A_{j \smin1,\smin1}, B_{j }]_{q} ]_q. 
\end{equation} 
if we specialize $\KK_{j-1} \mapsto - q$. 
\enlem 

\Proof
Observe that 
\begin{align*}
-q\mi \KK_{j-1} A_{j ,\smin1} &= (B_{j \smin1}^2A_{j ,\smin1} - [2]B_{j \smin1}A_{j ,\smin1}B_{j \smin1} + A_{j ,\smin1}B_{j \smin 1}^2) \\ 
&= [B_{j \smin1}, [B_{j \smin1}, A_{j ,\smin1}]_{q\mi} ]_q 
= q\mi[B_{j \smin1}, [A_{j \smin1,\smin1}, B_{j}]_{q} ]_q, 
\end{align*}
where the first equality follows from the Serre relation \eqref{eq: grel6}, and the third one from \eqref{eq: grel4}. 
Specializing $\KK_{j-1} \mapsto - q$ recovers the desired formula. 
\enproof

We will also need the following technical result from \cite[Proposition A]{GI97}. 
Let us first introduce some notation. 
Suppose that we are given a sequence of triples $\trip = (l_r,l'_r,l''_r)_{r=1}^p$ which satisfy the Gelfand--Tsetlin constraints (via the translation between $m$- and $l$-coordinates given in \eqref{eq: l-coord}). We allow the final triple to be `truncated', e.g., $(l_p,l'_p,\cdot)$ is allowed. Given $1 \leq i \leq p$, let $\trip_{\smin i}$ be the sequence obtained from $\trip$ by replacing $l'_i$ with $l'_i - 1$.
Moreover, define the following functions: 
\begin{align*}
\varphi_r&(\trip) 
= \frac{\prod_{s=1}^p f(l_{s};l'_{r})
\prod_{s=1}^{p-1} f(l''_{s};l'_{r})}
{\prod_{s\ne r}^p f(l'_{s};l'_{r})f(l'_{s}+1;l'_{r})}, \qquad
\Phi(\trip) =  \sum_{r=1}^{p} \frac{\varphi_r(\trip_{\smin i}) - \varphi_r(\trip)}{[2l_{r,2p}]}, 
\end{align*}
where $f(x;y) = [x+y][x-y-1]$. The function $f(x;y)$ has the property 
\eq
\label{eq: f property} 
f(x;z) - f(y;z) = f(x;y-1). 
\eneq

\Lem \label{pro: gen Appendix pro}
One has $\Phi(\trip) = 1$. 
\enlem

\Proof 
See \cite[Appendix]{GI97}. 
\enproof

We will also need another technical lemma of similar flavour.  
Define   
\begin{align*} 
\Psi(\trip) &= 1 - \sum_{i=1}^p \frac{za(q-q\mi)^2}{[2l'_i]} \left( \psi_i(\trip) - \psi_i(\trip_{\smin i}) 
\right), \\ 
\psi_i (\trip) &= 
\frac{[l'_i][l'_i + 1]}{(1-q^{2l'_i{+}1}az)(1-q^{\smin2l'_i\smin1}az)}
\frac{\prod_{r=1}^p [l_r+l'_i][|l_r-l'_i-1|] \prod_{r=1}^{p-1}[l''_r +l'_i][|l''_r -l'_i-1|] }
{\prod_{r\ne i}^p[l'_r+l'_i][|l'_r-l'_i|][l'_r+l'_i+1][|l'_r-l'_i-1|]}. 
\end{align*} 
In the cases where a triple is truncated, we remove all factors depending on the missing $l$-coordinate in the formula for 
$\psi_i (\trip)$. 

\Lem \label{pro: gen Appendix pro2}
We have 
\[
\Psi(\trip) = 
\prod_{i=1}^p \frac{(1-q^{2l_{i}-1}az)(1-q^{1-2l_{i}}az)}{(1-q^{2l'_{i}+1}az)(1-q^{-1-2l'_{i}}az)} \frac{(1-q^{2l''_{i}-1}az)(1-q^{1-2l''_{i}}az)}{(1-q^{2l'_{i}-1}az)(1-q^{1-2l'_{i}}az)}, 
\]
where, on the RHS, $l_{p,2p-1} = 0$ by convention. 
\enlem 

\Proof
The proof is by induction on the length $p$ of a sequence of triples. The length~$1$ case can be checked by a direct calculation. Then the inductive step is split into two parts. Let us first introduce some notation for various sequences of triples obtained from $\trip$. Let $\trip^{\bullet} = (l_r,l'_r,l''_r)_{r=2}^p$; and let $\trip^{\circ}$ be the sequence obtained from $\trip$ by replacing $l_1$ with $l'_1+1$. In the first part of the inductive step, one shows that 
\eq 
\label{eq: first step ind}
\Psi(\trip^{\circ}) = \frac{(1-q^{2l''_{1}-1}az)(1-q^{1-2l''_{1}}az)}{(1-q^{2l'_{1}-1}az)(1-q^{1-2l'_{1}}az)} \Psi(\trip^{\bullet}). 
\eneq
Then, using \eqref{eq: first step ind}, one shows in the second step that 
\[
\Psi(\trip) = 
\frac{(1-q^{2l_{1}-1}az)(1-q^{1-2l_{1}}az)}{(1-q^{2l'_{1}+1}az)(1-q^{-1-2l'_{1}}az)} \frac{(1-q^{2l''_{1}-1}az)(1-q^{1-2l''_{1}}az)}{(1-q^{2l'_{1}-1}az)(1-q^{1-2l'_{1}}az)} 
\Psi(\trip^{\bullet}). 
\] 
The corresponding proofs are similar, so we will only prove the second step explicitly. 

It follows from \eqref{eq: f property} that 
\[
\Psi(\trip) = \Psi(\trip^\circ) + [l_1+l'_1][l_1-l'_1-1] (\Psi(\trip^\sharp)-1), \qquad
\Psi(\trip^\sharp)-1 = \frac{1}{[2l'_1]} \left( \Psi(\trip^\circ) - \Psi(\trip^*)  \right), 
\]
where $\trip^\sharp$ is the sequence obtained from $\trip$ by replacing the first triple with $(\cdot,l'_1,l''_2)$; and $\trip^*$ by replacing it with $(l'_1,l'_1,l''_2)$. Arguing as in \cite[(A4)]{GI97}, and using \eqref{eq: first step ind}, we get 
\[
\Psi(\trip^*) =  \Psi(\trip^\circ)|_{l'_1 \to l'_1 + 1} = \frac{(1-q^{2l''_{1}-1}az)(1-q^{1-2l''_{1}}az)}{(1-q^{2l'_{1}+1}az)(1-q^{-1-2l'_{1}}az)} \Psi(\trip^{\bullet}). 
\]
Therefore, 
\begin{align*}
\Psi(\trip) &= \Bigg(
\frac{(1-q^{2l''_{1}-1}az)(1-q^{1-2l''_{1}}az)}{(1-q^{2l'_{1}-1}az)(1-q^{1-2l'_{1}}az)} + \frac{[l_1+l'_1][l_1-l'_1-1]}{[2l'_1]} \cdot \\
&\cdot \left( \frac{(1-q^{2l''_{1}-1}az)(1-q^{1-2l''_{1}}az)}{(1-q^{2l'_{1}-1}az)(1-q^{1-2l'_{1}}az)}  -   \frac{(1-q^{2l''_{1}-1}az)(1-q^{1-2l''_{1}}az)}{(1-q^{2l'_{1}+1}az)(1-q^{-1-2l'_{1}}az)} 
\right)
\Bigg) \Psi(\trip^{\bullet}). 
\end{align*}
An easy calculation shows that the expression in the outer brackets evaluates to 
\[
\frac{(1-q^{2l_{1}-1}az)(1-q^{1-2l_{1}}az)}{(1-q^{2l'_{1}+1}az)(1-q^{-1-2l'_{1}}az)} \frac{(1-q^{2l''_{1}-1}az)(1-q^{1-2l''_{1}}az)}{(1-q^{2l'_{1}-1}az)(1-q^{1-2l'_{1}}az)}, 
\]
completing the proof. 
\enproof

\subsection{Horizontal induction: $A_{2p,-1}$} 
\label{subsec: hor ind}

In this subsection, set $l_k =l_{k,2p+1}$, $l'_k = l_{k,2p}$ and $l''_k = l_{k,2p-1}$; and $\ell_k = l_{k,2p}$, $\ell'_k = l_{k,2p-1}$ and $\ell''_k = l_{k,2p-2}$. 
It follows from \eqref{lem: Serrcor} and induction that 
\begin{align}
 A_{2p ,\smin1} \cdot \sigma &= \sum_{j,i \neq h} \tau_{\star j,\star i,\star h} \sigma_{\star j,\star i,\star h}^{2p,2p-1,2p-1}
 + \sum_{j,i} \tau_{\star j,\pm i,\pm i} \sigma_{\star j,\pm i,\pm i}^{2p,2p-1,2p-1} \\ 
&+ \sum_{j,i} \tau_{\star j,\star i} \sigma_{\star j,\star i}^{2p,2p-1} 
+ \sum_{j} \tau_{\pm j} \sigma_{\pm j}^{2p}, 
\end{align}
where $\star \in \{ +, - \}$ is used to indicate that we do note assume any relation between the respective signs on $j,i$ and $h$. The multiple sub/superscript notation on $\sigma$ means that, e.g., $\sigma_{+ j, - i ,+ h}^{2p,2p-1,2p-1}$ is obtained from $\sigma$ by replacing $m_{j,2p}$ by $m_{j,2p}+1$, $m_{i,2p-1}$ by $m_{i,2p-1}-1$, and $m_{h,2p-1}$ by $m_{h,2p-1}+1$. 

We claim that all the coefficients $\tau$ vanish except for $\tau_{\pm j}$. The proof for coefficients belonging to each of the summations above is formally the same, regardless of the choice of signs, so we will choose, without loss of generality, positive signs throughout. 

Observe that the final Gelfand--Tsetlin pattern, e.g., $\sigma_{+ j, - i ,+ h}^{2p,2p-1,2p-1}$, can be obtained from $\sigma$ by changing the entries in different orders. Some of these may produce an intermediate result which is not a Gelfand--Tsetlin pattern. Whenever this does \emph{not} happen, while applying the three operators on the RHS of \eqref{lem: Serrcor}, in any of the four possible orders, we say that we are dealing with a \emph{generic case}. 
Below we restrict ourselves to providing the required calculations in the generic cases. The non-generic cases can be obtained through a straightforward modification.

\subsubsection{Case: $\tau_{+ j, + i, + h}$ and $i \neq h$} 
\label{subsec: jih coeff} 

The total coefficient $\tau_{+j,+i,+h}$ on $\sigma_{+j,+i,+h}^{2p,2p-1,2p-1}$ is equal to the sum of the coefficients $\tau', \tau''$ on 
\[ 
q\mi[B_{2p \smin1}^{+h}, [A_{2p \smin1,\smin1}^{+i}, B_{2p}^{+j}]_{q} ]_q \cdot \sigma, \qquad 
q\mi[B_{2p \smin1}^{+i}, [A_{2p \smin1,\smin1}^{+h}, B_{2p}^{+j}]_{q} ]_q \cdot \sigma, 
\]
respectively. 
Define  
\begin{align*}
F(j,i,h) &= a\mi  \big( P^j_{2p}(\sigma) Q_{2p\smin1}^i(\sigma^{2p}_{+j}) Q_{2p\smin1}^h(\sigma^{2p,2p\smin1}_{+j,+i}) 
-  qP^j_{2p}(\sigma_{+i}^{2p\smin1}) Q_{2p\smin1}^i(\sigma) Q_{2p\smin1}^h(\sigma^{2p,2p\smin1}_{+j,+i}) \\
&- q P^j_{2p}(\sigma_{+h}^{2p\smin1}) Q_{2p\smin1}^i(\sigma_{+j,+h}^{2p,2p\smin1}) Q_{2p\smin1}^h(\sigma) 
+ q^2 P^j_{2p}(\sigma_{+i,+h}^{2p\smin1,2p\smin1}) Q_{2p\smin1}^i(\sigma_{+h}^{2p\smin1}) Q_{2p\smin1}^h(\sigma) \big).
\end{align*}
Then $\tau' = q^{\gamma_i}F(j,i,h)$ and $\tau'' = q^{\gamma_h}F(j,h,i)$, where 
\[ \gamma_i = -(\col(m_{i,2p-1})+2m_{i,2p-1}+1). \]

Abbreviate 
\[
d_j = l_{j,2p} \quad
d_i = l_{i,2p-1} \quad 
d_h = l_{h,2p-1}, \quad 
\langle \cdot \rangle = [ \cdot ]^{\frac12}. 
\]
Let 
\begin{align}
X_1 &= - {\rm i} a\mi  \Bigg( \frac{1}{\{ l'_j \} \{ l'_j + 1 \}}
\frac{\prod_{r=1}^p [l_r+l'_j][|l_r-l'_j-1|]
      \prod_{r\neq i,h}^{p-1}[l''_r +l'_j][|l''_r -l'_j-1|] }
     {\prod_{r\ne j}^p[l'_r+l'_j][|l'_r-l'_j|]
         [l'_r+l'_j+1][|l'_r-l'_j-1|]} \\ 
&\cdot 
 \frac{1}{[\ell'_i]^2[2\ell'_i+1][2\ell'_i-1]}\frac{\prod_{r\neq j}^{p}
[\ell_r+\ell'_i][|\ell_r-\ell'_i|]
      \prod_{r=1}^{p-1}[\ell''_r +\ell'_i][|\ell''_r -\ell'_i|] }
     {\prod_{r\ne i}^{p-1}[\ell'_r+\ell'_i][\ell'_r+\ell'_i-1]
     \prod_{r \neq i,h}^{p-1} [|\ell'_r-\ell'_i|][|\ell'_r-\ell'_i-1|]} \\
&\cdot 
 \frac{1}{[\ell'_h]^2[2\ell'_h+1][2\ell'_h-1]}\frac{\prod_{r\neq j}^{p}
[\ell_r+\ell'_h][|\ell_r-\ell'_h|]
      \prod_{r=1}^{p-1}[\ell''_r +\ell'_h][|\ell''_r -\ell'_h|] }
     {\prod_{r\ne h}^{p-1}[\ell'_r+\ell'_h][\ell'_r+\ell'_h-1]
     \prod_{r \neq i,h}^{p-1} [|\ell'_r-\ell'_h|][|\ell'_r-\ell'_h-1|]}\Bigg)^\frac12. 
\end{align} 
Then $X_1\mi F(j,i,h)$ is equal to 
\begin{align}
&\bigg( \langle d_i + d_j\rangle  \langle |d_i - d_j -1|\rangle  \langle d_h + d_j\rangle  \langle |d_h - d_j - 1|\rangle  \\
&\ \ \ \cdot \langle d_j + d_i + 1\rangle  \langle |d_j - d_i +1|\rangle  \langle d_j + d_h + 1\rangle  \langle |d_j - d_h +1|\rangle  \\
&- q \langle d_i + d_j + 1\rangle  \langle |d_i - d_j |\rangle  \langle d_h + d_j\rangle  \langle |d_h - d_j - 1|\rangle  \\
&\ \ \ \cdot \langle d_j + d_i \rangle  \langle |d_j - d_i|\rangle  \langle d_j + d_h + 1\rangle  \langle |d_j - d_h +1|\rangle \bigg) 
\bigg( [|d_h - d_i|][|d_i - d_h + 1|] \bigg)\mi \\ 
&+ \bigg(- q \langle d_i + d_j\rangle  \langle |d_i - d_j -1|\rangle  \langle d_h + d_j+1\rangle  \langle |d_h - d_j|\rangle  \\
&\ \ \ \cdot \langle d_j + d_i + 1\rangle  \langle |d_j - d_i +1|\rangle  \langle d_j + d_h\rangle  \langle |d_j - d_h|\rangle  \\
&+ q^2 \langle d_i + d_j +1\rangle  \langle |d_i - d_j|\rangle  \langle d_h + d_j+1\rangle  \langle |d_h - d_j|\rangle  \\
&\ \ \ \cdot \langle d_j + d_i\rangle  \langle |d_j - d_i|\rangle  \langle d_j + d_h\rangle  \langle |d_j - d_h|\rangle \bigg) \bigg(
[|d_h - d_i|][|d_i - d_h - 1|]
\bigg)\mi. 
\end{align}
Let us write $b = d_i - d_j$ and $c = d_h - d_j$. Also set 
\eq \label{eq: X_2 defi}
X_2 = X_1 \cdot \langle d_i + d_j\rangle  \langle d_h + d_j\rangle  \langle d_j + d_i +1\rangle  \langle d_j + d_h + 1\rangle
([|b-c|] [|b-c +1|] [|b-c - 1|])\mi. 
\eneq
Assuming, without loss of generality, that $b-c \geq 1$ and $c \geq 1$, it follows that  
\begin{align}
X_2\mi F(j,i,h) &= [b\smin c\smin1] \big([b\smin1][c\smin1] - q[b][c\smin1]) + [b\smin c+1] \big( - q[b\smin1][c] +q^2 [b][c]\big) \\
&= - q^{b}[b \smin c \smin 1][c \smin 1] + q^{b+1}[b \smin c + 1] [c] \\
&= q^{2b-1}(q^{3} + q^{1-2c} + q^{1-2b} - q^{-1} - q^{3-2b} - q^{3-2c})(q-q\mi)^{-2} \\
&=: q^{2b-1} Y. 
\end{align}
Observing that $Y$ is symmetric in $b$ and $c$, and that $X_2\mi F(j,h,i)$ can be obtained 
from $-X_2\mi F(j,i,h)$ by interchanging $b \leftrightarrow c$, we deduce that $X_2\mi F(j,h,i) = -q^{2c} Y$. 
Noting that 
\[
b-c = d_i - d_h = l_{i,2p-1} - l_{h,2p-1} = (m_{i,2p-1} - m_{h,2p-1}) - (i - h) = \tfrac12 (\gamma_h - \gamma_i), 
\]
we conclude
\[
\tau' + \tau'' = q^{\gamma_i}F(j,i,h) + q^{\gamma_h}F(j,h,i) = 0. 
\]

\subsubsection{Case: $\tau_{+ j,+ i,+ i}$} 
\label{subsec: jii coeff}

The total coefficient $\tau_{+j,+i,+i}$ on $\sigma_{+j,+i,+i}^{2p,2p-1,2p-1}$ is equal to the coefficient $\tau'$ on 
$q\mi[B_{2p \smin1}^{+i}, [A_{2p \smin1,\smin1}^{+i}, B_{2p}^{+j}]_{q} ]_q \cdot \sigma.$ 
Define 
\begin{align}
F(j,i) &= a\mi  \big( P^j_{2p}(\sigma) Q_{2p\smin1}^i(\sigma^{2p}_{+j}) Q_{2p\smin1}^i(\sigma^{2p,2p\smin1}_{+j,+i}) 
- [2] P^j_{2p}(\sigma_{+i}^{2p\smin1}) Q_{2p\smin1}^i(\sigma) Q_{2p\smin1}^i(\sigma^{2p,2p\smin1}_{+j,+i}) \\ 
&+  P^j_{2p}(\sigma_{+i,+i}^{2p\smin1,2p\smin1}) Q_{2p\smin1}^i(\sigma_{+i}^{2p\smin1}) Q_{2p\smin1}^i(\sigma) \big).
\end{align}
Then $\tau' = q^{-\zeta_{i,2p-1}}F(j,i)$. 
Let 
\begin{align}
X_1 &= - {\rm i} a\mi  \Bigg( \frac{1}{\{ l'_j \} \{ l'_j + 1 \}}
\frac{\prod_{r=1}^p [l_r+l'_j][|l_r-l'_j-1|]
      \prod_{r\neq i}^{p-1}[l''_r +l'_j][|l''_r -l'_j-1|] }
     {\prod_{r\ne j}^p[l'_r+l'_j][|l'_r-l'_j|]
         [l'_r+l'_j+1][|l'_r-l'_j-1|]} \\ 
&\cdot 
 \frac{1}{[\ell'_i]^2[2\ell'_i+1][2\ell'_i-1]}\frac{\prod_{r\neq j}^{p}
[\ell_r+\ell'_i][|\ell_r-\ell'_i|]
      \prod_{r=1}^{p-1}[\ell''_r +\ell'_i][|\ell''_r -\ell'_i|] }
     {\prod_{r\ne i}^{p-1}[\ell'_r+\ell'_i][\ell'_r+\ell'_i-1]
      [|\ell'_r-\ell'_i|][|\ell'_r-\ell'_i-1|]} \\
&\cdot 
 \frac{\prod_{r\neq j}^{p}
[\ell_r+\ell'_i+1][|\ell_r-\ell'_i-1|]
      \prod_{r=1}^{p-1}[\ell''_r +\ell'_i+1][|\ell''_r -\ell'_i-1|] }
     {[\ell'_i+1]^2[2\ell'_i+3][2\ell'_i+1] \prod_{r\ne i}^{p-1}[\ell'_r+\ell'_i+1][\ell'_r+\ell'_i]
      [|\ell'_r-\ell'_i-1|][|\ell'_r-\ell'_i-2|]}\Bigg)^\frac12. 
\end{align} 
Then $X_1\mi F(j,i)$ is equal to 
\begin{align}
&\langle d_i + d_j\rangle  \langle |d_i - d_j -1|\rangle   \langle d_j + d_i + 1\rangle  \langle |d_j - d_i +1|\rangle  \langle d_j + d_i + 2\rangle  \langle |d_j - d_i|\rangle  \\
- [2] &\langle d_i + d_j + 1\rangle  \langle |d_i - d_j |\rangle   \langle d_j + d_i \rangle  \langle |d_j - d_i|\rangle  \langle d_j + d_i + 2\rangle  \langle |d_j - d_i |\rangle  \\ 
+  &\langle d_i + d_j +2\rangle  \langle |d_i - d_j +1|\rangle 
 \langle d_j + d_i\rangle  \langle |d_j - d_i|\rangle  \langle d_j + d_i + 1\rangle  \langle |d_j - d_i - 1|\rangle . 
\end{align}
Assuming, without loss of generality, $b = d_i - d_j \geq 0$ and setting
\eq \label{eq: X_2 defi}
X_2 = X_1 \cdot \langle  b  \rangle \langle d_i + d_j\rangle  \langle d_j + d_i +1\rangle  \langle d_j + d_h + 2\rangle, 
\eneq
it follows that  
\begin{align}
X_2\mi F(j,i) &= [b-1] - [2][b] + [b+1] = 0. 
\end{align}

\subsubsection{Case: $\tau_{+ j,+ i}$}

The total coefficient $\tau_{+j,+i}$ on $\sigma_{+j,+i,+h}^{2p,2p-1}$ is equal to the sum of the coefficients $\tau', \tau''$ on 
\[ 
q\mi[B_{2p \smin1}^{(0)}, [A_{2p \smin1,\smin1}^{+i}, B_{2p}^{+j}]_{q} ]_q \cdot \sigma, \qquad 
q\mi[B_{2p \smin1}^{+i}, [A_{2p \smin1,\smin1}^{(0)}, B_{2p}^{+j}]_{q} ]_q \cdot \sigma, 
\]
respectively. 
Define 
\begin{align}
F_1(j,i) &= a\mi  \big( P^j_{2p}(\sigma) Q_{2p\smin1}^i(\sigma^{2p}_{+j}) R_{2p\smin1}(\sigma^{2p,2p\smin1}_{+j,+i}) 
-  qP^j_{2p}(\sigma_{+i}^{2p\smin1}) Q_{2p\smin1}^i(\sigma) R_{2p\smin1}(\sigma^{2p,2p\smin1}_{+j,+i}) \\
&-  qP^j_{2p}(\sigma) Q_{2p\smin1}^i(\sigma_{+j}^{2p}) R_{2p\smin1}(\sigma) 
+ q^2 P^j_{2p}(\sigma_{+i}^{2p\smin1}) Q_{2p\smin1}^i(\sigma) R_{2p\smin1}(\sigma) \big), \\
F_2(j,i) &= a\mi  \big( P^j_{2p}(\sigma) Q_{2p\smin1}^i(\sigma^{2p}_{+j}) R_{2p\smin1}(\sigma^{2p}_{+j}) 
-  qP^j_{2p}(\sigma) Q_{2p\smin1}^i(\sigma_{+j}^{2p}) R_{2p\smin1}(\sigma) \\
&-  qP^j_{2p}(\sigma^{2p\smin1}_{+i}) Q_{2p\smin1}^i(\sigma) R_{2p\smin1}(\sigma^{2p,2p\smin1}_{+j,+i}) 
+ q^2 P^j_{2p}(\sigma_{+i}^{2p\smin1}) Q_{2p\smin1}^i(\sigma) R_{2p\smin1}(\sigma_{+i}^{2p\smin1}) \big). 
\end{align}
Then $\tau' = q^{-\zeta_{i,2p-1}}F_1(j,i)$ and $\tau'' = F_2(j,i)$. 
Let 
\begin{align}
X_1 &= -  a\mi  \Bigg( \frac{1}{\{ l'_j \} \{ l'_j + 1 \}}
\frac{\prod_{r=1}^p [l_r+l'_j][|l_r-l'_j-1|]
      \prod_{r\neq i}^{p-1}[l''_r +l'_j][|l''_r -l'_j-1|] }
     {\prod_{r\ne j}^p[l'_r+l'_j][|l'_r-l'_j|]
         [l'_r+l'_j+1][|l'_r-l'_j-1|]} \\ 
&\cdot 
 \frac{1}{[\ell'_i]^2[2\ell'_i+1][2\ell'_i-1]}\frac{\prod_{r\neq j}^{p}
[\ell_r+\ell'_i][|\ell_r-\ell'_i|]
      \prod_{r=1}^{p-1}[\ell''_r +\ell'_i][|\ell''_r -\ell'_i|] }
     {\prod_{r\ne i}^{p-1}[\ell'_r+\ell'_i][\ell'_r+\ell'_i-1]
      [|\ell'_r-\ell'_i|][|\ell'_r-\ell'_i-1|]}\Bigg)^\frac12  \\
&\cdot 
\prod_{r\neq j}^{p} [l_{r,2p}] \prod_{r=1}^{p-1} [l_{r,2p-2}]
\left(\prod_{r\neq i}^{p-1} [l_{r,2p-1}] [l_{r,2p-1}-1]\right)^{-1}. 
\end{align} 
Then $X_1\mi F_1(j,i)$ is equal to 
\begin{align}
&\langle d_i + d_j\rangle  \langle |d_i - d_j -1|\rangle \langle d_j + d_i + 1\rangle  \langle |d_j - d_i +1|\rangle [d_j+1] \left( [d_i+1][d_i]\right)\mi  \\
-q &\langle d_i + d_j + 1\rangle  \langle |d_i - d_j |\rangle \langle d_j + d_i \rangle  \langle |d_j - d_i|\rangle  [d_j+1]\left( [d_i+1][d_i]\right)\mi   \\ 
- q &\langle d_i + d_j\rangle  \langle |d_i - d_j -1|\rangle \langle d_j + d_i + 1\rangle  \langle |d_j - d_i +1|\rangle [d_j] \left( [d_i][d_i-1]\right)\mi \\
+q^2 &\langle d_i + d_j + 1\rangle  \langle |d_i - d_j |\rangle \langle d_j + d_i \rangle  \langle |d_j - d_i|\rangle  [d_j] \left( [d_i][d_i-1]\right)\mi . 
\end{align}
Assuming, without loss of generality, that $b = d_i - d_j \geq 0$, and setting
\eq \label{eq: X_2 defi'}
X_2 = X_1 \cdot \langle d_i + d_j\rangle  \langle d_j + d_i +1\rangle \left([d_i][d_i+1][d_i-1]\right)\mi, 
\eneq
it follows that  
\begin{align}
q^{-\zeta_{i,2p-1}}X_2\mi F_1(j,i) &= q^{-d_i-d_j}(q[d_j][d_i+1] - [d_j+1][d_i-1]) \\
&= (q-q\mi)\mi (q + q^{-2d_i+1} - [2]q^{-2d_j}).  
\end{align}
%denom?

On the other hand, $X_1\mi F_2(j,i)$ is equal to 
\begin{align}
&\langle d_i + d_j\rangle  \langle |d_i - d_j -1|\rangle   \langle d_j + d_i + 1\rangle  \langle |d_j - d_i +1|\rangle [d_j+1] \left( [d_i][d_i-1]\right)\mi  \\
-q &\langle d_i + d_j \rangle  \langle |d_i - d_j-1 |\rangle  \langle d_j + d_i +1 \rangle  \langle |d_j - d_i +1|\rangle  [d_j]\left( [d_i][d_i-1]\right)\mi   \\ 
- q &\langle d_i + d_j+1\rangle  \langle |d_i - d_j |\rangle \langle d_j + d_i \rangle  \langle |d_j - d_i |\rangle [d_j+1] \left( [d_i+1][d_i]\right)\mi \\
+q^2 &\langle d_i + d_j + 1\rangle  \langle |d_i - d_j |\rangle \langle d_j + d_i \rangle  \langle |d_j - d_i|\rangle  [d_j] \left( [d_i+1][d_i]\right)\mi, 
\end{align}
which implies that 
\begin{align}
X_2\mi F_2(j,i) &= q^{-d_j}[b-1][d_i+1] - [b]q^{1-d_j}[d_i-1] \\
&= -(q-q\mi)\mi (q + q^{-2d_i+1} - [2]q^{-2d_j}).  
\end{align}
Hence $q^{-\zeta_{i,2p-1}}F_1(j,i) + F_2(j,i) = 0$. 

\subsubsection{Case: $\tau_{+j}$} 

The coefficient $\tau_{+j}$ can be written as the sum of the coefficients $\tau'_i, \tau''_i$ on 
\[ 
-q\mi[B_{2p-1}^{-i}, [A_{2p-1,\smin1}^{+i}, B_{2p}^{+j}]_{q} ]_q \cdot \sigma, \qquad 
-q\mi[B_{2p-1}^{+i}, [A_{2p-1,\smin1}^{-i}, B_{2p}^{+j}]_{q} ]_q \cdot \sigma, 
\]
respectively, for $1 \leq i \leq p-1$, and the coefficient $\tau'''$ on 
\[
q\mi[B_{2p-1}^{(0)}, [A_{2p-1,\smin1}^{(0)}, B_{2p}^{+j}]_{q} ]_q \cdot \sigma. 
\]
Let us compute them. We have 
\begin{align*}
-a\tau'_i &= q^{-2d_i-1}P^j_{2p}(\sigma) (Q^i_{2p\smin1}(\sigma^{2p}_{+j}))^2
- q^{-2d_i}P^j_{2p}(\sigma^{2p\smin1}_{+i}) Q^i_{2p\smin1}(\sigma^{2p}_{+j}) Q^i_{2p\smin1}(\sigma) \\ 
&- q^{-2d_i+2}P^j_{2p}(\sigma^{2p\smin1}_{-i}) Q^i_{2p\smin1}(\sigma^{2p,2p\smin1}_{+j,-i}) Q^i_{2p\smin1}(\sigma^{2p\smin1}_{-i}) 
+ q^{-2d_i+3}P^j_{2p}(\sigma) (Q^i_{2p\smin1}(\sigma^{2p\smin1}_{-i}))^2, \\ 
%%%%%%%%%%%%%%%%%%%
-a\tau''_i &= q^{2d_i-3}P^j_{2p}(\sigma) (Q^i_{2p\smin1}(\sigma^{2p,2p\smin1}_{+j,-i}))^2
- q^{2d_i-2}P^j_{2p}(\sigma^{2p\smin1}_{-i}) Q^i_{2p\smin1}(\sigma^{2p,2p\smin1}_{+j,-i}) Q^i_{2p\smin1}(\sigma^{2p\smin1}_{-i}) \\ 
&- q^{2d_i}P^j_{2p}(\sigma^{2p\smin1}_{+i}) Q^i_{2p\smin1}(\sigma^{2p}_{+j}) Q^i_{2p\smin1}(\sigma) 
+ q^{2d_i+1}P^j_{2p}(\sigma) (Q^i_{2p\smin1}(\sigma))^2. 
\end{align*}
Let us refer to each line on the RHS above as $D_1, \cdots, D_4$, respectively. We have 
\begin{align*}
D_1 {+} D_4  &= 
\Bigg( q^{1+2d_i} - (q^{2d_i} + q^{-2d_i})\frac{P^j_{2p}(\sigma^{2p\smin1}_{+i})}{P^j_{2p}(\sigma)} \frac{Q^i_{2p\smin1}(\sigma^{2p}_{+j})}{Q^i_{2p\smin1}(\sigma)} 
+ q^{-2d_i-1}  \left(\frac{Q^i_{2p\smin1}(\sigma^{2p}_{+j})}{Q^i_{2p\smin1}(\sigma)}\right)^2 \Bigg)  \\
&\ \ \ \cdot P^j_{2p}(\sigma)(Q^i_{2p\smin1}(\sigma))^2, \\ 
D_2 {+} D_3  &= 
\Bigg( q^{3-2d_i} - (q^{2d_i\smin2} {+} q^{2\smin2d_i})\frac{P^j_{2p}(\sigma^{2p\smin1}_{-i})}{P^j_{2p}(\sigma)} \frac{Q^i_{2p\smin1}(\sigma^{2p,2p\smin1}_{+j,-i})}{Q^i_{2p\smin1}(\sigma^{2p\smin1}_{-i})} 
+ q^{2d_i\smin3}  \left(\frac{Q^i_{2p\smin1}(\sigma^{2p,2p\smin1}_{+j,-i})}{Q^i_{2p\smin1}(\sigma^{2p\smin1}_{-i})}\right)^2 \Bigg)  \\
&\ \ \ \cdot P^j_{2p}(\sigma)(Q^i_{2p\smin1}(\sigma^{2p\smin1}_{-i}))^2, \\ 
a\tau''' &= \Bigg( (q - 2 \frac{R_{2p-1}(\sigma^{2p}_{+j})}{R_{2p-1}(\sigma)} + q\mi \left(\frac{R_{2p-1}(\sigma^{2p}_{+j})}{R_{2p-1}(\sigma)} \right)^2 \Bigg) \cdot P^j_{2p}(\sigma) (R_{2p-1}(\sigma))^2. 
\end{align*}
It follows from the definitions that 
\begin{align*}
\frac{P^j_{2p}(\sigma^{2p\smin1}_{+i})}{P^j_{2p}(\sigma)} \frac{Q^i_{2p\smin1}(\sigma^{2p}_{+j})}{Q^i_{2p\smin1}(\sigma)} 
&= \frac{[d_i{+}d_j{+}1]}{[d_i{+}d_j]}, \quad &
\left(\frac{Q^i_{2p\smin1}(\sigma^{2p}_{+j})}{Q^i_{2p\smin1}(\sigma)}\right)^2 &= \frac{[d_i{+}d_j{+}1][|d_i{-}d_j{-}1|]}{[d_i{+}d_j][|d_i{-}d_j|]}, \\ 
\frac{P^j_{2p}(\sigma^{2p\smin1}_{-i})}{P^j_{2p}(\sigma)} \frac{Q^i_{2p\smin1}(\sigma^{2p,2p\smin1}_{+j,-i})}{Q^i_{2p\smin1}(\sigma^{2p\smin1}_{-i})} 
&= \frac{[|d_i{-}d_j{-}2|]}{[|d_i{-}d_j{-}1|]}, \quad &
\left(\frac{Q^i_{2p\smin1}(\sigma^{2p,2p\smin1}_{+j,-i})}{Q^i_{2p\smin1}(\sigma^{2p\smin1}_{-i})}\right)^2 &= \frac{[d_i{+}d_j][|d_i{-}d_j{-}2]}{[d_i{+}d_j{-}1][|d_i{-}d_j{-}1|]}, \\ 
\frac{R_{2p-1}(\sigma^{2p}_{+j})}{R_{2p-1}(\sigma)} 
&= \frac{[d_j+1]}{[d_j]}. 
\end{align*} 
A direct calculation now yields 
\begin{align*}
D_1 + D_4 &= \frac{q^{-2d_j-1}[2d_i+1]}{[d_i+d_j][|d_i-d_j|]}P^j_{2p}(\sigma)(Q^i_{2p\smin1}(\sigma))^2, \\ 
D_2 + D_3 &= \frac{-q^{-2d_j-1}[2d_i-3]}{[d_i+d_j-1][|d_i-d_j-1|]} P^j_{2p}(\sigma)(Q^i_{2p\smin1}(\sigma^{2p\smin1}_{-i}))^2, \\ 
a\tau''' &= \frac{q^{-2d_j-1}}{[d_j]^2} P^j_{2p}(\sigma) (R_{2p-1}(\sigma))^2. 
\end{align*} 

We have thus shown that $\tau_{+j}$ is equal to $a\mi q^{-2d_j-1} P^j_{2p}(\sigma)$ times the formula on the LHS of \cite[(27)]{GI97}, which is proven to be equal to $1$ in \emph{loc.\ cit.} The proof relies on Lemma \ref{pro: gen Appendix pro}, with replacements 
\begin{align*}
l''_{i} &\mapsto l_{i,2p-2}+\tfrac{1}{2}, \quad (i=1,\cdots,p-1), \\
l'_{i} &\mapsto l_{i,2p-1}-\tfrac{1}{2}, \quad (i=1,\cdots,p-1), \\
l_{i} &\mapsto l_{i,2p}+\tfrac{1}{2}, \quad (i=1,\cdots,p; \ i \not= j), \\ 
l_{j} &\mapsto \tfrac{3}{2}, \qquad l'_{p} \mapsto -\tfrac{1}{2}. 
\end{align*}

\subsection{Horizontal induction: $A_{2p+1,-1}$} 
\label{subsec: hor ind 2}

It follows from \eqref{lem: Serrcor} and induction that 
\begin{align}
 A_{2p+1 ,\smin1} \cdot \sigma &= \sum_{j,i \neq h} \tau_{\star j,\star i,\star h} \sigma_{\star j,\star i,\star h}^{2p+1,2p,2p}
 + \sum_{j,i} \tau_{\star j,\pm i,\pm i} \sigma_{\star j,\pm i,\pm i}^{2p+1,2p,2p} 
+ \sum_{j} \tau_{\pm j} \sigma_{\pm j}^{2p+1} \\ 
&+ \sum_{i \neq h} \tau_{\star i,\star h} \sigma_{\star i,\star h}^{2p,2p}
 + \sum_{i} \tau_{\pm i,\pm i} \sigma_{\pm i,\pm i}^{2p,2p} 
+ \tau \sigma. 
\end{align}
The proof that the coefficients $\tau_{\star j,\star i,\star h}$, $\tau_{\star j,\pm i,\pm i}$, $\tau_{\star i,\star h}$ and $\tau_{\pm i,\pm i}$ vanish is analogous to the corresponding proofs in \S \ref{subsec: jih coeff}--\ref{subsec: jii coeff}. Let us compute the remaining two coefficients. In this subsection, we set $d_j = l_{j,2p+1}, d_i = l_{i,2p}$.

\subsubsection{Case: $\tau_{+j}$} 

The coefficient $\tau_{+j}$ can be written as the sum of the coefficients $\tau'_i, \tau''_i$ on 
\[ 
-q\mi[B_{2p}^{-i}, [A_{2p,\smin1}^{+i}, B_{2p+1}^{+j}]_{q} ]_q \cdot \sigma, \qquad 
-q\mi[B_{2p}^{+i}, [A_{2p,\smin1}^{-i}, B_{2p+1}^{+j}]_{q} ]_q \cdot \sigma, 
\]
respectively, for $1 \leq i \leq p$. Let us compute them. We have 
\begin{align*}
-a\tau'_i &= q^{-2d_i-2}Q^j_{2p+1}(\sigma) (P^i_{2p}(\sigma^{2p+1}_{+j}))^2
- q^{-2d_i-1}Q^j_{2p+1}(\sigma^{2p}_{+i}) P^i_{2p}(\sigma^{2p+1}_{+j}) P^i_{2p}(\sigma) \\ 
&- q^{-2d_i+1}Q^j_{2p+1}(\sigma^{2p}_{-i}) P^i_{2p}(\sigma^{2p+1,2p}_{+j,-i}) P^i_{2p}(\sigma^{2p}_{-i}) 
+ q^{-2d_i+2}Q^j_{2p+1}(\sigma) (P^i_{2p}(\sigma^{2p}_{-i}))^2, \\ 
%%%%%%%%%%%%%%%%%%%
-a\tau''_i &= q^{2d_i-2}Q^j_{2p+1}(\sigma) (P^i_{2p}(\sigma^{2p+1,2p}_{+j,-i}))^2
- q^{2d_i-1}Q^j_{2p+1}(\sigma^{2p}_{-i}) P^i_{2p}(\sigma^{2p+1,2p}_{+j,-i}) P^i_{2p}(\sigma^{2p}_{-i}) \\ 
&- q^{2d_i+1}Q^j_{2p+1}(\sigma^{2p}_{+i}) P^i_{2p}(\sigma^{2p+1}_{+j}) P^i_{2p}(\sigma) 
+ q^{2d_i+2}Q^j_{2p+1}(\sigma) (P^i_{2p}(\sigma))^2. 
\end{align*}
Let us refer to each line on the RHS above as $D_1, \cdots, D_4$, respectively. We have 
\begin{align*}
D_1 {+} D_4  &= 
\Bigg( q^{2+2d_i} - (q^{2d_i+1} {+} q^{-2d_i-1})\frac{Q^j_{2p+1}(\sigma^{2p}_{+i})}{Q^j_{2p+1}(\sigma)} \frac{P^i_{2p}(\sigma^{2p+1}_{+j})}{P^i_{2p}(\sigma)} 
+ q^{-2d_i-2}  \left(\frac{P^i_{2p}(\sigma^{2p+1}_{+j})}{P^i_{2p}(\sigma)}\right)^2 \Bigg)  \\
&\ \ \ \cdot Q^j_{2p+1}(\sigma)(P^i_{2p}(\sigma))^2, \\ 
D_2 {+} D_3  &= 
\Bigg( q^{2-2d_i} - (q^{2d_i\smin1} {+} q^{1\smin2d_i})\frac{Q^j_{2p+1}(\sigma^{2p}_{-i})}{Q^j_{2p+1}(\sigma)} \frac{P^i_{2p}(\sigma^{2p+1,2p}_{+j,-i})}{P^i_{2p}(\sigma^{2p}_{-i})} 
+ q^{2d_i\smin2}  \left(\frac{P^i_{2p}(\sigma^{2p+1,2p}_{+j,-i})}{P^i_{2p}(\sigma^{2p}_{-i})}\right)^2 \Bigg)  \\
&\ \ \ \cdot Q^j_{2p+1}(\sigma)(P^i_{2p}(\sigma^{2p}_{-i}))^2. 
\end{align*}
It follows from the definitions that 
\begin{align*}
\frac{Q^j_{2p+1}(\sigma^{2p}_{+i})}{Q^j_{2p+1}(\sigma)} \frac{P^i_{2p}(\sigma^{2p+1}_{+j})}{P^i_{2p}(\sigma)} 
&= \frac{[d_i{+}d_j{+}1]}{[d_i{+}d_j]}, \quad &
\left(\frac{P^i_{2p}(\sigma^{2p+1}_{+j})}{P^i_{2p}(\sigma)}\right)^2 &= \frac{[d_i{+}d_j{+}1][|d_i{-}d_j|]}{[d_i{+}d_j][|d_i{-}d_j+1|]}, \\ 
\frac{Q^j_{2p+1}(\sigma^{2p}_{-i})}{Q^j_{2p+1}(\sigma)} \frac{P^i_{2p}(\sigma^{2p+1,2p}_{+j,-i})}{P^i_{2p}(\sigma^{2p}_{-i})} 
&= \frac{[|d_i{-}d_j{-}1|]}{[|d_i{-}d_j|]}, \quad &
\left(\frac{P^i_{2p}(\sigma^{2p+1,2p}_{+j,-i})}{P^i_{2p}(\sigma^{2p}_{-i})}\right)^2 &= \frac{[d_i{+}d_j][|d_i{-}d_j{-}1]}{[d_i{+}d_j{-}1][|d_i{-}d_j|]}. 
\end{align*} 
A direct calculation now yields 
\begin{align*}
D_1 + D_4 &= \frac{-q^{-2d_j}[2d_i+2]}{[d_i+d_j][|d_i-d_j+1|]}Q^j_{2p+1}(\sigma)(P^i_{2p}(\sigma))^2, \\ 
D_2 + D_3 &= \frac{q^{-2d_j}[2d_i-2]}{[d_i+d_j-1][|d_i-d_j|]} Q^j_{2p+1}(\sigma)(P^i_{2p}(\sigma^{2p}_{-i}))^2. 
\end{align*}

We have thus shown that $\tau_{+j}$ is equal to $a\mi q^{-2d_j} Q^j_{2p+1}(\sigma)$ times the formula on the LHS of \cite[(24)]{GI97}, which is proven to be equal to $1$ in \emph{loc.\ cit.}, using Lemma \ref{pro: gen Appendix pro}. 
The proof for the coefficient $\tau_{-j}$ is analogous. 

%%%%%%%%%

\subsubsection{Case: $\tau$} 

The coefficient $\tau_{+j}$ can be written as the sum of the coefficients $\tau'_i, \tau''_i$ on 
\[ 
q\mi[B_{2p}^{-i}, [A_{2p,\smin1}^{+i}, B_{2p+1}^{(0)}]_{q} ]_q \cdot \sigma, \qquad 
q\mi[B_{2p}^{+i}, [A_{2p,\smin1}^{-i}, B_{2p+1}^{(0)}]_{q} ]_q \cdot \sigma, 
\]
respectively, for $1 \leq i \leq p$. Let us compute them. We have 
\begin{align*}
a\tau'_i &= \big(q^{-2d_i-2}R^j_{2p+1}(\sigma) - q^{-2d_i-1}R^j_{2p+1}(\sigma^{2p}_{+i}) \big) (P^i_{2p}(\sigma))^2
 \\ 
&+ \big(q^{-2d_i+2}R^j_{2p+1}(\sigma) - q^{-2d_i+1}R^j_{2p+1}(\sigma^{2p}_{-i}) \big) (P^i_{2p}(\sigma^{2p}_{-i}))^2, \\ 
%%%%%%%%%%%%%%%%%%%
a\tau''_i &= \big( q^{2d_i-2}R^j_{2p+1}(\sigma) - q^{2d_i-1}R^j_{2p+1}(\sigma^{2p}_{-i}) \big) 
 (P^i_{2p}(\sigma^{2p}_{-i}))^2 \\ 
&+ \big(q^{2d_i+2}R^j_{2p+1}(\sigma) - q^{2d_i+1}R^j_{2p+1}(\sigma^{2p}_{+i})   \big) (P^i_{2p}(\sigma))^2. 
\end{align*}
Let us refer to each line on the RHS above as $D_1, \cdots, D_4$, respectively. We have 
\begin{align*}
D_1 {+} D_4  &= 
\Bigg( (q^{2+2d_i}+ q^{-2d_i-2})  - (q^{2d_i+1} + q^{-2d_i-1})\frac{R^j_{2p+1}(\sigma^{2p}_{+i})}{R^j_{2p+1}(\sigma)} \Bigg) \cdot R^j_{2p+1}(\sigma)(P^i_{2p}(\sigma))^2, \\ 
D_2 {+} D_3  &= 
\Bigg( (q^{2-2d_i} + q^{2d_i-2}) - (q^{2d_i\smin1} {+} q^{1\smin2d_i})\frac{R^j_{2p+1}(\sigma^{2p}_{-i})}{R^j_{2p+1}(\sigma)} \Bigg) \cdot R^j_{2p+1}(\sigma)(P^i_{2p}(\sigma^{2p}_{-i}))^2. 
\end{align*}
It follows from the definitions that 
\[
\frac{R^j_{2p+1}(\sigma^{2p}_{+i})}{R^j_{2p+1}(\sigma)} 
= \frac{[d_i + 1]}{[d_i]}, \qquad 
\frac{R^j_{2p+1}(\sigma^{2p}_{-i})}{R^j_{2p+1}(\sigma)} 
= \frac{[d_i - 1]}{[d_i]}. 
\]
A direct calculation now yields 
\begin{align*}
D_1 + D_4 &= - \frac{\{d_i + 1\}}{[d_i]} R^j_{2p+1}(\sigma)(P^i_{2p}(\sigma))^2, \\ 
D_2 + D_3 &= \frac{\{d_i - 1\}}{[d_i]}  R^j_{2p+1}(\sigma)(P^i_{2p}(\sigma^{2p}_{-i}))^2. 
\end{align*}

We have thus shown that $\tau$ is equal to $-a\mi R^j_{2p+1}(\sigma)$ times the formula on the LHS of \cite[(21)]{GI97}, which is proven to be equal to $1$ in \emph{loc.\ cit.}, using Lemma \ref{pro: gen Appendix pro}.

%%%%%%%%%%%%

\subsection{Vertical induction} 
\label{subsec: ver ind}

Since the special case of $A_{1,s}$ is easy to compute, we will focus on $A_{k,s}$ for $k \geq 2$. Having computed the action of $A_{k,-1}$ on $\sigma \in \GT(\boldsymbol\lambda)$ in \S \ref{subsec: hor ind}--\ref{subsec: hor ind 2}, and knowing the action of $A_{k,0} = B_k$ from Theorem \ref{thm: GIK B-action on GT}, we derive the action of $A_{k,s}$, for arbitrary $s \in \Z$, through downward and upward induction. Since the proofs are the same, we will only do the upward case explicitly. 
We will need the following special case of defining relation \eqref{eq: grel2}: 
\[
\Ag_{k,s+1} =  - [\Hg_{k-1,1}, \Ag_{k,s}] + a^2\Ag_{k,s-1}. 
\] 
By the `horizontal induction', we can assume that Theorem \ref{thm: Theta-action on GT} holds for $k-1$. In particular, we deduce that $H_{k-1,1}$ acts on~$\sigma$ by the following scalar: 
\begin{align*}
 \frac{a}{q-q\mi} \sum_{r=1}^2 \sum_{i =1}^{\lceil (k-1)/2 \rceil} \left( q^{(-1)^r \zeta_{i ,k-1} } +  q^{(-1)^r (\zeta_{i ,k-1}-2 )} - q^{(-1)^r (\zeta_{i ,k}-1)} - q^{(-1)^r (\zeta_{i ,k-2}-1)} \right). 
\end{align*}
Let us denote this scalar, divided by $a$, as $\eta(k{\smin}1;\sigma)$. 

It suffices to show that the action of $\left( [\Hg_{k-1,1}, \Ag_{k,s}^\star] +a^2\Ag_{k,s-1}^\star \right)$ on $\sigma$ is the same way as the action of $\Ag_{k,s+1}^\star$ described in Theorem \ref{thm: A-action on GT}, 
for each $\star \in \{ \pm j \}$ and $\star = (0)$ if $k$ is odd. Let us do the `$+j$' case for $k=2p$ explicitly. One needs to show that 
\eq \label{eq: vert ind calc}
\left(-  [\Hg_{2p-1,1}, \Ag_{2p,s}^{+j}] +a^2\Ag_{2p,s-1}^{+j} \right) \cdot \sigma = a^{s+1}q^{(s+1) \cdot \zeta_{j,2p}} P^j_{2p}(\sigma)\sigma^{2p}_{+j}. 
\eneq
By induction, the LHS of \eqref{eq: vert ind calc} is equal to 
\begin{align*}
&\big( 
1 - (q-q\mi)\mi \left(\eta(2p{\smin}1;\sigma^{2p}_{+j}) - \eta(2p{\smin}1;\sigma) \right) q^{\zeta_{j,2p}}\big) \cdot a^{s+1} q^{(s-1) \cdot \zeta_{j,2p}} P^j_{2p}(\sigma) \sigma^{2p}_{+j} \\
&= a^{s+1}q^{(s+1) \cdot \zeta_{j,2p}} P^j_{2p}(\sigma)\sigma^{2p}_{+j}, 
\end{align*} 
as desired. 
The other cases are similar.

%%%%%%%%%%%

\subsection{$\thvar{k}$-eigenvalues} 
\label{subsec: Theta on GT}

We are now ready to prove Theorem \ref{thm: Theta-action on GT}. 
In our setting, identity \eqref{eq: Theta from A} specializes to the following formula: 
\begin{equation} \label{eq: a->Theta spec}
\thvar{k} = 1 - \frac{q(q-q\mi)a^2 z^2}{1 - a^2 z^2}\left(z\mi[A_{k,-1}, \avar{k}]_{q^{-2}} - q^{-2}[A_{k,0}, \avar{k}]_{q^2} \right). 
\end{equation} 

\subsubsection{Case: $k=2p$} 
\label{subsec: 2p eigenv}

It follows from Theorem \ref{thm: A-action on GT} that $\thvar{k} \cdot \sigma$ can be written as a linear combination of $\sigma$ and, possibly, some other Gelfand--Tsetlin patterns. It is straightforward (though somewhat laborious) to check that, in fact, the coefficients on these other Gelfand--Tsetlin patterns vanish. Given that we have performed similar calculations in \S \ref{subsec: hor ind}--\ref{subsec: hor ind 2}, we will not repeat them. Instead, let us compute the eigenvalue on $\sigma$. 
Using \eqref{eq: a->Theta spec}, we can compute the different components contributing to the eigenvalue: 
\begin{align*}
A_{2p,-1}^{-i} \ivar{2p}{+i} \cdot \sigma &= \frac{q^{\zeta_{i,2p}}(P^i_{2p}(\sigma))^2}{-a(1-q^{\zeta_{i,2p}}az)} \sigma, &\quad 
\ivar{2p}{+i} A_{2p,-1}^{-i} \cdot \sigma &= \frac{q^{\zeta_{i,2p}-2}(P^i_{2p}(\sigma_{-i}^{2p}))^2}{-a(1-q^{\zeta_{i,2p}-2}az)} \sigma, \\ 
%%%%%%%%
A_{2p,-1}^{+i} \ivar{2p}{-i} \cdot \sigma &= \frac{q^{2-\zeta_{i,2p}}(P^i_{2p}(\sigma_{-i}^{2p}))^2}{-a(1-q^{2-\zeta_{i,2p}}az)} \sigma, &\quad 
\ivar{2p}{-i} A_{2p,-1}^{+i} \cdot \sigma &= \frac{q^{-\zeta_{i,2p}}(P^i_{2p}(\sigma))^2}{-a(1-q^{-\zeta_{i,2p}}az)} \sigma, \\ 
%%%%%%%%%%%%%%%%%%%%%%%%%%%%%
A_{2p,0}^{-i} \ivar{2p}{+i} \cdot \sigma &= \frac{(P^i_{2p}(\sigma))^2}{-(1-q^{\zeta_{i,2p}}az)} \sigma, &\quad 
\ivar{2p}{+i} A_{2p,0}^{-i} \cdot \sigma &= \frac{(P^i_{2p}(\sigma_{-i}^{2p}))^2}{-(1-q^{\zeta_{i,2p}-2}az)} \sigma, \\ 
%%%%%%%%
A_{2p,0}^{+i} \ivar{2p}{-i} \cdot \sigma &= \frac{(P^i_{2p}(\sigma_{-i}^{2p}))^2}{-(1-q^{2-\zeta_{i,2p}}az)} \sigma, &\quad 
\ivar{2p}{-i} A_{2p,0}^{+i} \cdot \sigma &= \frac{(P^i_{2p}(\sigma))^2}{-(1-q^{-\zeta_{i,2p}}az)} \sigma. 
\end{align*}
Plugging these into \eqref{eq: a->Theta spec}, we compute that $\thvar{2p}$ acts on $\sigma$ by a scalar equal to  
\begin{align*}
1 + \sum_{i=1}^p za(q-q\mi) \left(
\frac{(q^{\zeta_{i,2p}+1}-q^{-\zeta_{i,2p}-1})(P^i_{2p}(\sigma))^2}{(1-q^{\zeta_{i,2p}}az)(1-q^{-\zeta_{i,2p}}az)} 
- \frac{(q^{\zeta_{i,2p}-3}-q^{3-\zeta_{i,2p}})(P^i_{2p}(\sigma_{-i}^{2p}))^2}{(1-q^{\zeta_{i,2p}-2}az)(1-q^{2-\zeta_{i,2p}}az)}
\right). 
\end{align*} 
Using the definitions of $\zeta_{i,2p}$ and $P^i_{2p}(\sigma)$, one recognizes immediately that the formula above can be converted into $\Psi(\trip)$, with 
\[
\trip = 
((l_{1,2p+1},l_{1,2p},l_{1,2p-1}),\cdots,
(l_{p-1,2p+1},l_{p-1,2p},l_{p-1,2p-1}),(l_{p,2p+1},l_{p,2p},\cdot)).  
\] 
Theorem \ref{thm: Theta-action on GT} now follows directly as an application of Lemma \ref{pro: gen Appendix pro2}.

\subsubsection{Case: $k=2p-1$}

Similar calculations to those in \S \ref{subsec: 2p eigenv} show that $\thvar{2p-1}$ acts on $\sigma$ by a scalar equal to    
\begin{align*}
1 &+ \sum_{i=1}^{p-1} za(q-q\mi) \Big(
\frac{(q^{\zeta_{i,2p-1}+1}-q^{-\zeta_{i,2p-1}-1})(Q^i_{2p-1}(\sigma))^2}{(1-q^{\zeta_{i,2p-1}}az)(1-q^{-\zeta_{i,2p-1}}az)} \\
&- \frac{(q^{\zeta_{i,2p-1}-3}-q^{3-\zeta_{i,2p-1}})(Q^i_{2p-1}(\sigma_{-i}^{2p-1}))^2}{(1-q^{\zeta_{i,2p-1}-2}az)(1-q^{2-\zeta_{i,2p-1}}az)}
\Big) - \frac{za(q-q\mi)^2}{(1-az)^2} (R_{2p-1}(\sigma))^2. 
\end{align*} 
Theorem \ref{thm: Theta-action on GT} again follows as an application of Lemma \ref{pro: gen Appendix pro2}. Given the analogy with the even case, we omit further details.

\subsection{Non-classical evaluation modules} 
\label{subsec: nc Theta on GT}

Let $\boldsymbol\lambda_+ \in \Lambda_n^+$ be a regularly dominant weight, and consider the non-classical evaluation module $L_{\boldsymbol\lambda_+,\vec{\varepsilon}}(a) = \eva^*(L_{\boldsymbol\lambda_+,\vec{\varepsilon}})$. Below we state the non-classical analogues of Theorems  \ref{thm: A-action on GT} and \ref{thm: Theta-action on GT}. Since the proofs are similar to the classical case, we do not replicate them. 

\Thm 
\label{thm: A-action on GT nc}
Let $\sigma \in \GTn(\boldsymbol\lambda_+)$. Then 
\begin{align*} 
A_{2p,s} \cdot \sigma &= \sum^p_{j=1} (-a)^sq^{s \cdot \zeta_{j,2p}} \dot{P}^j_{2p}(\sigma)\sigma^{2p}_{+j}
-\sum^p_{j=1} (-a)^sq^{s \cdot (2-\zeta_{j,2p})} \dot{P}^j_{2p}(\sigma^{2p}_{-j})\sigma^{2p}_{-j} \\ 
&+  \delta_{m_{p,2p},\frac{1}{2}} \frac{(-a)^s \varepsilon_{2p+1}}{q^{\frac12} - q^{-\frac12}} \dot{S}_{2p}(\sigma) \sigma, \\
%%%%%%%%%%%%%
A_{2p-1,s} \cdot \sigma &= \sum^{p-1}_{j=1} (-a)^sq^{s \cdot \zeta_{j,2p-1}} \dot{Q}^j_{2p-1}(\sigma) \sigma^{2p-1}_{+j}  - \sum^{p-1}_{j=1} (-a)^s q^{s \cdot (2-\zeta_{j,2p-1})} \dot{Q}^j_{2p-1}({\sigma}^{2p-1}_{-j}) {\sigma}^{2p-1}_{-j} \\
&+ a^s \varepsilon_{2p} \dot{R}_{2p-1}(\sigma) \sigma. 
\end{align*} 
for each $s \in \Z$. 
\enthm 

\Thm 
\label{thm: Theta-action on GT nc}
Let $\sigma \in \GTn(\boldsymbol\lambda_+)$. Then $\sigma$ is an eigenvector for the action of the commutative algebra $\mathcal{H}$. The corresponding eigenvalues are given by 
\begin{align*}
\thvar{2p} 
\cdot \sigma &= 
\prod_{j=1}^{p-1} \frac{(1+ q^{2l_{j,2p+1}\smin 1}az)(1+ q^{1\smin 2l_{j,2p+1}}az)}{(1+ q^{2l_{j,2p}+1}az)(1+ q^{\smin 1\smin 2l_{j,2p}}az)} \frac{(1+ q^{2l_{j,2p\smin 1}\smin 1}az)(1+ q^{1\smin 2l_{j,2p\smin 1}}az)}{(1+ q^{2l_{j,2p}\smin 1}az)(1+ q^{1\smin 2l_{j,2p}}az)} \cdot \\
&\cdot \frac{(1-qaz)(1-q^{-1}az)(1+q^{2l_{p,2p+1}-1}az)(1+q^{1-2l_{p,2p+1}}az)}{ (1+q^{2l_{p,2p}-1}az)(1+q^{1-2l_{p,2p}}az) (1+q^{2l_{p,2p}+1}az) (1+q^{-1-2l_{p,2p}}az) } \sigma, \\
\thvar{2p\smin 1} 
\cdot \sigma &=\prod_{j=1}^{p-1} 
\frac{(1+q^{2l_{j,2p}}az) (1+q^{-2l_{j,2p}}az)}{(1+q^{2l_{j,2p-1}}az)(1+q^{-2l_{j,2p-1}}az)} 
\frac{(1+q^{2l_{j,2p-2}}az)(1+q^{-2l_{j,2p-2}}az)}{(1+q^{2l_{j,2p-1}-2}az)(1+q^{2-2l_{j,2p-1}}az)} \cdot  \\
&\cdot \frac{(1+q^{2l_{p,2p}}az) (1+q^{-2l_{p,2p}}az)}{(1-az)^2} \sigma. 
\end{align*} 
\enthm

\Rem \label{rem: nice formula Theta nc}
We can express the formulae above in the following uniform way, using $\zeta$-coordinates: 
\begin{align*}
\thvar{k} 
\cdot \sigma = 
\prod_{j=1}^{\lceil k/2 \rceil} 
&\frac{(1+(-1)^k(-q)^{\zeta_{j,k+1}-1}az)(1+(-1)^k(-q)^{1-\zeta_{j,k+1}}az)}{(1+(-1)^k(-q)^{\zeta_{j,k}}az)(1+(-1)^k(-q)^{-\zeta_{j,k}}az)} \\
&\frac{(1+(-1)^k(-q)^{\zeta_{j,k-1}-1}az)(1+(-1)^k(-q)^{1-\zeta_{j,k-1}}az)}{(1+(-1)^k(-q)^{\zeta_{j,k}-2}az)(1+(-1)^k(-q)^{2-\zeta_{j,k}}az)} \sigma, 
\end{align*}
where, for each $p > 0$, by convention, we set $\zeta_{p,2p-1}=0$ and $\zeta_{p,2p-2} = -1$. Note that when $k$ is odd, one of the fractions on the RHS reduces to $1$. \rmkend
\enrem

%%%%%%%%%%%%%%
%%%%%%%%%%%%%%
%%%%%%%%%%%%%%

\section{Boundary $q$-characters of evaluation representations} 
\label{sec: app to bq char}

Below we derive various consequences of Theorems \ref{thm: Theta-action on GT} and \ref{thm: Theta-action on GT nc}, and illustrate them with examples. 

\subsection{Character formulae} 

We begin by deducing a general formula for the boundary $q$-characters of classical and non-classical evaluation modules. 

\subsubsection{Classical modules} 

Let $\boldsymbol\lambda \in \Lambda_n$ be a regularly dominant weight. Given $\sigma \in \GT(\boldsymbol\lambda)$, 
define 
\[
\kappa_1(j,k) = (\zeta_{j,k+1}-\zeta_{j,k}-1)/2, \qquad 
\kappa_2(j,k) = (\zeta_{j,k}-\zeta_{j,k-1}-1)/2. 
\]
Note that if $k$ is odd then $\kappa_2(\lceil k/2 \rceil, k) = 0$. 
If $\boldsymbol\lambda$ is integral then $\kappa_1(j,k)$ and $\kappa_2(j,k)$ are always non-negative integers. If $\boldsymbol\lambda$ is half-integral then the same is true except for $\kappa_1(\lceil k/2 \rceil, k)$ (if $k$ is odd) and $\kappa_2(k/2, k)$ (if $k$ is even), which are in $\frac12 + \Z_{\geq 0}$. Accordingly, for $i \in \{1,2\}$, we set
\[
\tilde{\kappa}_i(j,k) = \lfloor \kappa_i(j,k) \rfloor. 
\]

\Cor \label{cor: iq char formula classical}
The boundary $q$-character of the classical simple evaluation module $L_{\boldsymbol\lambda}(a)$ is given~by
\[
\ichmap\left(  L_{\boldsymbol\lambda}(a) \right) = \sum_{\sigma \in \GT(\boldsymbol\lambda)} \Xi(\sigma),
\]
where: 
\begin{enumerate}
\item if $\boldsymbol\lambda \in \Lambda^{\Z}_n$ is integral then 
\eq \label{eq: ichar int}
\Xi(\sigma) = 
\prod_{k=1}^n \prod_{j=1}^{\lceil k/2 \rceil} \prod_{i=1}^{\kappa_1(j,k)} \ourY_{k,q^{\zeta_{j,k}+2i-1}}\mi 
\prod_{h=1}^{\kappa_2(j,k)} \ourY_{k,q^{\zeta_{j,k-1}+2h-2}}; 
\eneq
\item if $\boldsymbol\lambda \in \Lambda^{\frac12 + \Z}_n$ is half-integral then 
\eq \label{eq: ichar halfint}
\Xi(\sigma) =  
\prod_{k=1}^n
Z(\sigma,k)
 \prod_{j=1}^{\lceil k/2 \rceil} \prod_{i=1}^{\tilde{\kappa}_1(j,k)} \ourY_{k,q^{\zeta_{j,k}+2i-1}}\mi 
\prod_{h=1}^{\tilde{\kappa}_2(j,k)} \ourY_{k,q^{\zeta_{j,k-1}+2h-2}}, 
\eneq
\end{enumerate} 
where
\[
Z(\sigma,k) = \left\{ \begin{array}{l l}
\ourX\mi_{k, q^{\zeta_{\lceil k/2 \rceil, k+1} - \frac32}}  & \mbox{ if } k \mbox{ is odd}, \\[2pt]
\ourX_{k, q^{\zeta_{k/2, k} - \frac52}} & \mbox{ if } k \mbox{ is even}. 
\end{array} \right.
\]
\encor

\Proof
This follows directly from Theorem \ref{thm: Theta-action on GT}, Remark \ref{rem: nice formula Theta} and Lemma \ref{lem: monomial vs eigen}.   
\enproof

In particular, Corollary \ref{cor: iq char formula classical} implies that 
$
\ichmap\left(  L_{\boldsymbol\lambda}(a) \right) \in \Xmod_a. 
$ 
It is useful to rewrite \eqref{eq: ichar halfint} in the form 
\eq 
\Xi(\sigma) =  
\prod_{k=1}^n
Z'(\sigma,k)
 \prod_{j=1}^{\lceil k/2 \rceil-1} \prod_{i=1}^{\kappa_1(j,k)} \ourY_{k,q^{\zeta_{j,k}+2i-1}}\mi 
\prod_{h=1}^{\kappa_2(j,k)} \ourY_{k,q^{\zeta_{j,k-1}+2h-2}}, 
\eneq
with 
\[
Z'(\sigma,k) = \left\{ \begin{array}{l l}
\ourX\mi_{k, q^{\frac12}} 
\prod_{i=1}^{m_{(k+1)/2,k+1}-\frac12} \ourY_{k,q^{2i}}\mi   & \mbox{ if } k \mbox{ is odd}, \\[2pt]
\ourX\mi_{k, q^{\frac12}} \prod_{i=1}^{\kappa_1(k/2,k)} \ourY_{k,q^{\zeta_{k/2,k}+2i-1}}\mi \prod_{h=1}^{m_{k/2,k}-\frac12} \ourY_{k,q^{2h-1}}   & \mbox{ if } k \mbox{ is even}. 
\end{array} \right. 
\]

\subsubsection{Non-classical modules} 

Let $\boldsymbol\lambda_+ \in \Lambda_n^+$. Given $\sigma \in \GTn(\boldsymbol\lambda_+)$, we define $\tilde{\kappa}_i(j,k)$ in the same way as in the half-integral classical case. 

\Cor \label{cor: iq char formula non-classical}
The boundary $q$-character of the non-classical simple evaluation module $L_{\boldsymbol\lambda_+, {\vec{\varepsilon}}}(a)$ is given~by
\[
\ichmap\left( L_{\boldsymbol\lambda_+, {\vec{\varepsilon}}}(a) \right) = \sum_{\sigma \in \GTn(\boldsymbol\lambda_+)} \Xi(\sigma),
\]
where
\eq 
\Xi(\sigma) =  
\prod_{k=1}^n
Z'(\sigma,k)
 \prod_{j=1}^{\lceil k/2 \rceil-1} \prod_{i=1}^{\kappa_1(j,k)} \ourY_{k,-q^{\zeta_{j,k}+2i-1}}\mi 
\prod_{h=1}^{\kappa_2(j,k)} \ourY_{k,-q^{\zeta_{j,k-1}+2h-2}}, 
\eneq
with 
\[
Z'(\sigma,k) = \left\{ \begin{array}{l l}
\widetilde{\ourX}\mi_{k, q^{\frac12}} 
\prod_{i=1}^{m_{(k+1)/2,k+1}-\frac12} \ourY_{k,-q^{2i}}\mi   & \mbox{ if } k \mbox{ is odd}, \\[2pt]
\widetilde{\ourX}\mi_{k, q^{\frac12}} \prod_{i=1}^{\kappa_1(k/2,k)} \ourY_{k,-q^{\zeta_{k/2,k}+2i-1}}\mi \prod_{h=1}^{m_{k/2,k}-\frac12} \ourY_{k,-q^{2h-1}}   & \mbox{ if } k \mbox{ is even}. 
\end{array} \right. 
\]
In particular, the boundary $q$-character is independent of $\vec{\varepsilon}$. 
\encor

\Proof
This follows directly from Theorem \ref{thm: Theta-action on GT nc}, Remark \ref{rem: nice formula Theta nc} and Lemma \ref{lem: monomial vs eigen}.  
\enproof

\subsection{Symmetry} 
%notation for set of wts? 
\nc{\Orb}{\on{Orb}} 

Let $\boldsymbol\lambda \in \Lambda_n$. 
There is a $(\Z/2\Z)^{\lfloor n/2 \rfloor}$-action on $\GT(\boldsymbol\lambda)$, given by flipping the signs of the entries  $m_{p,2p}$. Let $\Orb(\boldsymbol\lambda)$ be the set of orbits of this action. 
Note that if $\boldsymbol\lambda \in \Lambda^{\frac12 + \Z}_n$ then the action is free. On the other hand, if $\boldsymbol\lambda \in \Lambda^{\Z}_n$ then the action has non-trivial stabilizers whenever one of the entries $m_{p,2p}$ is zero. In particular, Gelfand--Tsetlin patterns with all first column entries equal to zero (possibly except the first row) are fixed by the action. 

We will say that $\sigma \in \GT(\boldsymbol\lambda)$ is \emph{standard} if all the first column entries are non-negative. Let $\GT^{\mathrm{st}}(\boldsymbol\lambda) \subset \GT(\boldsymbol\lambda)$ denote the subset of standard Gelfand--Tsetlin patterns. Clearly, there is a bijection $\Orb(\boldsymbol\lambda) \longleftrightarrow \GT^{\mathrm{st}}(\boldsymbol\lambda)$. 

\Cor \label{cor: symmetry cnst orbits} 
The eigenvalues of $\thvar{k}$ are constant on $(\Z/2\Z)^{\lfloor n/2 \rfloor}$-orbits. 
Hence 
\[
\ichmap\left(  L_{\boldsymbol\lambda}(a) \right)  = \sum_{[\sigma] \in \Orb(\boldsymbol\lambda)} |[\sigma]| \cdot \Xi(\sigma). 
\]
\encor

\Proof 
By Theorem \ref{thm: Theta-action on GT}, 
the sign flip $m_{p,2p} \mapsto - m_{p,2p}$ only affects the eigenvalues of $\thvar{2p-1}$ and $\thvar{2p}$. Noting that 
$l_{p,2p} = m_{p,2p}$, it is clear that the polynomial factors in the corresponding formulae depending on $l_{p,2p}$ are invariant under the sign flip. 
\enproof

In the non-classical case, independence of the boundary $q$-character of $\vec{\varepsilon}$, observed in Corollary \ref{cor: iq char formula non-classical}, is an analogue of the symmetry present in the classical case. 

\subsection{Tableaux interpretation} 
\nc{\SST}{\on{SST}^{\mathrm{ort}}} 
\nc{\Ups}{\boldsymbol{\Upsilon}} 
\label{subsec: tableaux int}

In subsections \S \ref{subsec: tableaux int}--\S \ref{subsec: partial order} we assume that $\boldsymbol\lambda \in \Lambda_n^{\Z}$. The modifications needed in the half-integral and non-classical cases will be given in \S \ref{subsec: adjust}.

Given $\boldsymbol\lambda \in \Lambda_n^{\Z}$, let $\SST_n(\boldsymbol\lambda)$ be the set of semistandard Young tableaux with the following properties:
\begin{enumerate}
\item of shape $(\lambda_1, \cdots, \lambda_{\lceil n/2 \rceil})$;
\item with entries from the alphabet $\{1, \cdots, n\}$;
\item $k$-th row contains entries from the restricted alphabet $\{2k-1, \cdots, n\}$. 
\end{enumerate}

\Lem \label{lem: Orb SST}
There is a canonical bijection
\[
\Orb(\boldsymbol\lambda) \ \longleftrightarrow \ \SST_n(\boldsymbol\lambda). 
\]
\enlem

\Proof
This bijection is constructed in the same way as the usual correspondence between  Gelfand--Tsetlin patterns for $\mathfrak{gl}_{n+1}$ and semistandard Young tableaux (see, e.g., \cite[\S 4.5]{FrMukHopf} or \cite[Remark 2.2]{MolevGT}). Namely, given a standard $\sigma \in \GT(\boldsymbol\lambda)$, one constructs the associated tableau in the following way. We view $\sigma$ as a sequence of partitions $\lambda^{(1)} \subseteq \cdots \subseteq \lambda^{(n)}$, where $\lambda^{(k)} = (m_{1,k+1}, \cdots m_{\lceil k/2 \rceil, k+1})$. Then the skew diagram $\lambda^{(k)}$ is a horizontal strip, and the associated tableau is obtained by placing the entry $k$ into each box of $\lambda^{(k)}/\lambda^{(k-1)}$. It is clear that the result is a semistandard Young tableau satisfying properties (1)--(2). Property (3) follows from the fact that, for each odd $k$, the partitions $\lambda^{(k)}$ and $\lambda^{(k+1)}$ are of the same size. It is easy to see that this procedure can be reversed, yielding a bijection. 
\enproof 

%half int 

We follow the English convention for displaying Young diagrams. We identify a box in a Young diagram with a pair $(i,j)\in \Z_{\geq 1}^2$, where $i$ is the row-coordinate, and $j$ the column coordinate of the box (counting from the top and from the left). The content of a box $(i,j)$ is $c(i,j) = j-i$. If a Young diagram is filled with numbers to form a tableau, let $t(i,j)$ be the number assigned to the box $(i,j)$. 

For any $b \in \C^\times$ and $1 \leq i \leq n$, let 
\[ \boldsymbol{\Upsilon}_{i,b} = \ourY_{i,bq^{i}}\mi \ourY_{i+1,bq^{i-1}}.\]
Here we formally set $\ourY_{n+1,bq^{n}} = 1$. To each tableau $T \in \SST_n(\boldsymbol\lambda)$, we associate a monomial
\eq \label{eq: Ups rule}
\Xi(T) = \prod_{(i,j) \in T} \Ups_{t(i,j), q^{2c(i,j)}}. 
\eneq
Moreover, let $\gamma_T$ be the number of rows of $T$ which satisfy the condition: 
\begin{itemize} 
\item the first entry of the row equals $2k-1$ (for the $k$-th row). 
\end{itemize}

\begin{exam} \label{exa: GT pattern vs tableau}
Consider the Gelfand--Tsetlin pattern below and the associated semistandard Young tableau. 
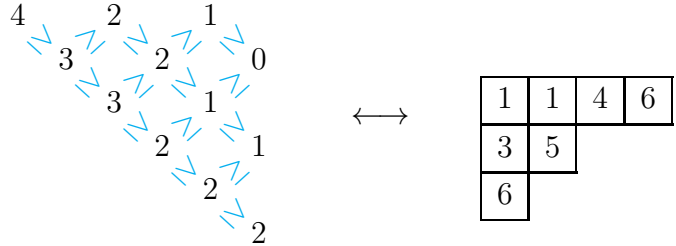
\[
\begin{tikzcd}[column sep=2pt,row sep=0pt]
4 \arrow[rd, phantom, sloped, "\geq", cyan]  & & 2 \arrow[rd, phantom, sloped, "\geq", cyan]  &  & 1 \arrow[rd, phantom, sloped, "\geq", cyan]  \\
& 3 \arrow[ru, phantom, sloped, "\geq", cyan] \arrow[rd, phantom, sloped, "\geq", cyan]  & & 2  \arrow[rd, phantom, sloped, "\geq", cyan] \arrow[ru, phantom, sloped, "\geq", cyan] & &   0  \\
 & & 3 \arrow[rd, phantom, sloped, "\geq", cyan] \arrow[ru, phantom, sloped, "\geq", cyan] & & 1 \arrow[rd, phantom, sloped, "\geq", cyan] \arrow[ru, phantom, sloped, "\geq", cyan] \\
 & & & 2 \arrow[rd, phantom, sloped, "\geq", cyan] \arrow[ru, phantom, sloped, "\geq", cyan]   & &  1  \\
 & & & & 2 \arrow[rd, phantom, sloped, "\geq", cyan] \arrow[ru, phantom, sloped, "\geq", cyan] \\
 & & & & &  2 
\end{tikzcd} 
\qquad \longleftrightarrow \qquad
\begin{ytableau}
1 & 1 & 4 & 6 \\
3 & 5 \\
6
\end{ytableau} 
\]
We have 
\[
\gamma_T = 2, \qquad 
\Xi(T) = \Ups_{1,1} \Ups_{1,q^{2}} \Ups_{3,q^{-2}} \Ups_{4,q^{4}} \Ups_{5, 1}, \Ups_{6, q^{-4}} \Ups_{6, q^{6}}. 
\]
In terms of the $\mathbf{Y}$-variables, one gets 
\[ \reqnomode
\Xi(T) = \ourY_{1,q}\mi \ourY_{1, q^3}\mi \ourY_{2, q^2} \ourY_{3,q}\mi \ourY_{4,q^8}\mi \ourY_{5,q^7} 
\ourY_{5,q^5}\mi \ourY_{6,q^4} \ourY_{6,q^2}\mi \ourY_{6,q^{12}}. \tag*{$\triangledown$}
\] 

\end{exam}

The following corollary allows us to read off the boundary $q$-character directly from tableaux combinatorics, in analogy to \cite[Lemma 4.7]{FrenMuk-comb}. 

\Cor \label{cor: tableaux formula for ich}
The boundary $q$-character of $L_{\boldsymbol\lambda}(a)$ is given~by
\[
\ichmap\left(  L_{\boldsymbol\lambda}(a) \right) = \sum_{T \in \SST(\boldsymbol\lambda)} 2^{\gamma_T} \cdot \Xi(T). 
\]
\encor

\Proof
This follows from Corollary \ref{cor: iq char formula classical} and Lemma \ref{lem: Orb SST}. 
\enproof

\subsection{Weight vectors}

There exists a unique Gelfand--Tsetlin pattern in $\GT(\boldsymbol\lambda)$ 
 whose entries are constant along the diagonals, i.e., it is of the form 
\[
\begin{tikzcd}[column sep=2pt,row sep=0pt]
\lambda_{1} \arrow[rd, phantom, sloped, "\geq", cyan]  & & \lambda_{2} \arrow[rd, phantom, sloped, "\geq", cyan]  &  & ... & &  \lambda_{r} \arrow[rd, phantom, sloped, "\geq", cyan]  \\
& \lambda_{1} \arrow[ru, phantom, sloped, "\geq", cyan] & & \lambda_{2}   & & & &  \lambda_{r}  \\
& & ... \\
& & & ... \\
& & & & \lambda_{1} \arrow[rd, phantom, sloped, "\geq", cyan]  & & \lambda_{2} \arrow[rd, phantom, sloped, "\geq", cyan]  \\
& & & & & \lambda_{1} \arrow[rd, phantom, sloped, "\geq", cyan] \arrow[ru, phantom, sloped, "\geq", cyan]   & &  \lambda_{2}  \\
& & & & & & \lambda_{1} \arrow[rd, phantom, sloped, "\geq", cyan] \arrow[ru, phantom, sloped, "\geq", cyan] \\
& & & & & & &  \lambda_{1}. 
\end{tikzcd}
\] 
We call it the \emph{highest weight} Gelfand--Tsetlin pattern and denote it by $\sigma_0$. 
This name is motivated by the fact that $\sigma_0$ is proportional to the highest weight vector of $L_{\boldsymbol\lambda}$ in the sense of Letzter and Wenzl \cite{Letzter-cartan, Wenzl}. However, it must be emphasized that, in general, Gelfand--Tsetlin pattern basis vectors are \emph{not} weight vectors, i.e., they are not eigenvectors for the action of Letzter's Cartan subalgebra $\mathcal{T} = \langle B_i \mid i \mbox{ is odd } \rangle$. This is to be expected since, in the classical limit, the Gelfand--Tsetlin basis is also not a weight basis as it arises from a chain of inclusions of alternating type $\mathsf{B}$ and type $\mathsf{D}$ Lie algebras.

More generally, consider Gelfand--Tsetlin patterns of the form 
\eq \label{eq: fixed GT pat}
\begin{tikzcd}[column sep=2pt,row sep=0pt]
\lambda_{1} \arrow[rd, phantom, sloped, "\geq", cyan]  & & \lambda_{2} \arrow[rd, phantom, sloped, "\geq", cyan]  &  & ... & &  \lambda_{r} \arrow[rd, phantom, sloped, "\geq", cyan]  \\
& \lambda_{1} \arrow[ru, phantom, sloped, "\geq", cyan] & & \lambda_{2}   & & & &  m_{r,2r}  \\
& & ... \\
& & & ... \\
& & & & \lambda_{1} \arrow[rd, phantom, sloped, "\geq", cyan]  & & \lambda_{2} \arrow[rd, phantom, sloped, "\geq", cyan]  \\
& & & & & \lambda_{1} \arrow[rd, phantom, sloped, "\geq", cyan] \arrow[ru, phantom, sloped, "\geq", cyan]   & &  m_{2,4}  \\
& & & & & & \lambda_{1} \arrow[rd, phantom, sloped, "\geq", cyan] \arrow[ru, phantom, sloped, "\geq", cyan] \\
& & & & & & &  m_{1,2}. 
\end{tikzcd}
\eneq
with $m_{p,2p}$ arbitrary. 
The following corollary exhibits a larger family of eigenvectors for Letzter's Cartan subalgebra $\mathcal{T}$. 

\Cor
Any Gelfand--Tsetlin pattern of the form \eqref{eq: fixed GT pat} is an eigenvector for the commutative algebra $\mathcal{T}$. 
\encor

\Proof
This follows from Theorem \ref{thm: GIK B-action on GT}. 
\enproof

\Cor \label{cor: hw gt mon}
The monomial corresponding to the highest weight Gelfand--Tsetlin pattern $\sigma_0$ is 
\eq 
\Xi(\sigma_0) = \prod_{i=1}^{\lceil n/2 \rceil} \prod_{j=0}^{\lambda_i-1} \Ups_{2i-1, q^{2j-i+1}} 
= \prod_{i=1}^{\lceil n/2 \rceil} \prod_{j=0}^{\lambda_i-1} \ourY_{2i-1, q^{2j+i}}\mi \ourY_{2i, q^{2j+i-1}}.  
\eneq
In particular, $\Xi(\sigma_0)$ is neither dominant nor anti-dominant. 
\encor

\Proof
The pattern $\sigma_0$ corresponds to a tableau of the form 
\[
\ytableausetup
{mathmode, boxframe=normal} 
\begin{ytableau}
1 & 1 & \none[\dots] & 1 & \none[\dots] & 1 & \none[\dots] & 1 \\
3 & 3 & \none[\dots] & 3 & \none[\dots] & 3 \\
\none[\vdots] & \none[\vdots] & \none & \none[\vdots] & \none[\iddots] \\
 n & n & \none[\dots] & n
\end{ytableau}
\]
(if $n$ is odd). We can read the corresponding monomial from the tableau using \eqref{eq: Ups rule}. 
\enproof

\subsection{Minimal GT patterns} 

Let $\boldsymbol\lambda \in \Lambda^{\Z}_n$. 
There exists a unique Gelfand--Tsetlin pattern in $\GT(\boldsymbol\lambda)$ such that, whenever any entry is decreased (except the top row), the Gelfand--Tsetlin constraints are violated. 
We call it the \emph{minimal} Gelfand--Tsetlin pattern and denote it by $\sigma_{\mathrm{min}}$. It is of the form 
\eq \label{eq: min GT pat}
\begin{tikzcd}[column sep=2pt,row sep=0pt]
\lambda_{1} \arrow[rd, phantom, sloped, "\geq", cyan]  & & \lambda_{2} \arrow[rd, phantom, sloped, "\geq", cyan]  &  & \lambda_3 & & ... & & \lambda_{r}\\ %\arrow[rd, phantom, sloped, "\geq", cyan]  \\
& \lambda_{2} \arrow[ru, phantom, sloped, "\geq", cyan] \arrow[rd, phantom, sloped, "\geq", cyan] & & \lambda_{3} \arrow[ru, phantom, sloped, "\geq", cyan]   & & ... & & \lvert \lambda_r \rvert & & 0  \\
& & \lambda_3 \arrow[ru, phantom, sloped, "\geq", cyan]  & & ... & & \lvert \lambda_r \rvert & & 0 \\
& & & ... & & \lvert \lambda_r \rvert & & 0 & & ... \\
& & & & \lvert \lambda_{r} \rvert  & & 0  & & ...   \\
& & & & & 0   & &  ... & & 0  \\
& & & & & & ... & & 0 \\
& & & & & & &  0 & & 0 \\ 
& & & & & & & &  0 &  \\ 
& & & & & & & & & 0.  
\end{tikzcd}
\eneq 
One may think of $\sigma_{\mathrm{min}}$ as the lowest weight type $\mathsf{A}$ Gelfand--Tsetlin pattern associated to $\boldsymbol\lambda$, extended be zero to a Gelfand--Tsetlin pattern of orthogonal type.

\Cor \label{cor: min gt mon}
The minimal Gelfand--Tsetlin pattern $\sigma_{\mathrm{min}}$ is, in general, \emph{not} a weight vector. 
The corresponding monomial is 
\[
\Xi(\sigma_{\mathrm{min}}) = \prod_{i=1}^{\lceil n/2 \rceil} \prod_{j=0}^{l_i-1} \ourY_{\lfloor n/2 \rfloor + i, q^{\lfloor n/2 \rfloor + i + 2(m_{i-1}+j)}}\mi, 
\] 
where $l_1 = \lambda_{\lceil n/2 \rceil}$, $l_2 = \lambda_{\lceil n/2 \rceil - 1} - \lambda_{\lceil n/2 \rceil}$, $\cdots$, $l_{\lceil n/2 \rceil - 1} = \lambda_1 - \lambda_2$; and $m_ i = \sum_{j=1}^i l_i$, with $m_0 = 0$. 
In particular, $\Xi(\sigma_0)$ is anti-dominant. 
\encor 

\Proof
Let $z=\lceil n/2 \rceil $. 
The pattern $\sigma_{\mathrm{min}}$ corresponds to a tableau of the form 
\[
\ytableausetup
{mathmode, boxframe=normal} 
\begin{ytableau}
 1 & \none[\dots] & 1 & 2 & \none[\dots] & 2 & \none[\dots] & \none[\dots]  & \scriptstyle z{-}1 &  \none[\dots] & \scriptstyle z{-}1 &  z &  \none[\dots] &  z \\
2 & \none[\dots] & 2 & 3 & \none[\dots] & 3 & \none[\dots] & \none[\dots]  & z &  \none[\dots] &  z  \\
 \none[\vdots] & \none & \none[\vdots] & \none[\iddots] \\
 z & \none[\dots] & z 
\end{ytableau}
\]
with all entries shifted by $\lfloor n/2 \rfloor$. Observe that all the positive powers of $\ourY$ cancel out in
the product of $\Ups$'s contributed by each column. In effect, $\Xi(\sigma_{\mathrm{min}})$ is determined only by the first row of the tableau, i.e., it is the product of the negative powers of $\ourY$ occurring in the $\Ups$'s contributed by the first row. 
\enproof 

Corollaries \ref{cor: hw gt mon} and \ref{cor: min gt mon} indicate an interesting trade-off between (anti-)dominance and the property of being a weight vector.

\subsection{Drinfeld polynomials} 

The Drinfeld polynomials of evaluation representations for quantum affine algebras were first computed and described in terms of $q$-strings in \cite[Theorem 3.5]{CP-small}. Below we propose a notion of a Drinfeld polynomial for evaluation representations in the split quantum symmetric pair setting.

Recall that, given a polynomial $P(z) \in \C[z]$ with constant term $1$, $P^\dag(z)$ denotes the polynomial with constant term~$1$ whose roots are obtained from those of $P(z)$ via the transformation $x \mapsto a^{-2} x\mi$. 
Corollary \ref{cor: hw gt mon} implies that the eigenvalue of $\thvar{k}$ on $\sigma_{0}$ is of the form  
\begin{align*}
\frac{P_{2i-1}(q z)}{P_{2i-1}(q\mi z)} \frac{P^\dag_{2i-1}(q\mi z)}{P^\dag_{2i-1}(q z)},& \ &
P_{2i-1}(z) &= \prod_{j=0}^{\lambda_i-1} (1-aq^{2j+i}z), \\ 
\frac{P_{2i}(q\mi z)}{P_{2i}(q z)} \frac{P^\dag_{2i}(q z)}{P^\dag_{2i}(q\mi z)},& \ &
P_{2i}(z) &= \prod_{j=0}^{\lambda_i-1} (1-aq^{2j+i-1}z), 
\end{align*} 
with the convention that the linear factors in $P_k(z)$ contain only non-negative powers of $q$. 
It is clear that the polynomials with odd and even indices contain the same information. 
We call $(P_{2i-1}(z))_{i=1}^{\lceil n/2 \rceil}$ the \emph{Drinfeld polynomials} associated to $L_{\boldsymbol\lambda}(a)$. 

Moreover, Corollary \ref{cor: min gt mon} implies that the eigenvalue of $\thvar{i}$ on $\sigma_{\mathrm{min}}$ equals $1$ for $1 \leq i \leq \lfloor n/2 \rfloor$. Otherwise, it is of the form  
\[
\frac{Q_i(q z)}{Q_i(q\mi z)}\frac{Q^\dag_i(q\mi z)}{Q^\dag_i(q z)}, \qquad 
Q_i(z) = \prod_{j=0}^{l_i-1} (1-aq^{\lfloor n/2 \rfloor + i + 2(m_{i-1}+j)}z), 
\] 
where $l_1 = \lambda_{\lceil n/2 \rceil}$, $l_2 = \lambda_{\lceil n/2 \rceil - 1} - \lambda_{\lceil n/2 \rceil}$, $\cdots$, $l_{\lceil n/2 \rceil - 1} = \lambda_1 - \lambda_2$; and $m_ i = \sum_{j=1}^i l_i$, with $m_0 = 0$. 
%same info 
We call $(Q_k(z))_{k=\lceil n/2 \rceil}^{n}$ the \emph{dual Drinfeld polynomials} associated to $L_{\boldsymbol\lambda}(a)$.

\subsection{Partial order}
\label{subsec: partial order}

One can define a partial order on $\GT(\boldsymbol\lambda)$ by declaring that $\sigma \leq \sigma'$ whenever each entry of $\sigma'$ is greater than or equal to the corresponding entry of $\sigma$. Analogously, one can define a partial order on $\SST(\boldsymbol\lambda)$ by declaring that $T < T'$  whenever each box of $T'$ contains a number greater  than or equal to the number in the corresponding box of $T$. It is easy to see that the bijection from Lemma \ref{lem: Orb SST} is order-reversing. Let us denote the tableau corresponding to $\sigma_{\mathrm{min}}$ by $T_{\mathrm{max}}$. Given another tableau $T \in \SST(\boldsymbol\lambda)$, we call $E(T) = T_{\mathrm{max}} - T$, with boxes containing $0$ truncated, the associated $\emph{excess diagram}$. Let $e(i,j)$ be the number assigned to the box $(i,j)$ in $E(T)$. 

\begin{exam} 
Take $n=6$ and $\boldsymbol\lambda = ([4],[2],[1])$. 
Then $\sigma_{\mathrm{min}}$, together with the associated tableau, is: 
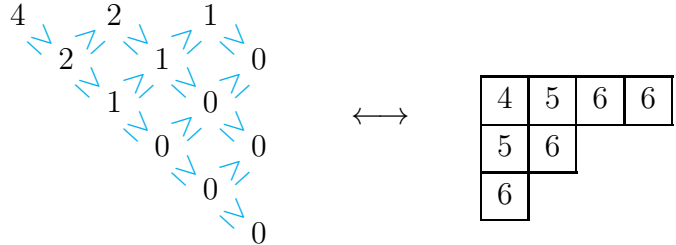
\[
\begin{tikzcd}[column sep=2pt,row sep=0pt]
4 \arrow[rd, phantom, sloped, "\geq", cyan]  & & 2 \arrow[rd, phantom, sloped, "\geq", cyan]  &  & 1 \arrow[rd, phantom, sloped, "\geq", cyan]  \\
& 2 \arrow[ru, phantom, sloped, "\geq", cyan] \arrow[rd, phantom, sloped, "\geq", cyan]  & & 1  \arrow[rd, phantom, sloped, "\geq", cyan] \arrow[ru, phantom, sloped, "\geq", cyan] & &   0  \\
 & & 1 \arrow[rd, phantom, sloped, "\geq", cyan] \arrow[ru, phantom, sloped, "\geq", cyan] & & 0 \arrow[rd, phantom, sloped, "\geq", cyan] \arrow[ru, phantom, sloped, "\geq", cyan] \\
 & & & 0 \arrow[rd, phantom, sloped, "\geq", cyan] \arrow[ru, phantom, sloped, "\geq", cyan]   & &  0  \\
 & & & & 0 \arrow[rd, phantom, sloped, "\geq", cyan] \arrow[ru, phantom, sloped, "\geq", cyan] \\
 & & & & &  0 
\end{tikzcd} 
\qquad \longleftrightarrow \qquad
\begin{ytableau}
4 & 5 & 6 & 6 \\
5 & 6 \\
6
\end{ytableau} 
\]
Let $\sigma$ be the Gelfand--Tsetlin pattern, and $T$ the corresponding tableau, from Example~\ref{exa: GT pattern vs tableau}. Then 
\[
\begin{ytableau} 
3 & 4 & 2  \\
2 & 1 
\end{ytableau} 
\] 
is the associated excess diagram. 
\rmkend
\end{exam}

Recall that in Definition \ref{defi: monomial ordering}, we defined a partial order on monomials. The following corollary shows that the monomial corresponding to the minimal Gelfand--Tsetlin pattern is the least element in this order. This is a (dual) analogue of a well-known property of $q$-characters of quantum affine algebras (see \cite[Theorem 4.1]{FrenMuk-comb}). 

\Cor 
\label{cor: root property}
Let $\boldsymbol\lambda$ be integral. Then 
\[
\ichmap\left(  L_{\boldsymbol\lambda}(a) \right) = \Xi(\sigma_{\mathrm{min}}) \left(1 + \sum_{T_{\mathrm{max}} \neq T \in \SST(\boldsymbol\lambda)} 2^{\gamma_T} \cdot \prod_{(i,j) \in E(T)} \prod_{k=0}^{e(i,j)-1} \ourA_{t(i,j) - k, q^{t(i,j) + 2c(i,j) - k -1}}  \right). 
\]
\encor

\Proof
This follows immediately from Corollary \ref{cor: tableaux formula for ich} and the fact that decreasing the $(i,j)$-th entry in a tableau by $1$ changes the corresponding monomial by 
\begin{align*}
\Ups_{t(i,j), q^{2c(i,j)}}\mi \Ups_{t(i,j)-1, q^{2c(i,j)}} &= 
\ourY_{t(i,j), q^{2c(i,j)+t(i,j)}} \ourY_{t(i,j)+1,q^{2c(i,j)+t(i,j)-1}}\mi  \\ 
&\ \ \cdot \ourY_{t(i,j)-1, q^{2c(i,j)+t(i,j)-1}}\mi \ourY_{t(i,j),q^{2c(i,j)+t(i,j)-2}} \\
&= \ourA_{t(i,j), q^{2c(i,j)+t(i,j)-1}}. \qedhere
\end{align*} 
\enproof

\subsection{Adjustments for the half-integral and non-classical cases} 
\label{subsec: adjust}

We first consider the classical half-integral case. 
Given $\boldsymbol\lambda \in \Lambda^{\frac12 + \Z}_n$, let $\widetilde{\boldsymbol\lambda} = ([\lambda_1 - \frac12], \cdots, [\lambda_{\lceil n/2 \rceil} - \frac12])$ be the integral weight arising through a uniform shift by $- \frac12$. 

\Cor
Let $\boldsymbol\lambda \in \Lambda^{\frac12 + \Z}_n$. 
\begin{enumerate}
\item There is a canonical bijection\footnote{If $\widetilde{\boldsymbol\lambda}$ is the zero partition, then $\SST_n(\widetilde{\boldsymbol\lambda})$ consists of the empty tableau.} $\Orb(\boldsymbol\lambda) \longleftrightarrow \SST_n(\widetilde{\boldsymbol\lambda})$. 
\item The boundary $q$-character of $L_{\boldsymbol\lambda}(a)$ is given~by
\[
\ichmap\left(  L_{\boldsymbol\lambda}(a) \right) = \sum_{T \in \SST(\widetilde{\boldsymbol\lambda})} 2^{\lfloor n/2 \rfloor} \cdot {\Xi}(T), \qquad 
{\Xi}(T) = \prod_{k=1}^n \ourX_{k,q^{\frac12}}\mi \prod_{(i,j) \in T}  \Ups_{t(i,j), q^{2c(i,j)+1}}. 
\]
\item $\ichmap\left(  L_{\boldsymbol\lambda}(a) \right)$ can be expressed as
\[
2^{\lfloor n/2 \rfloor} \Xi(\sigma_{\mathrm{min}}) \left(1 + \sum_{T_{\mathrm{max}} \neq T \in \SST(\widetilde{\boldsymbol\lambda})} \prod_{(i,j) \in E(T)} \prod_{k=0}^{e(i,j)-1} \ourA_{t(i,j) - k, q^{t(i,j) + 2c(i,j) - k}}  \right), 
\]
where 
\[
\Xi(\sigma_{\mathrm{min}}) = \prod_{k=1}^n \ourX_{k,q^{\frac12}}\mi  \prod_{i=1}^{\lceil n/2 \rceil} \prod_{j=0}^{l_i-1} \ourY_{\lfloor n/2 \rfloor + i, q^{\lfloor n/2 \rfloor + i + 2(m_{i-1}+j)+1}}\mi
\] 
and $l_1 = \lambda_{\lceil n/2 \rceil} - \frac12$, $l_2 = \lambda_{\lceil n/2 \rceil - 1} - \lambda_{\lceil n/2 \rceil}$, $\cdots$, $l_{\lceil n/2 \rceil - 1} = \lambda_1 - \lambda_2$; and $m_ i = \sum_{j=1}^i l_i$, with $m_0 = 0$. 
\end{enumerate}
\encor 

Next, consider the non-classical case. We define the shifted weight $\widetilde{\boldsymbol\lambda}_+$ in analogy to the half-integral case. 

\Cor
Let $\boldsymbol\lambda_+ \in \Lambda^+_n$. 
\begin{enumerate} 
\item The boundary $q$-character of $L_{\boldsymbol\lambda_+, {\vec{\varepsilon}}}(a)$ is given~by
\[
\ichmap\left(  L_{\boldsymbol\lambda_+, {\vec{\varepsilon}}}(a) \right) = \sum_{T \in \SST(\widetilde{\boldsymbol\lambda})} {\Xi}(T), \qquad 
{\Xi}(T) = \prod_{k=1}^n \widetilde{\ourX}_{k,q^{\frac12}}\mi \prod_{(i,j) \in T}  \Ups_{t(i,j), -q^{2c(i,j)+1}}. 
\]
\item $\ichmap\left( L_{\boldsymbol\lambda_+, {\vec{\varepsilon}}}(a) \right)$ can be expressed as 
\[
\Xi(\sigma_{\mathrm{min}}) \left(1 + \sum_{T_{\mathrm{max}} \neq T \in \SST(\widetilde{\boldsymbol\lambda})} \prod_{(i,j) \in E(T)} \prod_{k=0}^{e(i,j)-1} \ourA_{t(i,j) - k, -q^{t(i,j) + 2c(i,j) - k}}  \right), 
\]
where 
\[
\Xi(\sigma_{\mathrm{min}}) = \prod_{k=1}^n \widetilde{\ourX}_{k,q^{\frac12}}\mi \prod_{i=1}^{\lceil n/2 \rceil} \prod_{j=0}^{l_i-1} \ourY_{\lfloor n/2 \rfloor + i, -q^{\lfloor n/2 \rfloor + i + 2(m_{i-1}+j)+1}}\mi. 
\] 
\end{enumerate}
\encor

\subsection{Fundamental representations}

Below we calculate the boundary $q$-character of the evaluation $W_{\pi_i}(a)$ of the $i$-th fundamental representation. 

\Cor 
If $n$ is even and $1 \leq i \leq n/2 -1$, or $n$ is odd and $1 \leq i \leq (n+1)/2 - 2$, then 
\[
\ichmap(W_{\pi_i}(a)) = \sum_{\substack{1 \leq j_1 < j_2 \leq \cdots \leq j_i \leq n, \\ j_k \geq 2k-1 }} 2^{\sum \delta_{j_k = 2k-1}} \Ups_{j_1, 1} \cdots \Ups_{j_i, q^{2-2i}}. 
\]
Moreover, 
\[
\ichmap(W_{\pi_{\lceil n/2 \rceil}}(a)) = 2^{\lfloor n/2 \rfloor} \prod_{k=1}^n \ourX_{k,q^{\frac12}}\mi 
= \ichmap(W_{\pi_{\lceil n/2 \rceil}-1}(a)),
\]
where the latter equality holds if $n$ is odd. 
\encor

\Proof 
In the first case, to compute this $q$-character, we need to sum over tableaux with only one column of length $i$. 
The formula now follows from Corollary \ref{cor: tableaux formula for ich}. The other two cases are immediate from Corollary \ref{cor: iq char formula classical}.  
\enproof

\subsection{Example: rank $2$} 

Let us consider the case $n=2$ in more detail. Let $\boldsymbol\lambda = ([\lambda]) \in \Lambda_2$. The classical simple modules can be constructed explicitly as quotients of the Verma modules from \S \ref{subsec: Verma cnstr}. Letting $v_0$ be the image of $1$ in the quotient $\ui / I_{\boldsymbol\lambda}$, the set $\{ v_r := B_2^r \cdot v_0 \mid 0 \leq r \leq 2\mu \}$ is a basis of $L_{\boldsymbol\lambda}$. In particular, $\dim L_{\boldsymbol\lambda} = 2\lambda+1$. It is easy to see that each $v_r$ is a multiple of the Gelfand--Tsetlin pattern vector 
\[
\begin{tikzcd}[column sep=2pt,row sep=0pt]
\lambda \arrow[rd, phantom, sloped, "\geq", cyan] \\
& \lambda - r. 
\end{tikzcd}
\]

\Cor \label{thm: vi eigenvector for theta}
The series $\thvar{k}$ act diagonally on $v_r \in L_{\boldsymbol\lambda}(a)$ (for $k \in \{1,2\}$). The corresponding series of eigenvalues $\theta_{k,r}(z)$ is given by
\[
\theta_{k,r}(z) = 
\left\{ \begin{array}{l l}
\displaystyle \frac{(1 - q^{2\lambda-2r}az)(1 - q^{2r-2\lambda}az)}{(1- az)^2} & \mbox{if } k=1, \\[2pt]
\displaystyle \frac{(1-qaz)(1-q\mi az)(1-q^{2\lambda+1}az)(1-q^{-2\lambda-1}az)}{(1-q^{2\lambda-2r-1}az)(1-q^{1+2r-2\lambda}az)(1-q^{2\lambda+1-2r}az)(1-q^{2r-1-2\lambda}az)} & \mbox{if } k=2. \\[2pt] 
\end{array} \right. 
\]
\encor

\Cor
The boundary $q$-character of the classical evaluation module $L_{\boldsymbol\lambda}(a)$ has the following properties.  
\begin{enumerate}
\item If $\lambda \in \Z_{\geq 0}$, then the vectors $v_r$ and $v_{2\lambda - r}$ ($0 \leq r \leq \lambda-1$) lie in the same $\mathcal{H}$-eigenspace. The corresponding monomial is  
\[
\prod_{l=0}^{\lambda-1-r} \ourY_{1,q^{1+2l}}\mi \prod_{l=1}^{\lambda-1-r} \ourY_{2, q^{2l}} \prod_{l=0}^{r-1} \ourY_{2, q^{2\lambda - 2l}}\mi.   
\]
The monomial corresponding to the eigenvector $v_\lambda$ is  
$\prod_{l=1}^{\lambda} \ourY_{2, q^{2l}}\mi$. 
It is the lowest anti-dominant monomial in the ordering from Definition \ref{defi: monomial ordering}.  
\item If $\lambda \in \frac{1}{2} + \Z_{\geq 0}$, then the vectors $v_r$ and $v_{2\lambda - r}$ ($0 \leq r \leq \lambda-\frac{1}{2}$) lie in the same $\mathcal{H}$-eigenspace. The corresponding monomial is  
\[
\ourX_{1,q^{2\lambda-2r-\frac{1}{2}}}\mi \ourX_{2, q^{2\lambda - 2r - 1\frac{1}{2}}}   \prod_{l=0}^{\lambda-1\frac{1}{2}-r} \ourY_{1,q^{1+2l}}\mi \prod_{l=1}^{\lambda-1\frac{1}{2}-r} \ourY_{2, q^{2l}} \prod_{l=0}^{r-1} \ourY_{2, q^{2\lambda - 2l}}\mi.   
\] 
\item For any $ \lambda \in \frac{1}{2} \Z_{\geq 0}$, the boundary $q$-character of $L_{\boldsymbol\lambda}(a)$ is compatible with the Weyl group (i.e., $S_2$) action on the weights (with respect to $B_1$), in the sense that the monomial associated to $v_r$ (which is of weight $[\lambda-r]$) is the same as the one associated to $v_{2\lambda - r}$ (which is of weight $-[\lambda-r]$).  
\item The boundary $q$-character of $L_{\boldsymbol\lambda}(a)$ is given by 
\[
\chi_q^\imath \big(L_{\boldsymbol\lambda}(a) \big) = 
\left\{ \begin{array}{l l}
\displaystyle \prod_{l=1}^{\lambda} \ourY_{2, q^{2l}}\mi \left(1 + 2 \sum_{r=1}^\lambda \prod_{m=1}^r \ourA_{2, q^{2m-1}} \right) &    \mbox{if } \lambda \in \Z_{\geq 0}, \\[2pt]
\displaystyle 2 \cdot \ourX_{1,q^{\frac{1}{2}}}\mi \ourX_{2,q^{\frac{1}{2}}}\mi  \prod_{l=1}^{\lambda-\frac{1}{2}} \ourY_{2, q^{2l+1}}\mi\left( 1 + \sum_{r=1}^{\lambda-\frac{1}{2}} \prod_{m=1}^r \ourA_{2, q^{2m}} \right)  & \mbox{if } \lambda \in \frac{1}{2} + \Z_{\geq 0}. \\[2pt] 
\end{array} \right. 
\]
\item Each $v_r$ is also an eigenvector for $\mathcal{T} = \langle B_1 \rangle$. 
\end{enumerate}
\encor

Next, let $\boldsymbol\lambda_+ = ([\lambda]_+) \in \Lambda_2^+$. There are four pairwise non-isomorphic simple modules associated to this weight, each of dimension $\lambda + \frac12$. In each case, the boundary $q$-character is given by 
\[
\ichmap\left( L_{\boldsymbol\lambda_+, {\vec{\varepsilon}}}(a) \right) = \widetilde{\ourX}_{1,q^{\frac{1}{2}}}\mi \widetilde{\ourX}_{2,q^{\frac{1}{2}}}\mi  \prod_{l=1}^{\lambda-\frac{1}{2}} \ourY_{2, q^{2l+1}}\mi\left( 1 + \sum_{r=1}^{\lambda-\frac{1}{2}} \prod_{m=1}^r \ourA_{2, q^{2m}} \right). 
\]

%%%%%%%
%%%%%%%
%%%%%%%

\providecommand{\bysame}{\leavevmode\hbox to3em{\hrulefill}\thinspace}
\providecommand{\MR}{\relax\ifhmode\unskip\space\fi MR }
\providecommand{\MRhref}[2]{%
  \href{http://www.ams.org/mathscinet-getitem?mr=#1}{#2}
}
\providecommand{\href}[2]{#2}

\end{document}